\documentclass[11pt,fleqn]{article}
\usepackage[bookmarks=false,pdfpagemode=UseOutlines,plainpages=false,hypertexnames=false,pdfpagelabels,hyperindex=true,colorlinks=true,citecolor=blue]{hyperref}
\RequirePackage[OT1]{fontenc}
\RequirePackage{amsthm,amsmath}
\RequirePackage[round]{natbib}

\usepackage{amsmath,amsfonts,amssymb,euscript,mathrsfs}
\usepackage{slashbox}
\usepackage{float}
\usepackage[final]{pdfpages}
\usepackage{multirow}

\usepackage{txfonts}
\usepackage[T1]{fontenc} 
\usepackage{fullpage}
\usepackage{graphicx}
\usepackage{dsfont}
\makeatletter
\def\captionof#1#2{{\def\@captype{#1}#2}}
\makeatother

\def\som{\displaystyle\sum}
\def\R{\mathbb{R}}

\def\N{\mathbb{N}}
\def\P{\mathbb{P}}
\def\E{\mathbb{E}}
\def\L{\mathbb{L}}

\def\R{\mathbb{R}}

\def\H{\mathbb{H}}
\def\[{\llbracket}
\def\]{\rrbracket}
\def\Tr{\mbox{Tr\,}}
\def\V{\mbox{Var\,}}

\def\o{o_{\P}\left(\frac{1}{\sqrt{T}}\right)}

\numberwithin{equation}{section}
\theoremstyle{plain}
\newtheorem{theo}{Theorem}
\newtheorem{lem}{Lemma}
\newtheorem{Prop}{Proposition}
\newtheorem{Def}{Definition}
\newtheorem{cor}{Corollary}

\begin{document}
\author{Lionel Truquet \footnote{UMR 6625 CNRS IRMAR, University of Rennes 1, Campus de Beaulieu, F-35042 Rennes Cedex, France and}
\footnote{ENSAI, Campus de Ker-Lann, rue Blaise Pascal, BP 37203, 35172 Bruz cedex, France. {\it Email: lionel.truquet@ensai.fr}.}
 }
\title{Parameter stability and semiparametric inference in time-varying ARCH models}
\date{}
\maketitle

\begin{abstract}
In this paper, we develop a complete methodology for detecting time-varying/non time-varying parameters in ARCH processes. For this purpose, we estimate and test various semiparametric versions of the time-varying ARCH model (tv-ARCH) which include two well known non stationary ARCH type models introduced in the econometric literature.    
Using kernel estimation, we show that non time-varying parameters can be estimated at the usual parametric rate of convergence and for a Gaussian noise, we construct estimates that are asymptotically efficient in a semiparametric sense. Then we introduce two statistical tests which can be used for detecting non time-varying parameters or for testing the second order dynamic. An information criterion for selecting the number of lags is also provided.  We illustrate our methodology with several real data sets.\end{abstract}



\section{Introduction}
The modeling of financial data using nonstationary time series has recently received considerable attention both in econometrics and in statistics. For classical daily series such as stock market indices or currency exchange rates, the stationarity assumption seems often incompatible with a long history of data and the necessity of using non stationary ARCH models has been pointed out by several authors. See for instance \citet{MS}, \citet{Starica}, \citet{ER}, \citet{Fryz} and the references therein.
However, it is difficult to find in the literature a consensus for representing non-stationary ARCH models. 
A natural approach is to allow time-varying parameters in the classical ARCH model of \citet{Engle}.
Such an extension has been proposed by \citet{DR} with the so-called time-varying ARCH model (tv-ARCH). The tv-ARCH processes are  defined by the recursive equations 
\begin{equation}\label{tvARCH}
X_t=\xi_t\sigma_t,\quad \sigma^2_t=a_0\left(\frac{t}{T}\right)+\sum_{j=1}^p a_j\left(\frac{t}{T}\right)X_{t-j}^2,\quad p+1\leq t\leq T,
\end{equation}
where for $0\leq j\leq p$, $a_j$ is a smooth function and $\xi$ a strong white noise with variance $1$. Since they can be locally approximated by stationary ARCH processes, the tv-ARCH processes are called locally stationary (the notion of local stationarity is introduced in \citet{Dahlhaus} for linear processes but the meaning of local stationarity for the non linear tv-ARCH can be found in \citet{DR}). From this important feature, a nice asymptotic theory can be developed for estimation of parameters, in particular local inference methods such as the local Quasi-maximum likelihood estimation studied in \citet{DR}, the local weighted least-squares estimation developed in \citet{Fryz} or the recursive online algorithms considered by \citet{DaSu}. 
In \citet{Fryz}, it was shown that tv-ARCH processes provide good fits and accurate forecasts for some financial series. 
However, statistical inference in model (\ref{tvARCH}) is complex, even for large samples, because $p+1$ functions have to be estimated using nonparametric methods. Thus, in practice, reducing complexity can be interesting to improve model fit or forecasts accuracy. For example, \citet{Starica} have shown that the simple model
\begin{equation}\label{Star}
X_t=\sigma\left(\frac{t}{T}\right)\xi_t,\quad 1\leq t\leq T, 
\end{equation}
with a smooth deterministic function $\sigma:[0,1]\rightarrow \R_+$ can already produce significantly better forecasts for the returns of the SP$\&$500 index than the classical GARCH$(1,1)$ model. A process of type (\ref{Star}) can be seen as a tv-ARCH process with zero lag coefficients. Note that model (\ref{Star}) does not assume autocorrelation for the absolute values or the squares of the process (a correlation property often called second order dynamic in the literature) but 
only some changes in the unconditional variance. 
In \citet{Starica}, it is argued that most of dynamic of the S$\&$P index can be explained
with a time-varying unconditional variance.  

But nonstationarity and second order correlation can also be combined in a very simple way, assuming constant lag coefficients in (\ref{tvARCH}):
\begin{equation}\label{lagconst}
X_t= \xi_t\sqrt{a_0\left(\frac{t}{T}\right)+\sum_{j=1}^pa_jX_{t-j}^2}.
\end{equation}
Model (\ref{lagconst}) combine a time-varying unconditional variance compatible with the analysis of \citet{Starica} and a second order dynamic for the series with a single nonparametric component. 
Note also that a process $(X_t)_t$ defined by equations (\ref{lagconst}) can be written using the multiplicative form 
\begin{equation}\label{lagconst+}
X_t=\sqrt{a_0\left(\frac{t}{T}\right)}\cdot Y_t,
\end{equation}
where 
$$Y_t=\xi_t\sqrt{1+\sum_{j=1}^pa_j\frac{a_0\left(\frac{t-j}{T}\right)}{a_0\left(\frac{t}{T}\right)}Y_{t-j}^2}\approx \xi_t \sqrt{1+\sum_{j=1}^pa_jY_{t-j}^2}$$
if we neglect the ratio $a_0\left(\frac{t-j}{T}\right)/a_0\left(\frac{t}{T}\right)$ which is of order $1+1/T$ when the function $a_0$ is positive and Lipschitz continuous over $[0,1]$. One can also notice that writing the model with the latter approximation or not lead to two processes that are both approximated by the same stationary ARCH processes with parameters $a_0(u),a_1,\ldots, a_p$ (see Lemma $1$ for this kind of approximation). In the stationary case, we remind that multiplying an ARCH process by a positive constant is equivalent to multiply the initial intercept coefficient. 
Then for a large sample size $T$, the process $(Y_t)_t$ behaves as a stationary ARCH process and $a_0(\cdot)$ is (up to a constant) the time-varying unconditional variance of the process $(X_t)_t$.  Such a multiplicative form for ARCH models has been first considered by \citet{ER} with the so-called Spline-GARCH model which
writes as model (\ref{lagconst+}) but with a GARCH$(1,1)$ process $(Y_t)_t$. 

Since the previous models satisfy the inclusions (\ref{Star})$\subset$ (\ref{lagconst})$\subset$ (\ref{tvARCH}), 
a natural question for any real data set is to test some properties of the lag coefficients.    
Testing the constancy of the lag coefficients can help to decide between model (\ref{lagconst}) and model (\ref{tvARCH}) while testing
the second order dynamic in model (\ref{lagconst}) is useful to determine if model (\ref{Star}) provides a sufficient fit.
Statistical tools to help the practitioners to choose among the three important specifications described above seems not available in the literature except in a recent paper of \cite{Pat} which introduces a test for the second order dynamic in model (\ref{lagconst}). 

In this paper, we propose a general approach for estimating an arbitrary subset of non time-varying coefficients in tv-ARCH processes.
Our estimators are $\sqrt{T}-$consistent and we will also study the semi-parametric asymptotic efficiency of our method when the noise is Gaussian. Using these results, we construct two statistical tests. The first test can be used to decide whether a given subset of parameters is time-varying or not. This test is based on a $\L^2$ distance between a nonparametric kernel estimator of the coefficients and the semiparametric estimator introduced in this paper. The second test can be used for deciding if the constant parameters are different from zero.  
Various applications can be considered as a simple particular case of our methodology: testing model (\ref{lagconst}) versus model (\ref{tvARCH}), testing model 
(\ref{Star}) versus model (\ref{lagconst}), estimating parameters and selecting lag variables in models (\ref{lagconst}) or (\ref{tvARCH}). When some coefficients are assumed to be non time-varying in (\ref{tvARCH}), the decomposition $X_t^2=\sigma_t^2+\left(\xi_t^2-1\right)\sigma_t^2$ leads to semiparametric inference in a time-varying regression model.
Detecting and estimating a parametric component in general time-varying regression models has been considered recently by \citet{ZW}.
However these authors do not consider the case of tv-ARCH processes with optimal moment condition for the marginal distribution, asymptotic semi-parametric efficiency for the estimation and Lipschitz continuity for the time-varying coefficients. Moreover, our approach for estimating non-time varying coefficients is quite different. 

We applied our methodology to three real data sets: the daily exchange rates between the US Dollar and the Euro or between the US Dollar and the Indian Rupee and the FTSE index. For the three series of interest, a non time-varying intercept is clearly rejected over the considered period. The conclusion for the lag coefficients depends on the series. In fitting model (\ref{lagconst}), we also found that incorporating non stationarity reduces the values of lag parameters with respect to the stationary case. Then the time-varying unconditional variance has an important contribution to volatility.   

The paper is organized as follows. In Section \ref{S2}, we introduce our notations and we describe the basis of our method for statistical inference in a tv-ARCH model for which some coefficients are assumed to be non time-varying. In Section \ref{S3}, we give the asymptotic results for our estimators and we discuss the problem of semiparametric asymptotic efficiency using the LAN theory. Statistical testing and their practical implementation are considered in Section \ref{S4} and Section \ref{S5} is devoted to real data applications. All the proofs of our results are postponed to the supplementary material which also contains many simulation studies showing the good behavior of our methodology. The Matlab codes and the data sets discussed in Section \ref{S5} are available at the URL
\begin{center}
https://github.com/time-varying/tests-and-estimation-for-tv-ARCH-           
\end{center}

\section{Semiparametric volatility and tv-ARCH processes}\label{S2}

\subsection{Formulation and notations}

In this section, we consider semiparametric versions of model (\ref{tvARCH}), assuming that some of the ARCH coefficients are not time-varying. 
For $t\in\[p+1,T\]=[p+1,T]\cap \N$, let $M_t$ and $N_t$ be
two random vectors of size $m$ and $n$ respectively, 
with $m+n=p+1$ and defined as follows. We split the interval $\[0,p\]$ into two parts $\{q_1,q_2,\ldots,q_m\}$ and $\{r_1,r_2,\ldots,r_n\}$ with $q_1<\cdots<q_m$ and $r_1<\cdots<r_n$.  
If $q_1=0$, we set $M_t=\left(1,X^2_{t-q_2},\ldots,X^2_{t-q_m}\right)'$, with the convention $M_t=1$ if $m=1$. 
If $q_1>0$, we set $M_t=\left(X^2_{t-q_1},\ldots,X^2_{t-q_m}\right)'$. The vector $N_t$ is defined 
similarly, replacing the $q_{\ell}$'s with the $r_{\ell}$'s.
In particular, the coordinates of the random vectors $M_t$ and $N_t$ form a bipartition of the set $\left\{1,X^2_{t-1},\ldots,X_{t-p}^2\right\}$.   
Now, we assume that the coefficients vector $\beta=\left(a_{r_1},\ldots,a_{r_n}\right)'$ is constant.
We also set $\alpha_t=\left(a_{q_1}(t/T),\ldots,a_{q_m}(t/T)\right)'$. 
Then the model writes  
\begin{equation}\label{SParam}
X_t=\xi_t\sigma_t,\quad \sigma_t^2=M_t'\alpha_t+N_t'\beta,
\end{equation}
for $t=p+1,\ldots, T$.  
Throughout this paper, we will assume that for all realization $\omega$ in the probability space, $\left(X_t(\omega)\right)_{t\leq 0}$ is a path of a stationary ARCH process with noise $\xi$ and 
coefficients $a_j(0)$. From this convention, one can get a local approximation of a tv-ARCH process by stationary ARCH processes with parameters $\left(a_0(u),\ldots,a_p(u)\right)$. See Lemma $1$ in the supplementary material for details.      

In Section \ref{S4}, a statistical test will be given for testing $H_0$: $\left(a_{r_1},\ldots,a_{r_n}\right)$ is constant. 
The sequel of this section is devoted to the statistical inference in model (\ref{tvARCH}) under the null hypothesis $H_0$. The corresponding estimation of parameter $\beta$ will be necessary to construct the test. 
\subsection{Estimators of the parametric part $\beta$}
Considering the square of the process (\ref{SParam}), statistical inference in model (\ref{SParam}) can be viewed as a linear regression problem. More precisely, we have for $t\geq p+1$, 
\begin{equation}\label{reg}
X^2_t=M_t'\alpha_t+N_t'\beta+(\xi^2_t-1)\sigma_t^2.
\end{equation}
In the sequel, we consider a sequence of weights $\left(W_t\right)_{p+1\leq t\leq T}$ such that $W_t$ is a measurable function of $t,T$ and $X_{t-1}^2,\ldots,X_{t-p}^2$. For stating our results, 
we will only consider sequences of the form
\begin{equation}\label{poids}
W_t=\left(\gamma_0\left(t/T\right)+\sum_{\ell=1}^p \gamma_j\left(t/T\right)X_{t-j}^2\right)^{-2},
\end{equation}
where the $\gamma_j$'s are positive and Lipschitz continuous functions defined over $[0,1]$.
The use of this kind of weights is classical in weighted least squares estimation in order to relax moment conditions on the marginal distribution or to gain in efficiency. 
The first goal of our procedure is to estimate the parameter $\beta$. This is 
the most difficult part of our methodology since a $\sqrt{T}-$consistent estimate is expected. 
Once an estimate $\hat{\beta}$ with the classical parametric rate of convergence is available, a pointwise estimate of parameter $\alpha(\cdot)$ can easily be obtained. One can just plug $\hat{\beta}$
in (\ref{reg}) and apply standard nonparametric methods (in this paper we will use the local weighted least-squares method studied in \citet{Fryz}).   
Our aim here is to first eliminate the nonparametric component $M_t'\alpha_t$. 
Our approach is classical in the setting of partially linear models, for which the regression function involves a parametric component and a nonparametric component (see for example \citet{PLM}, Chapter $6$ for some results in the case of time series). However, our method is based
on nonparametric estimation of linear projections of $\sqrt{W_t}X_t^2$ and $\sqrt{W_t}N_t$ onto the $\L^2$ subspace generated by the components of $\sqrt{W_t}M_t$ instead of a nonparametric estimation of the conditional expectations. Moreover,
our two-step approach involving some weights and leading to semiparametric efficient estimates is not common for nonstationary time series and no existing results from the theory of partially linear models can be used here for our purpose. 
Here, our approach can be also interpreted as a partial regression.  For stationary ARCH processes, estimation of the whole set of parameters
using a regression model for the squares and least squares estimation has been studied by \citet{Bose} and \citet{Hor}.
  
Now we introduce our estimator. We first multiply the two members of equation (\ref{reg}) by $\sqrt{W_t}$. 
If $\mathcal{P}_t\left(\sqrt{W_t}X_t^2\right)$ and $\mathcal{P}_t\left(\sqrt{W_t}N_t\right)$ denote the (componentwise) orthogonal projection of $\sqrt{W_t}X_t^2$ and $\sqrt{W_t}N_t$ onto the $\L^2$ linear subspace generated by the coordinates of $\sqrt{W_t}M_t$, it is easily seen that 
\begin{equation}\label{proproj}
\sqrt{W_t}X_t^2-\mathcal{P}_t\left(\sqrt{W_t}X_t^2\right)=\left(\sqrt{W_t}N_t-\mathcal{P}_t\left(\sqrt{W_t}N_t\right)\right)'\beta+(\xi^2_t-1)\sqrt{W_t}\sigma_t^2.
\end{equation}
The use of these orthogonal projections are natural in order to eliminate the nuisance parameter $\alpha_t$ and to get a partial regression involving parameter $\beta$ only.  
Let us introduce some notations for expressing these projections.  
For $t\geq p+1$, we set 
$$q_{1,t}=s_{3,t}^{-1}s_{1,t},\quad q_{2,t}=s_{3,t}^{-1}s_{2,t}$$
where
$$s_{3,t}=\E\left(W_tM_tM_t'\right),\quad s_{1,t}=\E\left(W_tM_t X^2_t\right),\quad s_{2,t}=\E\left(W_tM_t N_t'\right),$$
Then setting  
$$V_t=X^2_t-M_t'q_{1,t},\quad O_t=N_t-q_{2,t}'M_t,$$
we have 
$$\sqrt{W_t}X_t^2-\mathcal{P}_t\left(\sqrt{W_t}X_t^2\right)=\sqrt{W_t}V_t,\quad \sqrt{W_t}N_t-\mathcal{P}_t\left(\sqrt{W_t}N_t\right)=\sqrt{W_t}O_t$$
and equation (\ref{proproj}) writes 
$$\sqrt{W_t}V_t=\sqrt{W_t}O_t\beta+\left(\xi_t^2-1\right)\sqrt{W_t}\sigma_t^2.$$
The idea is now to use a least squares estimator for $\beta$. Of course, in order to obtain a feasible estimator, it is necessary to first estimate the two quantities $q_{1,t}$ and $q_{2,t}$. To this end, we consider 
$$\hat{s}_{3,b,t}=\sum_{i=p+1}^T k_{t,i}(b)W_iM_iM_i',\quad \hat{s}_{1,b,t}=\sum_{i=p+1}^Tk_{t,i}(b)W_iM_iX_i^2,$$
$$\hat{s}_{2,b,t}=\sum_{i=p+1}^Tk_{t,i}(b)W_iM_iN_i',\quad k_{t,i}(b)=\frac{K\left(\frac{t-i}{Tb}\right)}{\som_{i=p+1}^TK\left(\frac{t-i}{Tb}\right)},$$
where for $K$ is a kernel and $b>0$ is a bandwidth parameter. Throughout this paper, the kernel $K$ is assumed to be absolutely continuous and 
with support $[-1,1]$. 
Then we set
\begin{equation}\label{quotientbase}
\hat{q}_{1,b,t}=\left(\hat{s}_{3,b,t}\right)^{-1}\hat{s}_{1,b,t},\quad \hat{q}_{2,b,t}=\left(\hat{s}_{3,b,t}\right)^{-1}\hat{s}_{2,b,t}.
\end{equation}
Now we introduce the following notations.
The quantities 
$$\hat{V}_t=X^2_t-M_t'\hat{q}_{1,b,t},\quad \hat{O}_t=N_t-\hat{q}_{2,b,t}'M_t$$
estimate $V_t$ and $O_t$ respectively.
Our estimator of parameter $\beta$ will be denoted by $\hat{\beta}$ and minimizes the function $\ell_W$ defined by 
$$\ell_W\left(\overline{\beta}\right)=\sum_{t=p+1}^TW_t\left(\hat{V}_t-\hat{O}_t'\overline{\beta}\right)^2.$$
We get
\begin{equation}\label{MC1}
\hat{\beta}=\left(\sum_{t=p+1}^TW_t\hat{O}_t\hat{O}'_t\right)^{-1}\sum_{t=p+1}^TW_t\hat{O}_t\hat{V}_t.
\end{equation}
It is now possible to define an estimate of parameter $\alpha(u)$ for $u\in [0,1]$ by minimizing the function 
$$\overline{\alpha}\mapsto \sum_{i=p+1}^TK\left(\frac{t-i}{Tb'}\right)W_i\left(X_i^2-M_i'\overline{\alpha}-N_i'\hat{\beta}\right)^2,$$
for an integer $t\in\llbracket 1,T\rrbracket$ such that $\left\vert u-\frac{t}{T}\right\vert\leq \frac{1}{T}$ (e.g $t=[Tu]$ where
$[x]$ denotes the integer part of a real number $x$).
Since the nonparametric estimation of the function $\alpha$ requires a less restrictive assumption on the bandwidth parameter than for the estimation of parameter $\beta$, we introduce a new bandwidth $b'$. 
This leads to the estimate
\begin{equation}\label{MC2}
\hat{\alpha}_t=\hat{q}_{1,b',t}-\hat{q}_{2,b',t}\cdot\hat{\beta},
\end{equation}

\paragraph{Note.}
Expression $\hat{q}_{1,b,t}$ and $\hat{q}_{2,b,t}$ involve the inverse of the matrix $S_{3,b,t}$. 
One can show that 
$$\P\left(\det\left(\hat{s}_{3,t}\right)=0\mbox{ for some } t\leq T\right)\rightarrow 0.$$ See Lemma $3$ and its proof. However invertibility problems can occur when the noise $\xi$ has a mass at point $0$, for instance.
For simplicity, we always assume that all these matrices are invertible. Studying our estimator on an event with probability tending to one 
only complicates the statements and proofs of our results by adding some indicator sets but does not change the used approach. One can also show that this distinction is unnecessary for a noise $\xi$ having a density.

Asymptotic normality of the estimates (\ref{MC1}) and (\ref{MC2}) will be derived in the next section as well as
some plug-in versions to get optimal asymptotic results.  

\section{Asymptotic results}\label{S3}

\subsection{Estimation of the parametric component}

Our first result shows that the estimator (\ref{MC1}) is $\sqrt{T}-$consistent under some conditions. Here are our main assumptions.
\begin{description}
\item[A1.] For $j=1,\ldots,p$, the function $a_j$ is non-negative and Lipschitz continuous. The function $a_0$ is positive and Lipschitz continuous. Moreover, 
$$c=\sup_{u\in [0,1]}\sum_{j=1}^pa_j(u)<1.$$
\item[A2(h).] For the integer $h\geq 2$, there exists a real number $\delta\in(0,1)$ such that $\E\vert\xi_1\vert^{h(1+\delta)}<\infty$.
\end{description}
Assumption {\bf A1} is the classical contraction condition used in \citet{Fryz} to define tv-ARCH processes. Assumption {\bf A2(h)}, used for different values of $h$ in the sequel, implies a restriction on the noise distribution. Let us mention that this condition does not restrict the moment condition for the marginal $X_t$ (in the stationary case, i.e when the $a_j$'s are deterministic, assumption {\bf A1} is the necessary and sufficient condition for the condition $\E\left(X_t^2\right)<+\infty$). 

In the sequel, $\left(X_t(u)\right)_t$ will denote the stationary ARCH process with coefficients $\left(\alpha(u),\beta\right)$.
Then we will use the notation $\left(M_t(u)\right)_t$ (resp. $\left(N_t(u)\right)_t$, $\left(W_t(u)\right)_t$) for the stationary approximation of $\left(M_t\right)_t$ (resp. $\left(N_t\right)_t$, $(W_t)_t$). 
For example, $$W_t(u)=\left(\gamma_0(u)+\som_{\ell=1}^p\gamma_{\ell}(u)X_{t-\ell}^2(u)\right)^{-2}.$$

\begin{theo}\label{Nas}
Assume that assumptions ${\bf A1}$ and ${\bf A2(4)}$ hold, $b\sqrt{T}\rightarrow \infty$ and\\ $b^2\sqrt{T}\rightarrow 0$.
Then we have the following convergence in distribution:
$$\sqrt{T}\left(\hat{\beta}-\beta\right)\rightarrow_{T\rightarrow \infty} \mathcal{N}_n\left(0,\Sigma_1^{-1}\Sigma_2\Sigma_1^{-1}\right),$$
with  
$$\Sigma_1=\E\int_0^1 W_1(u)\left(N_1(u)-q_2(u))'M_1(u)\right)\cdot\left(N_1(u)-q_2(u)'M_1(u)\right)'du,$$
$$\Sigma_2=\V\left(\xi_0^2\right)\cdot\E\int W_1(u)^2\sigma_1(u)^4\left(N_1(u)-q_2(u)'M_1(u)\right)\cdot\left(N_1(u)-q_2(u)'M_1(u)\right)'du,$$
and 
$$q_2(u)=\E^{-1}\left(W_1(u)M_1(u)M_1(u)'\right)\E\left(W_1(u)M_1(u)N_1(u)'\right).$$
\end{theo}
\paragraph{Notes.}
\begin{enumerate}
\item
The bandwidth conditions used in Theorem \ref{Nas} are classical for estimating the parametric component in partially linear models. With this restriction, the nonparametric estimation step involved in the expression of $\hat{\beta}$ becomes negligible (i.e for $i=1,2$, $\hat{q}_{i,b,t}$ can be replaced with $q_{i,t}$ without changing the asymptotic behavior of $\hat{\beta}$). Let us explain the rule of these conditions. Some nonparametric estimates are introduced to approximate the two ratio $\hat{q}_{1,t}$ and $\hat{q}_{2,t}$. But, up to $C/T$ ($C$ denotes a positive constant), this two ratio can be seen as some Lipschitz functions of $t/T$ (see Lemma $5$ in the supplementary material). The mean square error for the kernel estimation of a Lipschitz functional in a regression model with deterministic design is bounded by $b^2+\frac{1}{Tb}$ (up to a constant). Then our bandwidth conditions entail that this mean square error converges to zero with a faster rate than $\sqrt{T}$. 
\item
The goal of the proof of Theorem \ref{Nas} is to show that the asymptotic distribution of $\sqrt{T}\left(\hat{\beta}-\beta\right)$ 
is the same as if the two quantities $\hat{q}_{1,b,t}$, $\hat{q}_{2,b,t}$ are replaced by $q_{1,t}$, $q_{2,t}$ respectively. Hence, the control of sums involving differences between these quantities are shown to be negligible. To this end, we make Taylor expansions and bound the variance of some multiple weighted sums appearing in this expansion using Lemma $4$ given in the supplementary material.
\item
The asymptotic variance in Theorem $\ref{Nas}$ can be estimated consistently using the data.
Indeed, the proof of Theorem \ref{Nas} given in the supplementary material shows that 
$$\Sigma_1=\lim_{T\rightarrow +\infty}\frac{1}{T}\sum_{t=p+1}^T W_t\hat{O}_t\hat{O}_t' \mbox{ a.s}.$$ 
Moreover, we have
\begin{equation}
\Sigma_2=\E\int W_1(u)^2\left(X_1(u)^2-\sigma_1(u)^2\right)^2\left(N_1(u)-q_2(u)'M_1(u)\right)\cdot\left(N_1(u)-q_2(u)'M_1(u)\right)'du.
\end{equation}
Then an estimate of $\Sigma_2$ can be obtained if we replace $q_2$ with $\hat{q}_{2,b,\cdot}$, $\sigma_t^2$ with a pointwise estimate $\hat{\sigma}^2_t$ (see Theorem \ref{Nas2}) and using the same kind of empirical counterpart as for $\Sigma_1$. 

\end{enumerate} 
The asymptotic variance given in Theorem \ref{Nas} depends on some weights $W_t(u)$.
One can show that its minimal value (in the sense of non-negative definite matrices) is obtained for the choice $W^{*}_t=\frac{1}{\sigma_t^4}$.
Indeed, setting $O_1(u)=N_1(u)-q_2(u)'M_1(u)$ and $O_1^{*}(u)=N_1(u)-q_2^{*}(u)'M_1(u)$, where
$$q_2^{*}(u)=\E^{-1}\left(\frac{M_1(u)M_1(u)'}{\sigma_1(u)^4}\right)\E\left(\frac{M_1(u)N_1(u)'}{\sigma_1(u)^4}\right),$$
we have for all $u\in[0,1]$,
$$\Sigma_1= \E\int_0^1W_1(u)O_1(u)O_1^{*}(u)'du= \E\int_0^1W_1(u)\sigma^2_1(u)O_1(u)\frac{O_1^{*}(u)'}{\sigma^2_1(u)}du.$$
Then if $x,y\in\R^n$, we have from the Cauchy-Schwarz inequality
$$\left(x'\Sigma_1y\right)^2\leq \frac{x'\Sigma_2x}{\V\left(\xi_0^2\right)}\times y'\Sigma y.$$
where $\Sigma=\E\int_0^1 \frac{O_1^{*}(u)O_1^{*}(u)'}{\sigma_1^4(u)}du$.
Now setting $x=\Sigma_1^{-1}z$ and $y=\Sigma^{-1}z$ for $z\in \R^k$, we get
$$\left(z'\Sigma^{-1}z\right)^2\leq \frac{z'\Sigma_1^{-1}\Sigma_2\Sigma_1^{-1}z}{\V\left(\xi_0^2\right)}\times z'\Sigma^{-1}z.$$
Then we have proved the following result.

\begin{Prop}\label{borneinferieure}
The lower bound for the asymptotic variance given in Theorem \ref{Nas} is $\V\left(\xi_0^2\right)\Sigma^{-1}$, where 
$$\Sigma=\int_0^1\E\left(\frac{1}{\sigma_1(u)^4}\left(N_1(u)-q^{*}_2(u)'M_1(u)\right)\left(N_1(u)-q^{*}_2(u)'M_1(u)\right)'\right)du,$$
where 
$$q_2^{*}(u)=\E^{-1}\left(\frac{M_1(u)M_1(u)'}{\sigma_1(u)^4}\right)\E\left(\frac{M_1(u)N_1(u)'}{\sigma_1(u)^4}\right).$$
\end{Prop}

Now, we show that it is possible to construct an estimate of parameter $\beta$ which has the asymptotic variance given in Proposition \ref{borneinferieure}. A natural candidate is obtained by replacing the weights $W_t$ in (\ref{MC1}) with an estimation of the optimal weights $W_t^{*}=\frac{1}{\sigma_t^4}$.
We set 
$\hat{W}^{*}_t=\frac{1}{\hat{\sigma}_t^4+\nu_T},$ 
where
$$\hat{\sigma}_t^2=M_t'\left(\hat{q}_{1,b,t}-\hat{q}_{2,b,t}\hat{\beta}\right)+N_t'\hat{\beta}$$
is an estimator of $\sigma_t^2$. The sequence $(\nu_T)_T$ is a sequence of positive real numbers such that $\nu_T=o\left(\frac{1}{\sqrt{T}}\right)$. 
The use of this sequence is just technical and avoids possible small values for the fitted volatility $\hat{\sigma}^2_t$ which is not ensured to be bounded away from $0$ for finite samples. However, for a large sample size, our simulations show that the choice $\nu_T=0$ does not alter the performance of the plug-in estimate. Let us define the quantities
$$\hat{s}^{*}_{3,b,t}=\sum_{i=p+1}^Tk_{t,i}(b)\hat{W}^{*}_iM_iM_i',\quad \hat{s}^{*}_{1,b,t}=\sum_{i=p+1}^Tk_{t,i}(b)\hat{W}^{*}_iM_iX_i^2,$$ $$\hat{s}^{*}_{2,b,t}=\sum_{i=p+1}^Tk_{t,i}(b)\hat{W}^{*}_iM_iN_i'.$$
We also set
$$\hat{q}^{*}_{1,b,t}=\left(\hat{s}^{*}_{3,b,t}\right)^{-1}\hat{s}^{*}_{1,b,t},\quad \hat{q}^{*}_{2,b,t}=\left(\hat{s}^{*}_{3,b,t}\right)^{-1}\hat{s}^{*}_{2,b,t}.$$
Now we introduce the following notations in order to simplify the expression of our estimator.
$$\hat{V}^{*}_t=X^2_t-M_t'\hat{q}^{*}_{1,b,t},\quad \hat{O}^{*}_t=N_t-\left(\hat{q}^{*}_{2,b,t}\right)'M_t.$$
Our plug-in estimate of parameter $\beta$ is now defined by 
$$\hat{\beta}_{*}=\left(\sum_{t=p+1}^T\hat{W}^{*}_t\hat{O}^{*}_t\left(\hat{O}^{*}_t\right)'\right)^{-1}\sum_{t=p+1}^T\hat{W}^{*}_t\hat{O}^{*}_t\hat{V}^{*}_t.$$
With respect to Theorem \ref{Nas}, we impose more restrictive assumptions.

\begin{description}
\item[A3.] $\xi_1$ has moments of any order.
\item[A4.] For $j=0,\ldots,p$, the coefficient $a_j$ is a positive function.  
\end{description}

\begin{theo}\label{Nasopt}
Assume that the assumptions ${\bf A1}$, ${\bf A3}$ and ${\bf A4}$ hold and 
$b T^{\frac{1}{2}-\tau}\rightarrow \infty$, $b^2\sqrt{T}\rightarrow 0$ for some $\tau\in (0,1/4)$. Then we have
$$\sqrt{T}\left(\hat{\beta}_{*}-\beta\right)\rightarrow_{T\rightarrow \infty} \mathcal{N}_n\left(0,\left(\E\xi_0^4-1\right)\Sigma^{-1}\right),$$
\end{theo}
The proof of Theorem \ref{Nasopt} is similar to that of Theorem \ref{Nas} but involves more tedious arguments. A detailed proof of Theorem \ref{Nasopt} is given in the supplementary material. 

\subsection{Estimation of the nonparametric component}
Now let us investigate the asymptotic properties for time-varying coefficients estimate $\hat{\alpha}$ defined by (\ref{MC2}).
The estimator $\hat{\beta}$ appearing in the expression (\ref{MC2}) is constructed using the initial bandwidth parameter $b$ which satisfies the assumptions of Theorem \ref{Nas}. 

\begin{theo}\label{Nas2}
Let $u\in [0,1]$. Assume that the assumptions of Theorem \ref{Nas} hold. Then if $b'\rightarrow 0$, $b'T\rightarrow \infty$ and $t=t_T$ satisfies $\vert\frac{t}{T}-u\vert\leq \frac{1}{T}$,
$$\sqrt{Tb'}\left(\hat{\alpha}_t-\alpha(u)\right)+\sqrt{Tb'}\E^{-1}\left(W_1(u)M_1(u)M_1(u)'\right)\left(A_t(u)-A^{\#}_t(u)\right)\rightarrow \mathcal{N}_m\left(0,\mathcal{V}(u)\right),$$
where
\begin{description}
\item $\mathcal{V}(u)=\V\left(\xi_1^2\right)\cdot\int K(x)^2dx\\\cdot\E^{-1}\left(W_1(u)M_1(u)M_1(u)'\right)\E\left(W_1(u)^2\sigma_1(u)^4M_1(u)M_1(u)'\right)\E^{-1}\left(W_1(u)M_1(u)M_1(u)'\right)$.
\item
$A_t(u)=\hat{s}_{1,b',t}-\hat{s}_{2,b',t}\beta-\hat{s}_{3,b',t}\alpha(u),$
\item
$A^{\#}_t(u)=\sum_{i=p+1}^Tk_{t,i}(b')W_i(u)M_i(u)\left(X_i(u)^2-N_i'(u)\beta-M_i(u)'\alpha(u)\right)$. 
\end{description}
\end{theo}
This result is similar to that obtained in \citet{Fryz}, Proposition $3$. 
The part $A_t(u)-A^{\#}_t(u)$ can be interpreted as a term of deviation with respect to stationarity. As pointed out in \citet{Fryz}, this term satisfies $A_t(u)-A^{\#}_t(u)=O_{\P}(b)$. 
One can easily check that the optimal asymptotic variance in Theorem \ref{Nas2} corresponds to the choice $W^{*}_t=\frac{1}{\sigma_t^2}$ for the weights. Thus, a plug-in approach is natural. We set
$\check{W}^{*}_{t,i}=\frac{1}{\hat{\sigma}_{t,i}^4+\mu_T}$, where 
$$\hat{\sigma}^2_{t,i}=M_i'\hat{\alpha}_t+N_i'\hat{\beta}$$
and $(\mu_T)$ is a sequence of positive real numbers which now plays the rule of the sequence $(\nu_T)_T$  previously used for the optimal estimator $\hat{\beta}_{*}$. Once again, this choice is only technical and we impose here $\mu_T=O\left(b'+\frac{1}{\sqrt{Tb'}}\right)$. Then we define the following estimator of parameter $\alpha_t$.
$$\hat{\alpha}_{*,t}=\check{s}_{3,b',t}^{-1}\left(\check{s}_{1,b',t}-\check{s}_{2,b',t}\hat{\beta}\right),$$
where for $j=1,2,3$, $\check{s}_{j,b',t}$ is obtained as $\hat{s}_{j,b',t}$ but replacing $W_i$ with $\check{W}^{*}_{t,i}$.

\begin{theo}\label{Nas2opt}
Assume that assumptions ${\bf A1}$, ${\bf A2(4)}$ and ${\bf A4}$ hold, $b'T\rightarrow \infty$ and
$\sqrt{Tb'}\left(b'\right)^{\ell_0}\rightarrow 0$ for a given integer $\ell_0$.
Then, if $\vert\frac{t}{T}-u\vert\leq \frac{1}{T}$, we have
$$\sqrt{Tb'}\left(\hat{\alpha}_{*,t}-\alpha(u)\right)+\sqrt{Tb'}\E^{-1}\left(\frac{M_1(u)M_1(u)'}{\sigma_1(u)^4}\right)\left(B_t(u)-B^{\#}_t(u)\right)\rightarrow \mathcal{N}_m\left(0,\mathcal{V}_{*}(u)\right),$$
where
\begin{description}
\item 
$\mathcal{V}_{*}(u)=\V\left(\xi_1^2\right)\cdot\int K(x)^2dx\cdot\E^{-1}\left(\frac{M_1(u)M_1(u)'}{\sigma_1(u)^4}\right),$
\item
$B_t(u)=\check{s}_{1,b',t}-\check{s}_{2,b',t}\beta-\check{s}_{3,b',t}\alpha(u),$
\item
$B^{\#}_t(u)=\sum_{i=p+1}^Tk_{t,i}(b')\check{W}^{*}_i(u)M_i(u)\left(X_i(u)^2-N_i'(u)\beta-M_i(u)'\alpha(u)\right)$. 
\end{description}
Moreover, $B_t(u)-B^{\#}_t(u)=O_{\P}(b')$.
\end{theo}

\paragraph{Notes}
\begin{enumerate}
\item
Compared to Theorem \ref{Nas2}, Theorem \ref{Nas2opt} uses a more restrictive assumption $b'$. However for powers of the sample size, i.e $b'=CT^{-\ell}$ with constants $C,\ell>0$, the conditions are equivalent. 
\item
When all the coefficients of the volatility are time-varying, replacing $M_1$ by $\left(1,X_0,\ldots,X_{-p+1}\right)'$, we recover
the expression of the optimal asymptotic variance given in \citet{Fryz} for the (local) weighted least-squares estimation. 
This asymptotic variance coincides with that obtained with the local QML estimator studied in \citet{DR}.
However, a crucial assumption in Theorem \ref{Nas2opt} is the positivity of all the coefficients of the volatility.
To avoid this restriction, \citet{Fryz} consider the sequence of weights $\hat{W}_t=\left(\hat{d}_t+\som_{j=1}^pX_{t-j}^2\right)^{-2}$
where $\hat{d}_t$ is a nonparametric estimate of $\E X_t^2$. For the nonparametric estimation of the whole set of coefficients, they show that the corresponding weighted least squares estimator is asymptotically normal, even if the ARCH coefficients are only nonnegative, 
but at the price of a small loss of efficiency. We claim that the weights $\hat{W}_t$ can also be used in our context to obtain a result similar to \citet{Fryz}, Proposition $4$. Details are omitted.
\item
As usual for this kind of nonparametric estimation, a better finite sample approximation of the distribution of the parameter estimators can often be obtained using bootstrap methods. With straightforward modifications, it is possible to use the bootstrap method studied in \citet{Fryz} (see Section $5$ of that paper) to obtain pointwise confidence intervals for the components of $\alpha$. Since this paper
is mainly devoted to testing and estimating some non-time varying coefficients, we will not consider this bootstrap scheme. 
\item
Detailed proofs of Theorem \ref{Nas2} and Theorem \ref{Nas2opt} are given in the supplementary material.    
\end{enumerate}

\subsection{Asymptotic semiparametric efficiency}
For Gaussian inputs (i.e, $\xi\sim\mathcal{N}(0,1)$), it is possible to show that the matrix $2\Sigma^{-1}$ given in Proposition \ref{borneinferieure} is a lower bound in semiparametric estimation. 
We refer the reader to \citet{Bic} for a general introduction to semiparametric models and the problem of efficient estimation of 
a finite dimensional parameter in such models. In our case case, the problem of semiparametric efficiency for estimating the parameter $\beta$ involves triangular arrays. This is why we will use an abstract result using the classical formalism presented in \citet{VWW} (see Chapter $3.11$). Intuitively, one can see the matrix $2\Sigma^{-1}$ as the smallest asymptotic variance obtained for estimating $\beta$ in submodels for which the nuisance parameter $\alpha_t$ is projected onto a finite dimensional space of square integrable functions. 
Formally, the approach consists in writing a LAN expansion of the likelihood ratio and then using a general convolution theorem.  
In the sequel, we set for $m,n\geq 1$, $\H=\L^2\left([0,1]\right)^m\times \R^n$. Then $\H$ is an Hilbert space for the classical scalar product
$$<(g,h);(\bar{g},\bar{h})>_1=\sum_{i=1}^m \int_0^1 g_i(u)\bar{g}_i(u)du +\sum_{j=1}^n h_j\bar{h}_j.$$
However, in the sequel, the space $\H$ will be endowed with an equivalent scalar product defined by 
$$<(g,h)\vert (\bar{g},\bar{h})>_{\H}=\frac{1}{2}\int_0^1\begin{pmatrix}g(u)\\h\end{pmatrix}' E(u)\begin{pmatrix} g(u)\\h\end{pmatrix}du,$$
where 
$$E(u)=\begin{pmatrix} \E\left(\frac{M_1(u)M_1(u)'}{\sigma_1(u)^4}\right)&\E\left(\frac{M_1(u)N_1(u)'}{\sigma_1(u)^4}\right)\\
\E\left(\frac{N_1(u)M_1(u)'}{\sigma_1(u)^4}\right)&\E\left(\frac{N_1(u)N_1(u)'}{\sigma_1(u)^4}\right)\end{pmatrix}.$$
Now, we denote by $\mathcal{L}$ the set of Lipschitz functions $f:[0,1]\rightarrow \R$ and we set $H=\mathcal{L}^m\times \R^n$ where $m$ (resp. $n$) is the dimension of vector $M_t$ (resp. $N_t$).
Then $H$ is a linear subspace of $\H$. The set $\left(H,<\cdot,\cdot>_{\H}\right)$ will be referred to the tangent space. 
For Gaussian inputs, we first derive a LAN expansion for the (conditional) likelihood ratio. 
We denote by $P_{T,\alpha,\beta}$ the conditional distribution  $\left(X_{p+1},\ldots,X_T\right)\vert X_1=x_1,\ldots,X_p=x_p$.

\begin{Prop}\label{LAN}
Assume that $\left((\alpha,\beta);(g,h)\right)\in H^2$ where the coordinates of $\alpha$ and $\beta$ are positive. 
Then we have
$$\log\frac{d\P_{T,\alpha+\frac{g}{\sqrt{T}},\beta+\frac{h}{\sqrt{T}}}}{d\P_{T,\alpha,\beta}}\left(X_{p+1},\ldots,X_T\right)=\Delta_{T,g,h}-\frac{1}{2}\Vert (g,h)\Vert_{\H}^2+o_{\P_{T,\alpha,\beta}}(1),$$
where
$$\Delta_{T,g,h}=\frac{1}{2\sqrt{T}}\sum_{t=p+1}^T\frac{X_t^2-\sigma_t^2}{\sigma_t^4}\left(M_t'g\left(\frac{t}{T}\right)+N_t'h\right).$$
Moreover, 
$$\Delta_{T,g,h}\stackrel{\mathcal{D}}{\rightarrow}\mathcal{N}_n\left(0,\Vert (g,h)\Vert_{\H}^2\right).$$
\end{Prop}

From this LAN expansion, we derive a lower bound for the asymptotic variance of regular estimators of $\beta$. 

For $(g,h)\in H$, we set $\kappa_T(g,h)=\beta+\frac{h}{\sqrt{T}}$. Then, the sequence of parameters $\left\{\kappa_T(g,h): (g,h)\in H\right\}$ is regular: if $\dot{\kappa}:\H\rightarrow \R^n$ is the projection operator defined by $\dot{\kappa}(g,h)=h$, then
$$\sqrt{T}\left(\kappa_T(g,h)-\kappa_T(0,0)\right)=h.$$

\begin{cor}\label{optimal}
If the assumptions of Proposition \ref{LAN} hold, then
the adjoint operator $\dot{\kappa}^{*}:\R^n\rightarrow \H$ of $\dot{\kappa}$ is given by 
$$\dot{\kappa}^{*}v=\begin{pmatrix}-q_2(\cdot)\\I_n\end{pmatrix} \left(\frac{1}{2}\Sigma\right)^{-1}v,\quad v\in\R^n.$$
Consequently,  the limit distribution of a regular estimator of $\beta$ equals the distribution of a sum $L_1+L_2$ of independent random vectors of $\R^n$ and such that
$$L_1\sim\mathcal{N}_n\left(0,2\Sigma^{-1}\right).$$ 
\end{cor}

\section{Statistical testing}\label{S4}

\subsection{Testing parameter constancy}
For a real data set, it is necessary, before applying the methodology given in Section \ref{S2}, to test if 
a coefficients vector of the form $\beta=\left(a_{r_1},\ldots,a_{r_n}\right)$ is time-varying or not in model (\ref{tvARCH}).
This is equivalent to test model (\ref{SParam}) versus model (\ref{tvARCH}).
When $n=p$ and $\beta=\left(a_1,\ldots,a_p\right)$, such a statistical test is interesting for deciding if model (\ref{lagconst}) 
is a convenient restriction of the tv-ARCH model. This case is of particular interest for real data applications. 
In \citet{ZW}, a procedure is proposed for testing if some coefficients are constant in a general time-varying regression model.
The null hypothesis is $H_0$: $\beta(\cdot)$ constant. 
The test statistic used in \citet{ZW} is based on a $\L^2$ distance between an estimate under the alternative and an estimate under the null hypothesis. 
In this part, we derive asymptotic properties of this test for tv-ARCH processes. 
For simplicity, we will only consider some estimates without plug-in (i.e we fix a sequence of weights $(W_t)_t$ of the form 
(\ref{poids}) and use the corresponding least-squares estimates).    
Let us first introduce some additional notations.

For a function $f:[-1,1]\rightarrow \R$ , we set $\Vert f\Vert_2=\sqrt{\int_{-1}^1 f(u)^2du}$ and for $x\in[-1,1]$, $$K^{*}(x)=\int_{-1}^{1-2\vert x\vert}K(v)K\left(v+2\vert x\vert\right)dv.$$
Setting for $u\in[0,1]$, $e_i(u)=\frac{1}{Tb}K\left(\frac{uT-i}{Tb}\right)$ and for $p+1\leq t\leq T$, $\mathcal{X}_t=\left(M_t',N_t'\right)'$, the kernel estimate of the full vector of ARCH coefficients $a(u)=\left(\alpha(u)',\beta(u)'\right)'$ is given by (see \citet{Fryz}) 
$$\widetilde{a}(u)=S_u^{-1}\sum_{i=p+1}^T e_i(u)W_iX_i^2\mathcal{X}_i,$$
where $S_u=\sum_{i=p+1}^Te_i(u)W_i\mathcal{X}_i\mathcal{X}_i'$. Then we set $\widetilde{\beta}(u)=A\widetilde{a}(u)$ where $A$ is the matrix of size $(p+1)\times n$ defined by $Ax=\begin{pmatrix}x_{m+1}\\\vdots\\x_{p+1}\end{pmatrix}$. Note also that $\beta(u)=A a(u)$.   
We also set $\kappa_u=\E\left(W_1(u)\mathcal{X}_1(u)\mathcal{X}_1(u)'\right)$ and
\begin{eqnarray*} 
\mathcal{O}(u)&=&\kappa_u^{-1}\cdot\E\left(W_1(u)^2\left(X_1(u)^2-\sigma^2_1(u)\right)^2\mathcal{X}_1(u)\mathcal{X}_1(u)'\right)\cdot \kappa_u^{-1}\\
&=&\V\left(\xi_1^2\right)\kappa_u^{-1}\cdot \E\left(W_1(u)^2\sigma_1(u)^4\mathcal{X}_1(u)\mathcal{X}_1(u)'\right)\cdot \kappa_u^{-1}.\\
\end{eqnarray*}
Let $\left(\Gamma(u)\right)_{u\in [0,1]}$ be a family of positive definite matrices of size $(p+1)\times (p+1)$ such that $u\mapsto \Gamma(u)$ is a Lipschitz function.  
Finally, we set for $j=1,2$, 
$$\varpi_j=\int_0^1\Tr\left[\Gamma(u)^{1/2}A\mathcal{O}(u)A'\Gamma(u)^{1/2}\right]^jdu.$$
We define our test statistic $\mathcal{S}_T\left(\widetilde{a},\beta\right)$ by
\begin{equation}\label{statbase}
\mathcal{S}_T\left(\widetilde{a},\beta\right)=\int_0^1 \left(\widetilde{\beta}(u)-\beta\right)'\Gamma(u)\left(\widetilde{\beta}(u)-\beta\right)du.
\end{equation}
The proof of the following theorem is given in the supplementary material.
\begin{theo}\label{testing}
Assume that assumptions ${\bf A1}$, ${\bf A2(8)}$ hold and $\beta$ is non time-varying. Then if $Tb^2\rightarrow \infty$ and $Tb^{3.5}\rightarrow 0$, we have
\begin{equation}\label{result}
T\sqrt{b}\left\{\mathcal{S}_T\left(\widetilde{a},\beta\right)
-\frac{\Vert K\Vert_2^2\varpi_1}{Tb}\right\}\rightarrow \mathcal{N}\left(0,4\Vert K^{*}\Vert_2^2\varpi_2\right).
\end{equation}
Moreover (\ref{result}) holds if $\mathcal{S}_T\left(\widetilde{a},\beta\right)$ is replaced with $\mathcal{S}_T\left(\widetilde{a},\hat{\beta}\right)$ where $\hat{\beta}$ is the estimate of Theorem \ref{Nas}. 
\end{theo}

\paragraph{Notes}
\begin{enumerate}
\item
As pointed out in \citet{ZW}, if we are interested in prediction, the matrix $\Gamma$ can be chosen as the asymptotic variance of the kernel estimate $\hat{\beta}(\cdot)$ (which has to be estimated in practice). In our numerical studies, we will use the simple choice $\Gamma(u)=I_n$ where $I_n$ denotes the identity matrix of size $n$. 
\item
Quantities $\varpi_1$ and $\varpi_2$ involved in the bias and asymptotic variance in (\ref{result}) 
can be estimated consistently, taking empirical counterpart. Then we obtain a pivotal statistic \\
$\mathcal{E}_T=T\sqrt{b}\left\{\mathcal{S}_T\left(\widetilde{a},\hat{\beta}\right)
-\frac{\Vert K\Vert_2^2\hat{\varpi_1}}{Tb}\right\}/\left(2\Vert K^{*}\Vert_2\sqrt{\hat{\varpi_2}}\right)$ and one can reject the null hypothesis for large values of this statistic. However, in practice, such nonparametric tests suffer from the slow convergence in Theorem \ref{testing}.  As in \citet{ZW}, one can use a Monte-Carlo type procedure which can improve the finite-sample performance (a similar Monte Carlo procedure is also used in \citet{Pat}). Note that the result of Theorem \ref{testing}
is valid for i.i.d series with a standard Gaussian marginal distribution. In particular, if $\mathcal{E}_T^{*}$ denotes the pivotal statistics computed with an i.i.d sample of standard Gaussian variables, we have $\lim_{T\rightarrow\infty}\mathcal{E}_T^{*}=\lim_{T\rightarrow \infty}\mathcal{E}_T=\mathcal{N}(0,1)$ in distribution. Then one can use the quantiles of the distribution of $\mathcal{E}_T^{*}$ to compute the critical values for the test (instead of the Gaussian quantiles). Let us give the details of the method proposed in \citet{ZW}. We assume that the bandwidth $b$ has been already selected.
\begin{itemize}
\item
First simulate $B$ samples of size $T$ of i.i.d Gaussian random variables. For each sample, compute the values of the estimators $\widetilde{\beta}(\cdot)$ and $\hat{\beta}$ as well as the realization of the pivotal statistics $\mathcal{P}^{*}_T$.
\item
Then, from these $B$ realizations of the random variable $\mathcal{E}^{*}_T$, compute the empirical quantile $q_{MC}(\alpha)$ of order $1-\alpha$.
\item
Reject $H_0$ if $\mathcal{E}_T$ is greater than $q_{MC}(1-\alpha)$.
\end{itemize} 
\item
In \citet{ZW}, the power of the test under some local alternatives is studied (see Theorem $3.2$ of this paper). A similar result can be derived here under the assumptions of Theorem \ref{testing}. In particular, for some local alternatives of the form $\beta_1(t/T)=\beta+z_T f(t/T)$, with $1/\sqrt{T\sqrt{b}}=o(z_T)$ and $f$ a Lipschitz function defined over $[0,1]$, the power of the test still converges to $1$.  
\end{enumerate}

\subsection{Testing if a constant parameter is equal to zero}
In this part, we consider model (\ref{SParam}) and our goal is to test whether the vector $\beta$ is equal to zero.
Two approaches are discussed below.
\begin{enumerate}
\item
One possibility is to use the asymptotic normality of the estimator $\hat{\beta}$ given in Theorem \ref{Nas}. Under the null hypothesis $H_0$: $\beta=0$, 
the asymptotic distribution of $\sqrt{T}\hat{\beta}$ is that of a centered Gaussian vector with covariance matrix $\mathcal{V}=\Sigma_1^{-1}\Sigma_2\Sigma_1^{-1}$.
We have already discussed how to estimate 
the covariance matrix $\mathcal{V}$. If $\hat{\mathcal{V}}$ denotes such an estimate, 
the statistics $T\Vert\hat{\mathcal{V}}^{-\frac{1}{2}}\hat{\beta}\Vert^2$ is asymptotically distributed as a $\chi^2$ with $n$ degrees of freedom (here $\Vert\cdot\Vert$ denotes the euclidean norm on $\R^n$). 
As for the test given in the previous subsection, one can use a Monte Carlo method instead of using the quantiles of the asymptotic distribution (the convergence in distribution of the previous statistics is quite slow because of the incorporation of nonparametric kernel estimates). If a bandwidth $b$ is selected, one can simulate $B$ samples of Gaussian i.i.d random variables and compute the corresponding values of our statistics. From these values, we can compute the empirical quantile $q_{MC,\alpha}$ of order $1-\alpha$ and reject the null hypothesis if $T\Vert\hat{\mathcal{V}}^{-\frac{1}{2}}\hat{\beta}\Vert^2>q_{MC,\alpha}$. 
This test has an asymptotic level $\alpha$ and a power converging to $1$ under a fixed alternative. 
However such approach is not completely natural because the value $\beta=0$ is on the boundary of the parameter space, our test is similar to the bilateral test for testing the hypothesis $\beta=0$ in regression models and we ignore the sign of $\beta$. This will result in a loss of power and it is more natural to consider a statistics based on the random vector $\left(\sqrt{T}\max\left(\hat{\beta}_i,0\right)\right)_{1\leq i\leq n}$. The vector $\left(\max(\hat{\beta}_i,0)\right)_{1\leq i\leq n}$ will be called truncated least squares estimator. 
As discussed in \citet{FZak} for stationary ARCH processes, truncated least squares estimators are natural for testing if some lag coefficients are equal to zero.  
However, the limiting distribution of this truncated random vector is that of $\left(\max\left(Z_i,0\right)\right)_{1\leq i\leq n}$ where $Z=(Z_1,\ldots,Z_n)$ is a Gaussian vector, with dependent entries in general.  
Then, except if $\V(Z)$ is diagonal, it is not possible to get a pivotal statistics from truncated least squares estimators. Then the Monte-Carlo method used for testing parameter constancy cannot be applied. 
Note also that bootstrapping the model is not appropriated here because the bootstrap is generally inconsistent for testing a parameter on the boundary (see for instance \citet{And} for this problem).
However, when $\beta$ denotes the full vector of lag coefficients, it is possible to use truncated least squares and the Monte Carlo method, provided $W_t\equiv 1$.  This point is discussed below.
\item
When $\beta$ is the full vector of lag coefficients, the problem is to test model (\ref{Star}) versus model (\ref{lagconst}).
For testing if the lag coefficients are equal to zero in model (\ref{lagconst}), our test is based on the following result.
The following notation will be used. If $p+1\leq t\leq T$, we set $\hat{d}_t=\sum_{i=p+1}^Tk_{t,i}(b)X_i^2$. Note that $\hat{d}_t$ is an estimator of $\E X_t^2$.
\begin{Prop}\label{newtest}
Assume that {\bf A2(4)} holds and that $b\in [cT^{-h},CT^{-h}]$ with $\frac{1}{4}<h<\frac{3}{4}-\frac{1}{2(1+\delta)}$, where $c,C$ are positive constants. Then under $H_0$, 
$$\left(\hat{a}_1,\ldots,\hat{a}_p\right)'=\arg\min_{a_1,\ldots,a_p}\sum_{t=p+1}^T\left(X_t^2-\hat{d}_t-\sum_{j=1}^p a_j\left(X_{t-j}^2-\hat{d}_{t-j}\right)\right)^2$$
satisfies 
$T\sum_{j=1}^p\max\left(\hat{a}_j,0\right)^2\rightarrow_{\mathcal{D}} \sigma^2\sum_{j=1}^p\max\left(Z_j,0\right)^2,$
where $Z$ standard Gaussian vector and $$\sigma^2=\frac{\int_0^1 \V\left(X_0(u)^2\right)^2 du}{\left(\int_0^1 \V\left(X_0(u)^2\right)du\right)^2}.$$  
\end{Prop}
\paragraph{Testing the second order dynamic.}
From this result, we reject $H_0$ for large values of the statistics $\Psi_T=T\sum_{j=1}^p\max\left(\hat{a}_j,0\right)^2/\hat{\sigma}^2$ where $\hat{\sigma}^2$ is a consistent estimator of $\sigma^2$. Typically one can choose 
$$\hat{\sigma}^2=T\sum_{t=1}^T \hat{a}_0(t/T)^4/\left(\sum_{t=1}^T\hat{a}_0(t/T)^2\right)^2\mbox{ with } \hat{a}_0(t/T)=\hat{d}_t,$$
which gives a consistent estimate under $H_0$. A first solution is to reject $H_0$ if $\Psi_T$ is larger than the quantile of order $1-\alpha$ of the distribution of the random variable $\sum_{j=1}^p\max\left(Z_j,0\right)^2$.
But the statistics $\Psi_T$ is also asymptotically pivotal and its quantiles 
can be approximated by a Monte Carlo procedure similar to that used for testing parameter constancy. 
\end{enumerate}

\paragraph{Notes}
\begin{enumerate}
\item
 Note that the estimators of the lag coefficients corresponds to the estimators (\ref{MC1})
with $W_t\equiv 1$. The optimal rate of convergence $h=\frac{1}{3}$ for kernel estimation of Lipschitz functionals can be used 
if $\delta>\frac{1}{5}$ in assumption ${\bf A2(4)}$. A proof of Proposition \ref{newtest} is given in the supplementary material.
\item
In the stationary case, two benchmark tests are usually used for testing the second order dynamic: the Lagrange multiplier test of \citet{Engle} and the portmanteau test of \citet{McL}. 
\citet{Pat} have recently extended these two tests for model (\ref{lagconst}), taking in account of nonstationarity. 
Here we provide an alternative test based on a direct estimation of the lag coefficients. A comparison of the different approaches
is beyond the scope of this paper. Let us observe that in the stationary case, the constant $\sigma^2$ in Proposition \ref{newtest}    
is equal to $1$ whereas in the nonstationary case, this constant $\sigma^2$ is a correction factor which can be written 
$$\sigma^2=\int_0^1a_0(u)^4du/\left(\int_0^1 a_0(u)^2du\right)^2>1.$$
The correction factor $\sigma^2$ also appears in the asymptotic results of \citet{Pat} (see the ratio $\overline{\omega}_4^2/\overline{\omega}_8$ appearing in the two statistics used in that paper). Ignoring this factor leads to an oversized test and the null hypothesis will be often rejected when the data are independent but not identically distributed. 
Moreover, let us notice that our moment condition for the noise distribution is less restrictive for applying the test (a moment greater than $8$ is assumed in \citet{Pat}). 
\item
One can study the power of our test under local alternatives of type $\sqrt{T}(a_1,\ldots,a_p)'\rightarrow s$ where $s$ is a vector of nonnegative real numbers. One can show that if $\E\xi_0^{8(1+\delta)}<\infty$
and the bandwidth $b$ satisfies $Tb^2\rightarrow \infty$ and $Tb^4\rightarrow 0$, then 
$$\P\left(T\sum_{j=1}^p\max\left(\hat{a}_j,0\right)^2/\hat{\sigma}^2>q_{1-\alpha}\right)\rightarrow \P\left(\sum_{j=1}^p\max\left(\frac{s_j}{\sigma}+Z_j,0\right)^2>q_{1-\alpha}\right),$$
where $q_{1-\alpha}$ is the quantile of order $1-\alpha$ of the distribution of $\sum_{j=1}^p\max(Z_j,0)^2$.
Details are omitted.
\end{enumerate}
   
\section{Real data applications}\label{S5}
Before real data applications, let us provide some recommendations on the choice of tuning parameters.
\subsection{Choice of tuning parameters}
The practical implementation of our estimation and testing procedures requires the 
choice of some weights $W_t$, some bandwidth parameters as well as the number of lags $p$ in the model. 
In this subsection, we discuss the practical choices of these parameters.
In all our studies, the kernel $W$ will be the Epanechnikov kernel.
For simplicity, only one bandwidth parameter will be selected for the semiparametric models (this means that we set $b'=b$ in (\ref{MC2})). 
We use a cross-validation method as specified below. The selected bandwidth will be used for the tests.
For the tests, the Monte-Carlo procedure will be always applied with $B=2000$ samples
of i.i.d standard Gaussian random variables.
\begin{itemize}
\item
In practice, a sequence of weights $(W_t)_t$ has to be chosen for applying our method. One possibility is to use the weights $W_t=\left(1+\som_{j=1}^pX_{t-j}^2\right)^{-2}$ suggested in \citet{Hor}. 
In our implementation, we use the weights $W_t=\left(\hat{v}+\sum_{j=1}^pX_{t-j}^2\right)^{-2}$ where $\hat{v}=\frac{1}{T}\sum_{t=1}^TX_t^2$ is an estimate of the average of the variance $v=\int_0^1\E\left(X_0(u)^2\right)du$.
There are several advantages in using these weights. First, the lag estimates obtained in model \ref{lagconst} do not depend on the scale of the returns, a property always satisfied for the true lag coefficients (if $W_t$ denotes the price at time $t$, $X_t=\log(P_t)-\log(P_{t-1})$ or $100\times X_t$ are two different scales used in practice). Moreover, for stationary Arch processes, this choice is equivalent to the weights used in \citet{Fryz}. We also noticed better finite sample performances for our tests and inference procedures with this choice. However, the introduction of the random quantity $\hat{v}$ is not taken in account in our theoretical results. But, inspection of our proofs shows that the conclusions of the theorems remain unchanged if $\hat{v}-v=o_{\P}\left(T^{-1/4}\right)$. The latter condition is satisfied if $\E\vert X_t\vert^h<\infty$ for $h>\frac{8}{3}$ (this can be justified using the moment inequality given in \citet{Fryz}, see Lemma A2). A sufficient condition for the finiteness of this moment is $\E^{1/h}\left(\vert \xi_0\vert^h\right)\sup_{u\in[0,1]}\sum_{j=1}^p a_j(u)<1$ which is more restrictive than the initial condition given in {\bf A3}. Despite this slight restriction, we only consider the aforementioned sequence of weights in the sequel.   
\item
For model (\ref{tvARCH}), we use the cross-validation method considered in \citet{Fryz}.
The bandwidth parameter is selected by minimizing the function
$$b\mapsto \sum_{t=p+1}^TW_t\left(X_t^2-\mathcal{X}'_t\hat{a}_t^{(-t)}(b)\right)^2,$$
where $\hat{a}_t^{-(t)}(b)=\left(\som_{\substack{k=p+1\\ k\neq t,\ldots,t+p}}^TK\left(\frac{t-k}{Tb}\right)W_k\mathcal{X}^2_k\mathcal{X}'_k\right)^{-1}\som_{\substack{k=p+1\\ k\neq t,\ldots,t+p}}^TK\left(\frac{t-k}{Tb}\right)W_kX_k^2\mathcal{X}_k.$
\item
For model (\ref{lagconst}), we choose the bandwidth $b$ by minimizing
$$(\beta,b)\mapsto \sum_{t=p+1}^TW_t\left( X_t^2-\hat{q}_{1,b,t}^{(-t)}-\left(N_t-\hat{q}_{2,b,t}^{(-t)}\right)'\beta\right)^2.$$
Here for $n=1,2$, $\hat{q}^{(-t)}_{n,b,t}$ is the version
of $\hat{q}_{n,b,t}$ (see (\ref{quotientbase})) defined as $\hat{a}_t^{-(t)}(b)$, $N_t=\left(X_{t-1}^2,\ldots,X_{t-p}^2\right)'$ and $\beta=(a_1,\ldots,a_p)'$. This type of cross-validation is a weighted version of the method proposed by \citet{Hart} for AR models with a time-varying mean. Note that minimizing the latter function with respect to $\beta$, $b$ being fixed, leads to an estimate close to the estimate $\hat{\beta}$ defined in (\ref{MC1}). 
\item
Finally, we discuss the selection of the number of lags $p$ in model \ref{tvARCH}. An information criterion has been studied recently by \citet{ZW} for time-varying regression models. It is possible to adapt the approach used by these authors to our setting. 
To this end, we define for $p=0,1,\ldots,q$,
$$C(p)=\log\left[\sum_{t=p+1}^TW^{(q)}_t\left(X_t^2-\mathcal{X}_i'\hat{a}^{(p)}(t/T)\right)^2\right]+\zeta_T (p+1),$$  
where $\mathcal{X}_t=\left(1,X_{t-1}^2,\ldots,X_{t-p}^2\right)'$, $W^{(q)}_t=\left(\hat{v}+\sum_{j=1}^q X_{t-j}^2\right)^{-2}$, $\hat{a}^{(p)}$ are the coefficients estimates 
obtained for the tv-ARCH with $p$ lags but computed with the weights $W^{(q)}$ and $\zeta_T$ is a vanishing sequence of positive numbers. 
The goal of the selection procedure is to minimize $p\mapsto C(p)$ in the spirit of AIC or BIC criterion used for regression models. 
The bandwidth $b$ is selected by cross validation for the tv-ARCH model with $q$ lags. Of course, condition on the decrease of $\zeta_T$
has to be imposed to get consistency (i.e $\P\left(\arg\min_{0\leq p\leq q}C(p)=p_0\right)\rightarrow 1$ if $p_0\leq q$ is the true number of lags). Inspecting the proofs of Lemma A$6$ and Theorem $3.3$ in \citet{ZW}, we find that condition $T^{2/3}\zeta_T\rightarrow \infty$ guarantees consistency (using the arguments of the proof of Theorem $5$, one can show that the quantity $\phi_T(\phi_T+\rho_T)$ in Lemma A$6$ of \citet{ZW} can be replaced with $T^{-2/3}$ in our context). 
For real data applications we choose $\zeta_T=\log(\log(T))/(Tb)$ where $b$ is selected using cross-validation. 
Our intensive simulation study reported in the supplementary material shows that this choice gives reasonable performances. 
This choice can be also justified using the argument given in \citet{ZW2}: $(p+1)/b$ can be seen as the effective number of parameters     
in kernel smoothing. Hence our choice has a similarity with the Hannan-Quinn information, except that we gave up the constant $2c$ with $c>1$ used in \citet{Han}. We found that adding such a factor underestimates the order in our case. A precise justification of our choice using 
a version of the law of iterated logarithm is not the goal of this paper. However, one can notice that applying cross-validation on an interval of type $\left[c T^{-1/3},C T^{-1/3}\right]$ is compatible with our consistency condition. In the applications, the maximal number of lags is set to $q=10$.
\end{itemize}

Using the approach described above, we conducted an extensive simulation study. Some numerical experiments are reported in the supplementary material and show the good behavior of our method for various simulation setups. This simulation study shows that our estimators and tests have reasonable performances for the sample sizes considered in the sequel. Here we only report the results obtained with real data.  

For real data applications, the bandwidth will be always chosen using CV over a grid of the form $[c,C]\times T^{-1/3}$ where $c$ and $C$ are two positive constants (we recall that $T^{-1/3}$ is the optimal rate of convergence for the bandwidth in the nonparametric estimation of a Lipschitzian regression function). We will also use some acronyms for our models. Model tv$(p)$
denotes the tvARCH model with $p$ lags (the case $p=0$ refers to model (\ref{Star})) while model sptv$(p)$ denotes model (\ref{lagconst}) with $p\geq 1$. In the sequel, we consider two currency exchange rates and one stock market index. The log returns $X_t=\log\left(P_t\right)-\log\left(P_{t-1}\right)$ of the initial series $(P_t)_t$ will be modelized with ARCH processes.

\subsection{Exchange rate USD/Euro}

In this subsection, we study the exchange rate series between the US Dollar and the Euro.
We consider the period from January $03$, $2000$ to February $13$, $2015$. The sample size is $T=3799$.
As usual for this type of series, the autocorrelograms suggest correlation for the squares of the transformed series.   
\begin{figure}[H]
\begin{center}
\includegraphics[width=8cm,height=7cm]{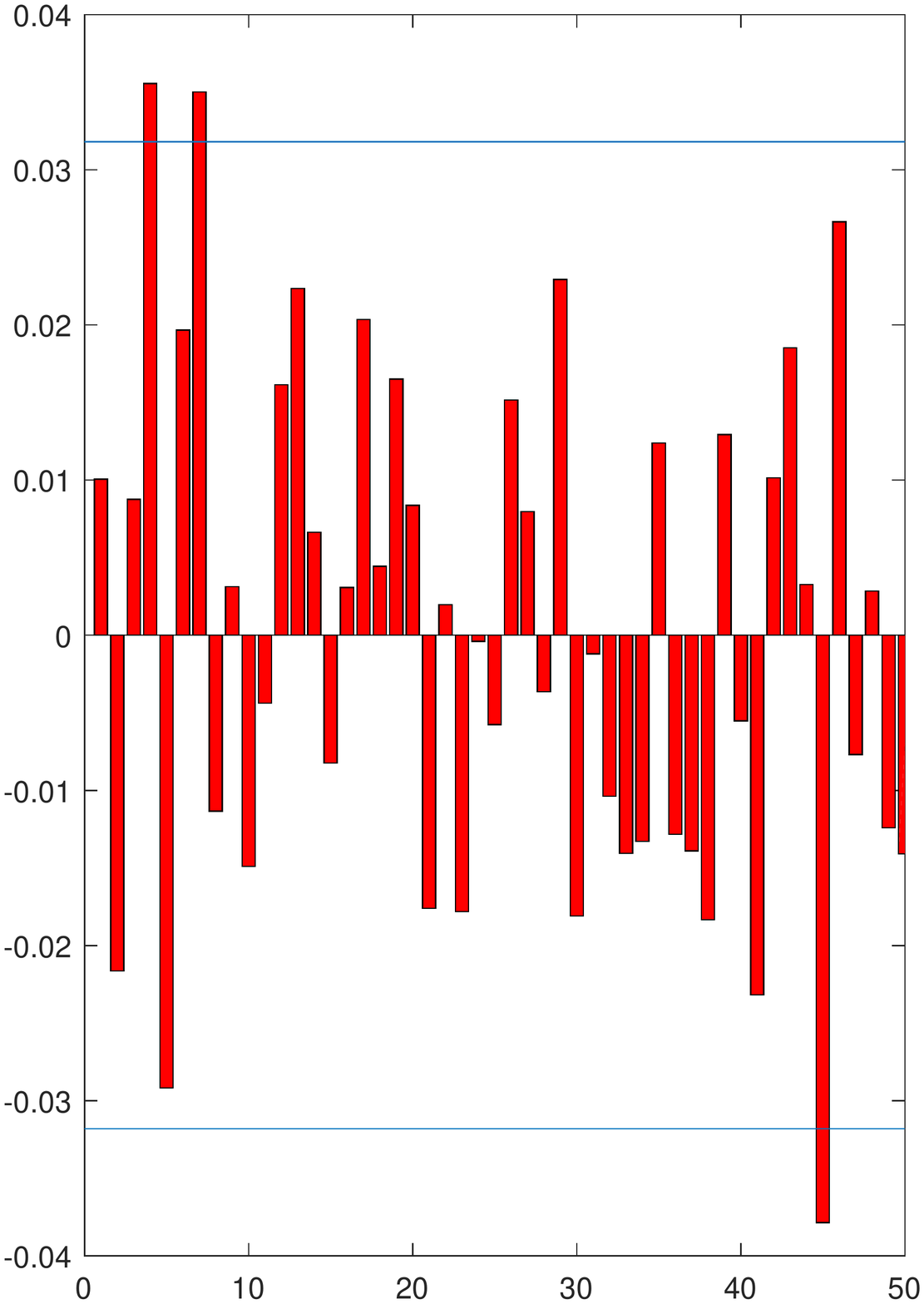}
\includegraphics[width=8cm,height=7cm]{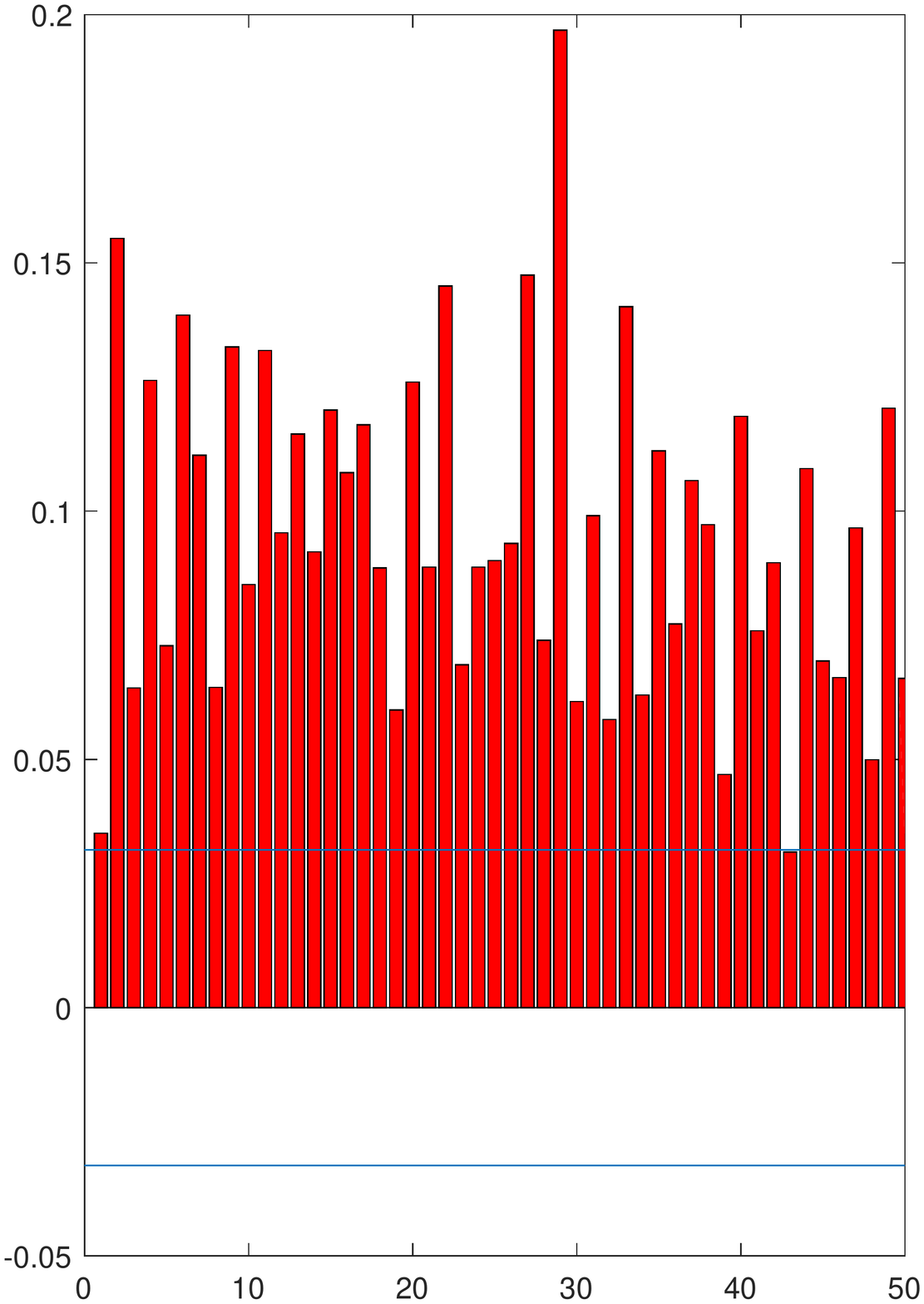}
\end{center}
\caption{Autocorrelogram and autocorrelogram of the squares for the logged and differenced daily exchange rates USD/Euro \label{Hulk}} 
\end{figure}
For this data set, the information criterion selects $p=0$.
To confirm this choice, we apply our procedure with $p=2$. The results for testing the hypothesis of non time-varying coefficients 
are reported in Table \ref{ttt}. The intercept function seems not constant in contrast to the lag coefficients for which it is not possible to reject the null hypothesis.
Fitting a sptv$(2)$ process gives small negative values for the lag coefficients and it is of course not necessary to test if they are equal to zero. This conclusion suggests an absence of second order dynamic for this series. We refer the reader to \citet{Starica} and \citet{Starica2} for other analysis suggesting a similar behavior for some financial time series.

\begin{figure}[H]
\begin{center}
\includegraphics[width=8cm,height=7cm]{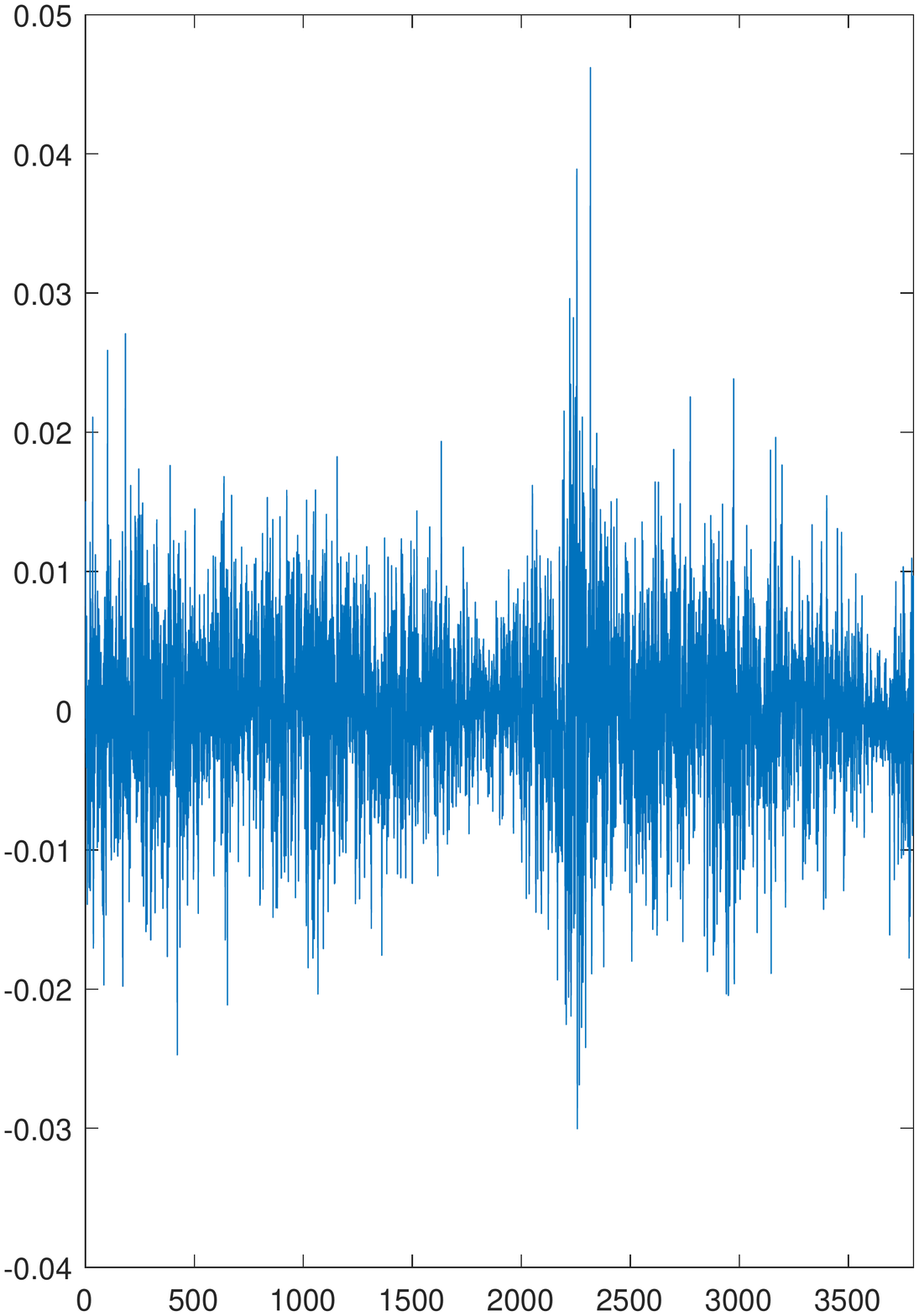}
\includegraphics[width=8cm,height=7cm]{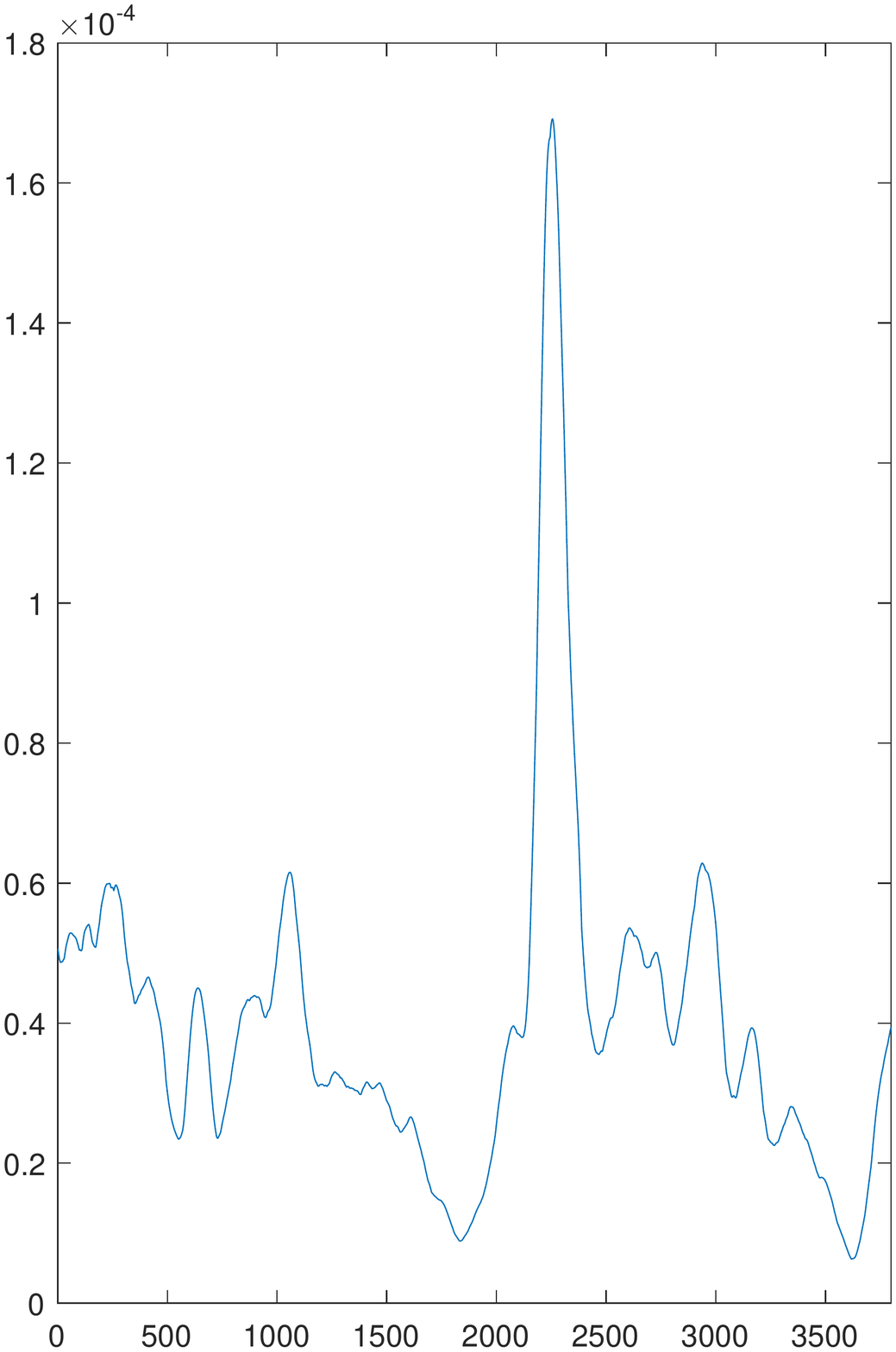}
\end{center}
\caption{Logged and differenced daily exchange rates between USD and Euro and estimation of the unconditional variance \label{euroro}} 
\end{figure}

\begin{table}[H]
\begin{center}
\begin{tabular}{|c|c|c|c|c|}
\hline
Non t-v $a_0$&Non t-v $a_1$& Non t-v $a_2$& Non t-v $(a_1,a_2)$& $\hat{b}$\\\hline
 $0.0005$& $0.138$&$0.1645$ &$0.6415$&$0.028$\\\hline
\end{tabular}
\caption{The $p-$values for testing the hypothesis of non time-varying coefficients (the first line gives the null hypothesis)\label{ttt}} 
\end{center}
\end{table}
A plot of the series with the final estimate of the intercept function is given Figure \ref{euroro}.

\subsection{A second example: the exchange rates between the US Dollar and the Indian Rupee}
In this subsection, we analyze the exchange rates series between the US Dollar and the Indian Rupee over the period starting from 
December $19$, $2005$ to February $18$, $2015$. The sample size is $T=2303$. 
The information criterion selects $p=1$ lag for this series.
From the $p-$values reported in Table \ref{rupipi}, the hypothesis of a constant intercept function is clearly rejected but it is not possible to reject the assumption of a non time-varying first lag coefficient. Fitting a sp$(1)$ process gives a small but significant
lag estimates (the $p-$ value for testing the second order dynamic is less than $10^{-4}$). In contrast, fitting a stationary ARCH process with one lag (one can simply use $b=1$ and our procedure) leads to $\hat{a}_1=0.3041$ (s.e $0.0717$) and several significant lag estimates are found for larger values of $p$ .

\begin{table}[H]
\begin{center}
\begin{tabular}{|c|c|c|c|c|}
\hline
Non t-v $a_0$&Non t-v $a_1$& $\hat{b}_{tv}$& $\hat{a}_1$& $\hat{b}_{sptv}$\\\hline
 $<10^{-4}$ & $0.4215$&$0.035$ &$0.1527$ (s.e $0.0688$)&$0.028$  \\\hline 
\end{tabular}
\end{center}
\caption{Test and estimation for the USD/Rupee series ($p-$values for the test, the selected bandwidths for fitting a tv$(1)$/sptv$(1)$ process and estimation of the first lag coefficients) \label{rupipi}}  
\end{table}

\begin{figure}[H]
\begin{center}
\includegraphics[width=8cm,height=7cm]{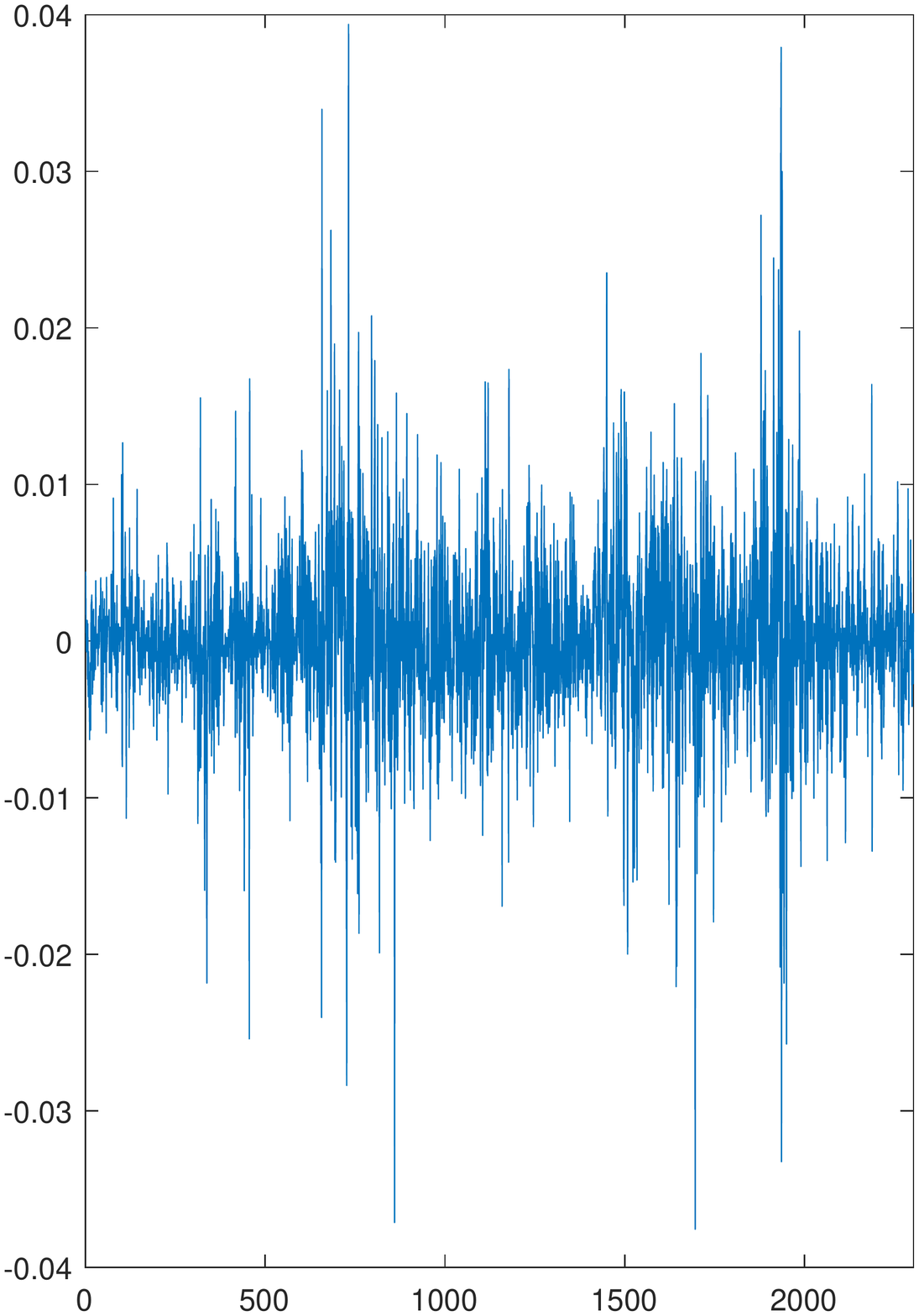}
\includegraphics[width=8cm,height=7cm]{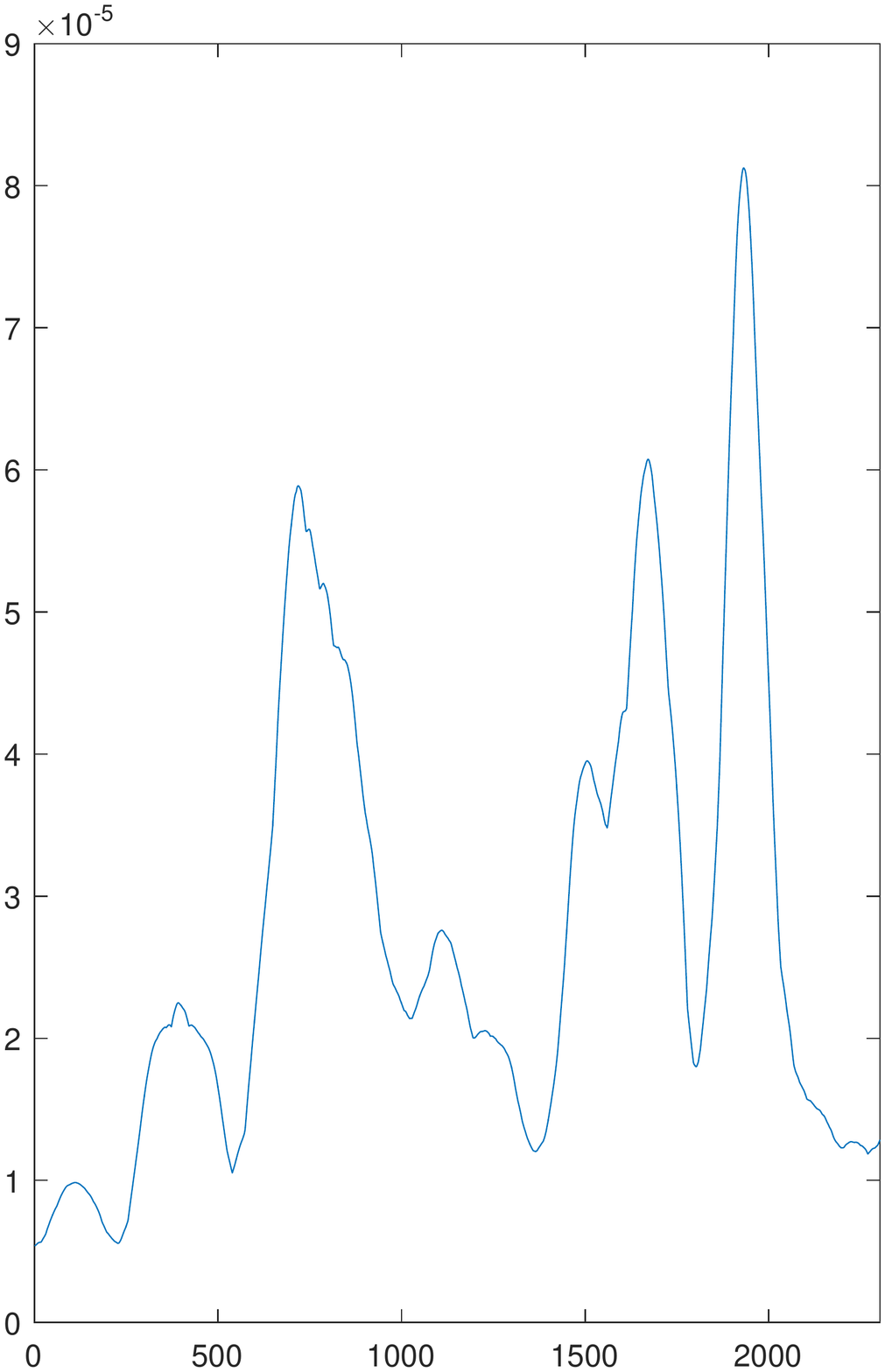}
\end{center}
\caption{Logged and differenced daily exchange rates between the US Dollar and the Indian Rupee and estimation of the intercept function}
\end{figure}

\subsection{A classical stock market index: the FTSE}
\begin{figure}[H]
\centering
\includegraphics[width=9cm,height=7cm]{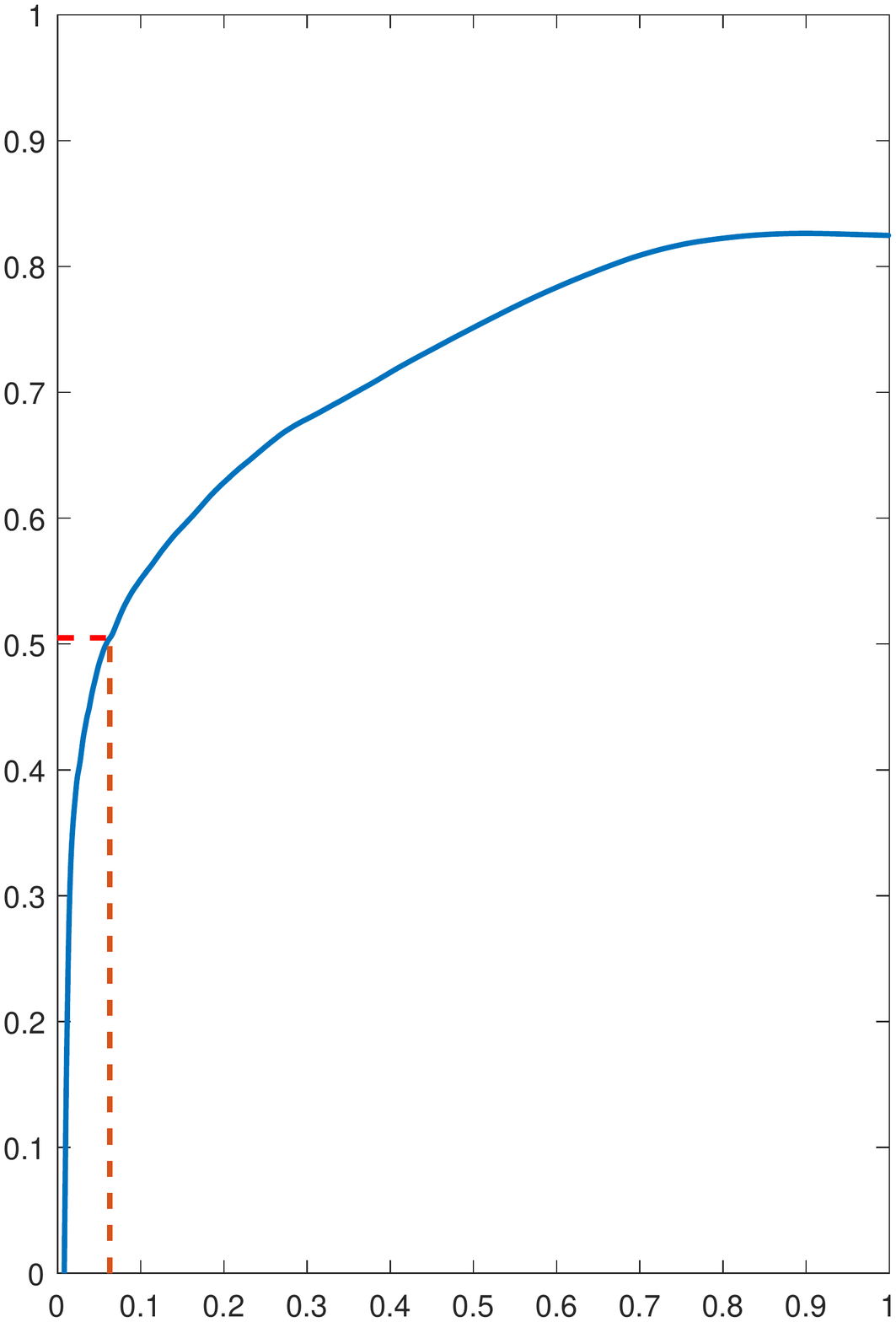}
\caption{Sum of the first five lag coefficients with respect to the value of the bandwidth $b$ (red dashed lines correspond to the bandwidth selected using our method)\label{FTFTFT}}
\end{figure}

Finally, we consider the closing values of the FTSE index from January $04$, $2005$ to March $04$, $2015$, taking as usual the logged and differenced daily returns. Our information criterion selects $p=5$ lags. Testing constancy of the intercept function gives a $p-$value 
of $3\times 10^{-3}$ and the $p-$value for testing constant lag coefficients is $0.066$. Hence, the assumption of constant lag coefficients 
is rejected at level $\alpha=10\%$ (testing constancy of the third and fourth coefficients gives the p-values $0.084$ and $0.067$ respectively, the other $p-$values exceed $10\%$) and considering a tv-ARCH process for this data set could be interesting. 
We also fit a sptv$(5)-$process. The estimated lag parameters and their standard errors are reported in Table \ref{ftft}. The $p-$value for testing the absence of second order dynamic is close to zero. 
\begin{table}[H]
\begin{center}
\begin{tabular}{|c|c|c|}
\hline
$\hat{a}_1$& $\hat{a}_2$& $\hat{a}_3$\\\hline
 $0.0547$ (s.e $0.0321$) & $0.1155$ (s.e $0.0320$) & $0.1204$ (s.e $0.0311$)\\\hline
$\hat{a}_4$& $\hat{a}_5$ & $\hat{b}_{SP}$\\\hline
$0.0942$ (s.e $0.0367$) &$0.1201$ (s.e $0.0324$) & $0.063$\\\hline
\end{tabular}
\end{center}
\caption{Estimated values of the lag coefficients and selected bandwidth for the FTSE \label{ftft}}
\end{table}
Here selecting the number of lags is important because fitting a sptv$(2)-$process for instance does not give significant lag estimates and the selected bandwidth for $p=2$ is very small.
In \citet{Fryz}, it is suggested that stationnary GARCH models give better forecasts for stock market indices than tvARCH processes and that this result could be explained by a more stationarity behavior of these series with respect to currency exchange rates.
This observation is compatible with our analyze of the FTSE index on this period of time which suggests that adding non linearity has a tendency to take away non stationarity, with larger selected bandwidths. However, Figure \ref{FTFTFT} shows that
incorporating a time-varying unconditional variance significantly reduces the values of ARCH parameters. In Figure \ref{FTFTFT}, two extreme cases are observed. When $b\rightarrow 0$, the sum of lag coefficients becomes arbitrary small whereas the value $b=1$ (which corresponds to the fitting of a stationnary ARCH process) leads to larger lag estimates. Moreover, in fitting a sptv$(5)$ process, the ratio $\sqrt{\hat{a}_0(u)}/\hat{\sigma}_t$ has an average of $0.75$ (s.e $0.14$) which means that the contribution of the time-varying intercept has a strong contribution to volatility.     

\paragraph{Acknowlegements.} We would like to thank two anonymous referees, the Associate Editor and the Joint Editor for their constructive comments which helped to improve the presentation of the results and some methodological aspects. We also thank Valentin Patilea for numerous suggestions about nonparametric kernel estimation and bandwidth selection.

\bibliographystyle{plainnat}
\bibliography{bibARCH}

\newpage
\appendix
\section*{\LARGE{Supplementary material}}

\section{Auxiliary results for the proofs}
In this subsection, we consider a general time-varying ARCH process $(X_t)_{1\leq t\leq T}$ defined by 
$$X_t=\xi_t\sqrt{a_0\left(\frac{t}{T}\right)+\sum_{j=1}^p a_j\left(\frac{t}{T}\right)X_{t-j}^2},\quad p+1\leq t\leq T.$$
The different coefficients $a_j$ can be time-varying or not and we assume that assumption {\bf A1} is satisfied. 
Let us introduce additional notations.
\begin{itemize}
\item
We will always denote by $\Vert\cdot\Vert$ the euclidean norm on $\R^g$ for an arbitrary positive integer $g$.
The corresponding operator norm on $\mathcal{M}_g$, the set of matrices of size $g\times g$ and with real coefficients, will be also denoted by $\Vert\cdot\Vert$.
\item   
If $X$ is an integrable random variable taking values in $\mathcal{M}_g$, we set $\bar{X}=X-\E(X)$.
\item
For a sequence $b=b_T\in (0,1)$ of bandwidths, we recall the notation
$$k_{t,i}(b)=\frac{\frac{1}{Tb}K\left(\frac{t-i}{Tb}\right)}{\frac{1}{Tb}\sum_{j=p+1}^T K\left(\frac{t-j}{Tb}\right)},\quad p+1\leq i,t\leq T.$$
Note that $\displaystyle\max_{p+1\leq i,t\leq T}k_{t,i}(b)=O\left((Tb)^{-1}\right)$. This bound will be extensively used in the sequel.
\item
For $1\leq t\leq T$, we set $\mathcal{F}_t=\sigma\left(\xi_s:s\leq t\right).$
\item
Finally, we set $Z_t=\xi_t^2-1$ for $t\in\N$. Then $Z_t$ is centered.
\item
{\bf Important notations}: for simplicity of notations, the quantities $\hat{s}_{j,b,t}$ appearing in the statements of Theorems $1$ and $3$ will be simply denoted by $S_{j,b,t}$ for $j=1,2,3$.
\end{itemize}

We first give a lemma about the regularity of tv-ARCH processes. The following result is crucial for deriving asymptotic properties of our estimators and it is a direct consequence of Theorem $1$ in \citet{DR} (see also \citet{SR}, Theorem $2.1$ and the discussion in Section $5.2$).
\begin{lem}\label{approxstat}
\begin{enumerate}
\item
There exists a constant $C>0$ such that for all $(u,v,T)\in[0,1]^2\times \N^{*}$, 
$$\E\vert X_1(u)^2-X_1^2(v)\vert \leq C\vert u-v\vert,\quad \max_{1\leq t\leq T}\E\vert X_t^2-X^2_t\left(\frac{t}{T}\right)\vert\leq \frac{C}{T}.$$
\end{enumerate}
\end{lem}
From this lemma, we get $\sup_{T\geq p+1}\max_{p+1\leq t\leq T}\E X_t^2<+\infty$.

In the sequel, we will use the following terminology.
\begin{Def}
We will say that a sequence of functions $f_T:\{1,\ldots,T\}\times \R^p_+\rightarrow \R$, $T\geq p+1$, is in the class $\mathcal{L}$ if there exists two positive real numbers $M$ and $L$, not depending on $T$, such that
$$f_T \mbox{ is bounded by } M,\quad \max_{1\leq t\leq T}\vert f_T(t,x)-f_T(t,y)\vert\leq L\Vert x-y\Vert.$$
\end{Def}
 
Now, we will consider two particular classes of processes.
\begin{Def}
\begin{enumerate}
\item
A process $\left(Y_t\right)_{p+1\leq t\leq T}$ is said to be of type I if 
$$Y_t=f_T\left(t,X^2_{t-1},\ldots,X^2_{t-p}\right),\quad p+1\leq t\leq T,$$
where $(f_T)_{T\geq p+1}$ is in the class $\mathcal{L}$.
\item
A process $\left(Y_t\right)_{p+1\leq t\leq T}$ is said to be of type II if 
$$Y_t=f_T\left(t,X^2_{t-1},\ldots,X^2_{t-p}\right)+Z_tg_T\left(t,X_{t-1}^2,\ldots,X_{t-p}^2\right),\quad p+1\leq t\leq T,$$
where $(f_T)_{T\geq p+1}$ and $(g_T)_{T\geq p+1}$ are both in the class $\mathcal{L}$ and $g_T\neq 0$.
\item
A process $\bar{S}$ defined by 
$$\bar{S}_t=\sum_{i=p+1}^Tk_{t,i}(b)\bar{Y}_i,\quad p+1\leq t\leq T,$$
with $Y$ of type I or II will be called a smoothing .
\end{enumerate}
\end{Def}
\paragraph{Notes} 
\begin{enumerate}
\item
An important example of processes of type II is $Y_t=W_tX_{t-j}^2X_t^2$,  for  $j\in\llbracket 1,p\rrbracket$.
Here $W_t$ is given by equation ($2.3$) in the paper. This is due to the decomposition $Y_t=W_t\sigma_t^2X_{t-j}^2+Z_tW_t\sigma_t^2X_{t-j}^2$ and to the particular form of the weights $W_t$ and of $\sigma_t^2$.  
\item
Some smoothings appear in the expression of $\hat{q}_{j,b,t}$ for $j=1,2$. Our method for proving Theorem $1$ is to make an asymptotic expansion of the estimator $\hat{\beta}$ and to show that the effect of the smoothings incorporated in $\hat{\beta}$ is negligible by computing some moments. The terminology type I or type II is just used for identifying the number of smoothings which impose a moment restriction.  
\end{enumerate}

\subsection{Covariance inequalities}
Here, we assume that the assumptions of Theorem $1$ are fulfilled. Sometimes, assumption {\bf A2(h)} can be used for a general value of the integer $h$, this will be precised in the statements of our Lemma/Propositions.

\begin{lem}\label{coupling}
Let $s$ and $t$ be two natural integers such that $T\geq t\geq s+1\geq p+1$. 
Now let $\left(\Upsilon_s,\Upsilon_{s-1},\ldots,\Upsilon_{s-p+1}\right)$ be a random vector independent from
the sequence $(\xi_t)_{1\leq t\leq T}$ and with the same distribution as\\ $\left(X_s,X_{s-1},\ldots,X_{s-p+1}\right)$.
For $s+1\leq k\leq t$, we define recursively 
$$\Upsilon_k=\xi_k\sqrt{a_0\left(\frac{t}{T}\right)+\sum_{j=1}^p a_j\left(\frac{t}{T}\right)\Upsilon_{k-j}^2}.$$
Then 
$$\E\vert X_t^2-\Upsilon_t^2\vert \leq 2d c^{\frac{t-s-1}{p}},$$
where $d=\sup_{T\geq p+1}\max_{1\leq t\leq T} \E X_t^2$ and $c$ is defined in assumption ${\bf A1}$.
\end{lem}

\paragraph{Proof of Lemma \ref{coupling}}
\begin{itemize}
\item
Assume first that $s+1\leq t\leq s+p$. Then
\begin{eqnarray*}
\E\vert X_t^2-\Upsilon_t^2\vert&=& \sum_{j=1}^p a_j\left(\frac{t}{T}\right)\E\vert X_{t-j}^2-\Upsilon_{t-j}^2\vert\\
&\leq& c \max_{1\leq j\leq p}\E\vert X_{t-j}^2-\Upsilon_{t-j}^2\vert\\
&\leq& 2dc. 
\end{eqnarray*}
Since $\frac{t-s-1}{p}\leq 1$, the result follows in this case.
\item
Suppose the inequality true for any $t\in\llbracket s+p, N\rrbracket$, where $N\in\llbracket s+p,T-1\rrbracket$. Then 
\begin{eqnarray*}
\E\vert X_{N+1}^2-\Upsilon_{N+1}^2\vert &\leq& c\times \max_{1\leq j\leq p}\E\vert X_{N+1-j}^2-\Upsilon_{N+1-j}^2\vert\\
&\leq& 2d \max_{1\leq j\leq p}c^{1+\frac{N+1-j-s-1}{p}}\\
&\leq& 2d c^{\frac{N-s}{p}}.
\end{eqnarray*}
\end{itemize}
Then the result of the lemma follows from a finite induction.$\square$

\begin{lem}\label{Esp}
Let $h,s,t$ be three integers such that $p\leq s\leq t\leq T$ and $h\geq 1$. Assume that $\E\xi_0^{2h(1+\delta)}<\infty$ with $0<\delta<1$.  Let $U_s$ be an integrable random
variable $\mathcal{F}_s-$measurable and $G:\R^{p+k+1}\rightarrow \R$ a bounded and Lipschitzian function. We set
$$U_t=Z_{t+\ell_1}\cdots Z_{t+\ell_o}G\left(X^2_{t-p},X^2_{t-p+1},\ldots,X^2_{t+k}\right),$$
with $0\leq \ell_1,\ldots,\ell_o\leq k$ and $o\leq h$. Assume that $\E\vert U_s Z_{t+\ell_1}\cdots Z_{t+\ell_o}\vert^{1+\delta}<\infty$.  Then, we have
$$\vert\mbox{Cov}\left(U_s,U_t\right)\vert\leq \left(C_1\vee C_2\right) c^{\frac{\kappa(t-s)}{p}},$$
where, setting $\kappa=\frac{\delta}{1+\delta}$, 
$$C_1= M(G)c^{-\kappa}\left(\E\vert U_sZ_{t+\ell_1}\cdots Z_{t+\ell_o}\vert+\E\vert U_s\vert\cdot \E\vert Z_{t+\ell_1}\cdots Z_{t+\ell_o}\vert\right).$$
and 
$$C_2=\frac{d^{\kappa}c^{-\kappa \frac{p+1}{p}}}{1-c^{\frac{\kappa}{p}}}L(G)^{\kappa} M(G)^{1-\kappa} \E^{\frac{1}{1+\delta}} \vert U_sZ_{t+\ell_1}\cdots Z_{t+\ell_o}\vert^{1+\delta},$$
where $d$ is given in Lemma \ref{coupling} and $M(G)$ (resp. $L(G)$) denotes the supremum (resp. the Lipschitz constant) of the function $G$.
\end{lem} 

\paragraph{Proof of Lemma \ref{Esp}.}
\begin{itemize}
\item
Assume first that $t-p\leq s$. Then, it is easy to get the bounds
\begin{eqnarray*}
\vert\mbox{Cov}\left(U_s,U_t\right)\vert &\leq& M(G)\left(\E\vert U_sZ_{t+\ell_1}\cdots Z_{t+\ell_o}\vert+\E\vert U_s\vert\cdot \E\vert Z_{t+\ell_1}\cdots Z_{t+\ell_o}\vert\right)\\
&\leq& C_1c^{\kappa\frac{t-s}{p}}.
\end{eqnarray*}
\item
Now, assume that $t\geq s+1+p$. 
We have the equality $\mbox{Cov}\left(U_s,U_t\right)=\E\left(U_s\left(U_t-U'_t\right)\right)$, where
$$U'_t=Z_{t+\ell_1}\cdots Z_{t+\ell_o}G\left(\Upsilon^2_{t-p},\ldots, \Upsilon^2_{t+k}\right)$$
is a random variable independent from $\mathcal{F}_s$.
Here $\Upsilon$ denotes the process introduced in Lemma \ref{coupling}.
Using Lemma \ref{coupling}, the following bounds are valid.
\begin{eqnarray*}
&&\vert \mbox{Cov}\left(U_s,U_t\right)\vert\\
&\leq& \E\vert U_sZ_{t+\ell_1}\cdots Z_{t+\ell_o}\cdot\vert G\left(X^2_{t-p},\ldots, X^2_{t+k}\right)-G\left(\Upsilon^2_{t-p},\ldots, \Upsilon_{t+k}^2\right)\vert\\
&\leq&  \left(2 M(G)\right)^{1-\kappa}\E\vert U_s Z_{t+\ell_1}\cdots Z_{t+\ell_o}\vert\cdot\vert G\left(X^2_{t-p},\ldots, X^2_{t+k}\right)-G\left(\Upsilon^2_{t-p},\ldots, \Upsilon_{t+k}^2\right)\vert^{\kappa}\\
&\leq&  L(G)^{\kappa}\left(2 M(G)\right)^{1-\kappa} \sum_{i=-p}^k \E \vert U_s Z_{t+\ell_1}\cdots Z_{t+\ell_o}\vert\cdot\vert X^2_{t+i}-\Upsilon_{t+i}^2\vert^{\kappa}\\
&\leq& L(G)^{\kappa}\left(2 M(G)\right)^{1-\kappa} \sum_{i=-p}^k \E^{\frac{1}{1+\delta}} \vert U_sZ_{t+\ell_1}\cdots Z_{t+\ell_o}\vert^{1+\delta}\times \E^{\kappa}\vert X^2_{t+i}-\Upsilon_{t+i}^2\vert\\
&\leq& (2d)^{\kappa} L(G)^{\kappa}\left(2 M(G)\right)^{1-\kappa} \E^{\frac{1}{1+\delta}} \vert U_sZ_{t+\ell_1}\cdots Z_{t+\ell_o}\vert^{1+\delta}\times \sum_{i=-p}^k c^{\kappa\frac{t+i-s-1}{p}}\\
&\leq& \frac{d^{\kappa}c^{-\kappa \frac{p+1}{p}}}{1-c^{\frac{\kappa}{p}}}L(G)^{\kappa} M(G)^{1-\kappa} \E^{\frac{1}{1+\delta}} \vert U_sZ_{t+\ell_1}\cdots Z_{t+\ell_o}\vert^{1+\delta} c^{\kappa\frac{t-s}{p}}.\\
\end{eqnarray*}
\end{itemize}
Then the result announced in Proposition \ref{Esp} is a direct consequence of the two previous points.$\square$
\bigskip

The following corollary is a direct consequence of Lemma \ref{Esp}.

\begin{cor}\label{incov}
Let $h\geq 1$ be an integer. Assume that $\E\xi_0^{2h(1+\delta)}<\infty$ with $0<\delta<1$.
Let $s,t,q,o$ four non-negative integers such that $p\leq s\leq t$ and $q+o\leq h$. Let $U_s$ and $U_t$ two random variables defined by 
$$U_s=Z_{h_1}\cdots Z_{h_q}H\left(X_1^2,\ldots,X^2_s\right),\quad U_t=Z_{t+\ell_1}\cdots Z_{t+\ell_o}G\left(X^2_{t-p},\ldots,X^2_{t+k}\right),$$
where $H$ and $G$ are two elements of $\mathcal{L}$ and $1\leq h_1\leq\cdots\leq h_q\leq s$ and $0\leq\ell_1\leq \cdots\leq \ell_o\leq k$.
 Then, we have the bound
$$\mbox{Cov}(U_s,U_t)\leq \left(C_1\vee C_2\right) c^{\kappa\frac{t-s}{p}},$$
where $C_1=M(G)M(H)c^{-\kappa}\left(\E\vert Z_1\vert^{q+o} +\E\vert Z_1\vert^o\cdot\E\vert Z_1\vert^q\right)$ and 
$$C_2=\frac{d^{\kappa}c^{-\kappa \frac{p+1}{p}}}{1-c^{\frac{\kappa}{p}}}L(G)^{\kappa} M(G)^{1-\kappa}M(H) \E^{\frac{1}{1+\delta}}\left(\vert Z_1\vert^{(q+o)(1+\delta)}\right).$$
\end{cor}

\subsection{Moment bounds}

 The two next results are crucial for proving the asymptotic normality of our estimators. In particular, Proposition \ref{cumulantcontrol} gives conditions under which some partial sums involving smoothings are convergent to zero with a faster rate than $\sqrt{T}$.

\begin{lem}\label{inter}
Assume that $\E\xi_0^{2h(1+\delta)}<\infty$ for a positive integer $h$. 
Let $Y^{(1)},\ldots,Y^{(q)}$ be $q\geq 1$ process of type I or II with at most $h$ processes of type II.
Then for a family $\left\{z_{T,i}: p+1\leq i\leq T, T\geq p+1\right\}$ of deterministic positive weights, we have
$$\sum_{p+1\leq i_1,\ldots,i_q\leq T}z_{T,i_1}\cdots z_{T,i_q}\vert \E\left(\bar{Y}^{(1)}_{i_1}\cdots \bar{Y}^{(q)}_{i_q}\right)\vert=O\left((\phi_Ts_T)^{\frac{q}{2}}\right),$$
where $s_T=\sum_{i=p+1}^T z_{T,i}$ and $\phi_T=\max_{p+1\leq i\leq T}z_{T,i}$.
\end{lem}

\paragraph{Proof of Lemma \ref{inter}} 
We set $\beta=c^{\frac{\delta}{p(1+\delta)}}$.
The result is clear for $q=1$. Assume that $q\geq 2$. 
First, we observe that from Corollary \ref{incov}, we have for $p+1\leq t_1\leq\cdots\leq t_q\leq T$ and $1\leq j\leq q-1$,
\begin{equation}\label{cocov}
\vert \mbox{Cov}\left(\bar{Y}^{(1)}_{t_1}\cdots \bar{Y}^{(j)}_{t_j},\bar{Y}^{(j+1)}_{t_{j+1}}\cdots \bar{Y}^{(q)}_{t_q}\right)\vert \leq C \beta^{t_{j+1}-t_j},
\end{equation}
where $C>0$ does not depend on $T$ and on $t_1,\ldots,t_q$. 
Inequality (\ref{cocov}) follows from the fact that the covariance given in (\ref{cocov}) can be decomposed as a sum of covariances of the
form $\mbox{Cov}\left(U_{t_j},U_{t_{j+1}}\right)$ given in Corollary \ref{incov} (replacing $s$ and $t$ with $t_j$ and $t_{j+1}$ respectively). Inequality (\ref{cocov}) is crucial for the sequel.

We set for $1\leq j\leq q$,
$$A^{(j)}_T\left(\bar{Y}^{(1)},\ldots,\bar{Y}^{(j)}\right)=
\sum_{p+1\leq i_1\leq \cdots\leq i_j\leq T}z_{T,i_1}\cdots z_{T,i_j}\vert \E\left(\bar{Y}_{i_1}^{(1)}\cdots\bar{Y}_{i_j}^{(j)}\right)\vert.$$
We use a classical method for bounding sums of cross moments using bounds on covariances (see \citet{D} p. $78$).
For a $q-$uplet ${\bf i}=(i_1,\ldots,i_q)\in \llbracket p+1,T\rrbracket$ such that $i_1\leq \cdots\leq i_q$, we define
$$s({\bf i})=\min\left\{j\leq q: i_{j+1}-i_j=\max_{1\leq \ell\leq q-1}(i_{\ell+1}-i_{\ell})\right\}.$$
Then, using the bound (\ref{cocov}), we have
\begin{eqnarray*}
&&A^{(q)}_T\left(\bar{Y}^{(1)},\ldots, \bar{Y}^{(q)}\right)\\
&\leq & \sum_{j=1}^{q-1}A^{(j)}_T\left(\bar{Y}^{(1)},\ldots, \bar{Y}^{(j)}\right)\cdot A^{(q-j)}_T\left(\bar{Y}^{(j+1)},\ldots, \bar{Y}^{(q)}\right)\\
&+&C\sum_{j=1}^{q-1}\sum_{r=0}^{T-1}\beta^r\sum_{{\bf i}: s({\bf i})=j, i_{j+1}=i_j+r}z_{T,i_1}\cdots z_{T,i_q}.
\end{eqnarray*}
Since, 
$$\sum_{{\bf i}: s({\bf i})=j, i_{j+1}=i_j+r}z_{T,i_1}\cdots z_{T,i_q}\leq s_T\left(\phi_T\right)^{q-1}(r+1)^q,$$
we conclude that
$$A^{(q)}_T\left(\bar{Y}^{(1)},\ldots, \bar{Y}^{(q)}\right)\leq \sum_{j=1}^{q-1}A^{(j)}_T\left(\bar{Y}^{(1)},\ldots, \bar{Y}^{(j)}\right)\cdot A^{(q-j)}_T\left(\bar{Y}^{(j+1)},\ldots, \bar{Y}^{(q)}\right)+O\left(s_T(\phi_T)^{q-1}\right).$$
Since $\phi_T\leq s_T$, we have $s_T(\phi_T)^{q-1}\leq (s_T\phi_T)^{\frac{q}{2}}$. Then using an induction on $q$, it easy to prove that 
$$A^{(q)}_T\left(\bar{Y}^{(1)},\ldots, \bar{Y}^{(q)}\right)=O\left((s_T\phi_T)^{\frac{q}{2}}\right).$$
This proves Lemma $2$.$\square$

\begin{Prop}\label{cumulantcontrol}
Assume that $\E \xi_0^{4h(1+\delta)}<\infty$ for a positive integer $h$ and that $b\sqrt{T}\rightarrow\infty$.
Let $Y^{(1)},\ldots,Y^{(q)}$ be $q\geq 2$ processes of type I or II with at most $h$ processes of type II. 
We denote by $\bar{S}^{(1)},\ldots,\bar{S}^{(q)}$ the corresponding smoothings.
\begin{enumerate}
\item
If $\left\{\mu_{t,T}: p+1\leq t\leq T, T\geq p+1\right\}$ denotes a family of real numbers such that $$\sup_{T\geq p+1}\max_{p+1\leq t\leq T}\vert\mu_{t,T}\vert<\infty,$$ we have, using the notations of point $2$, 
$$\frac{1}{\sqrt{T}}\sum_{t=p+1}^T \mu_{t,T}\bar{S}^{(1)}_t\cdots \bar{S}^{(q)}_t=o_{\P}(1).$$
\item
We have also $\frac{1}{\sqrt{T}}\sum_{t=p+1}^T \bar{Y}^{(1)}_t\bar{S}^{(2)}_t\cdots \bar{S}^{(q)}_t=o_{\P}(1)$.
\item
If $Y$ is a process of type $I$ (resp. $II$), we have for all positive integer $h'$ (resp. $h=h'$), $\max_{1\leq t\leq T}\E \bar{S}_t^{2h'}=O\left((Tb)^{-h'}\right)$.
\item
Assume that $Y^{(1)}$ is a process of type $I$. Then $\frac{1}{T}\sum_{t=p+1}^T \overline{Y}^{(1)}_t$ converges to $0$ a.s.
\end{enumerate}
\end{Prop}

\paragraph{Proof of Proposition \ref{cumulantcontrol}} 
\begin{enumerate}
\item
Assume that $\sup_{t,T}\vert\mu_{t,T}\vert\leq C$.
Taking the second order moment, we get 
\begin{eqnarray*}
&&\E\vert \frac{1}{\sqrt{T}}\sum_{t=p+1}^T\mu_{t,T} \bar{S}^{(1)}_t\cdots\bar{S}^{(q)}_t\vert^2\\
&\leq& \frac{C^2}{T}\sum_{s,t=p+1}^T\sum_{i_1,M_1,\ldots,i_q,j_q=p+1}^Tk_{s,i_1}k_{t,M_1}\cdots k_{s,i_q}k_{t,j_q}\cdot\\
&& \vert\E\left(\bar{Y}^{(1)}_{i_1}\bar{Y}^{(1)}_{M_1}\cdots \bar{Y}^{(q)}_{i_q}\bar{Y}^{(q)}_{j_q}\right)\vert.
\end{eqnarray*}

Using Lemma $2$ (replacing $h$ with $2h$) with $s_T=1$ and $\phi_T=O\left(\frac{1}{Tb}\right)$, the last bound is $O\left(\frac{1}{T^{q-1}b^q}\right)$. The result follows using the bandwidth assumptions.

\item
The second order moment writes
\begin{eqnarray*}
&&\frac{1}{T}\sum_{s,t=p+1}^T\sum_{i_1,M_1,\ldots,i_{q-1},j_{q-1}=1}^Tk_{s,i_1}k_{t,j_1}\cdots k_{s,i_{q-1}}k_{t,j_{q-1}}\cdot\\
&&\times\E\left(\bar{Y}^{(1)}_s \bar{Y}^{(1)}_t \bar{Y}^{(2)}_{i_1}\bar{Y}^{(2)}_{j_1}\cdots \bar{Y}^{(q)}_{i_{q-1}}\bar{Y}^{(q)}_{j_{q-1}}\right).
\end{eqnarray*}
Using the bound $\max_{p+1\leq t,i\leq T}k_{t,i}=O\left(\frac{1}{Tb}\right)$ and applying Lemma $2$ with $z_{T,i}=1$, it is easy to show that this second order 
moment is $O\left(\frac{1}{T^{q-1}b^{2(q-1)}}\right)$. This leads to the result.$\square$

\item
This is a consequence of Lemma $2$, using the inequality 
$$\E\left(\bar{S}_t^{2h'}\right)\leq \som_{p+1\leq i_1,\ldots i_{2h'}\leq T}k_{t,i_1}\cdots k_{t,i_{2h'}}\left\vert \E\left(\bar{Y}_{i_1}\cdots \bar{Y}_{i_{2h'}}\right)\right\vert.$$
\item
We have for $\epsilon>0$, 
\begin{eqnarray*}
&&\P\left(\frac{1}{T}\vert \sum_{t=p+1}^T\bar{Y}^{(1)}_t\vert>\epsilon\right)\\
&\leq& \frac{1}{\epsilon^4T^4}\E\vert \sum_{t=p+1}^T\bar{Y}_t^{(1)}\vert^4\\
&\leq&  \frac{1}{\epsilon^4T^4}\sum_{p+1\leq i_1,\ldots i_4\leq T}\vert \E\left(\bar{Y}^{(1)}_{i_1}\cdots \bar{Y}^{(1)}_{i_4}\right)\vert
\end{eqnarray*}
Using Lemma $2$, the last bound is $O\left(\frac{1}{T^2}\right)$. Then the result follows from the Borel-Cantelli Lemma. $\square$
\end{enumerate}

\subsection{Control of deterministic quantities using local stationarity}

\begin{lem}\label{biais}
\begin{enumerate}
\item
For $u\in [0,1]$, we set $s_2(u)=\E\left(W_1(u)M_1(u)N_1(u)'\right)$ and $s_3(u)=\E\left(W_1(u)M_1(u)M_1(u)'\right)$.
Then we have
\begin{equation}\label{bor1}
\max_{p+1\leq t\leq T}\left\{ \Vert s_{3,t}-s_3\left(\frac{t}{T}\right)\Vert +\Vert s_{2,t}-s_2\left(\frac{t}{T}\right)\Vert\right\}=O\left(\frac{1}{T}\right),
\end{equation}
\begin{equation}\label{bor2}
\inf_{u\in [0,1]}\det\left(s_3(u)\right)>0,
\end{equation}
\begin{equation}\label{bor3}
\sup_{u\in[0,1]}\Vert q_2(u)\Vert<\infty,
\end{equation}
\begin{equation}\label{bor4}
\max_{p+1\leq t\leq T}\Vert q_{2,t}-q_2\left(\frac{t}{T}\right)\Vert=O\left(\frac{1}{T}\right).
\end{equation}
\item
We have 
$$\sup_{T\geq p+1}\max_{p+1\leq t\leq T}\left\{\Vert \E^{-1}(S_{3,t})\Vert +\Vert \E\left(S_{1,t}\right)\Vert +\Vert \E\left(S_{2,t}\right)\Vert\right\}<\infty.$$
\item
Setting $\eta_{j,t}=\E^{-1}(S_{3,t})\E\left(S_{j,t}\right)-q_{j,t}$ for $j=1,2$, we have
$$\max_{p+1\leq t\leq T}\left\{\Vert \eta_{1,t}\Vert +\Vert \eta_{2,t}\Vert\right\}=O(b).$$
\item
We have 
\begin{equation}\label{borbor}
\max_{p+1\leq t\leq T}\left\{\Vert S_{3,t}\Vert +\Vert S_{3,t}^{-1}\Vert\right\}=O_{\P}(1).
\end{equation}
\end{enumerate}
\end{lem}

\paragraph{Proof of Lemma \ref{biais}}
\begin{enumerate}
\item
We prove the four assertions successively.
\begin{itemize}
\item
Since $W_tM_tM_t'=f\left(t/T,X_{t-1}^2,\ldots,X_{t-p}^2\right)$ where $f:[0,1]\times \R_+^p\rightarrow\mathcal{M}_{m,m}$ satisfies 
$$\Vert f(u,x_1,\ldots,x_p)-f(u,y_1,\ldots,y_p)\Vert \leq C\sum_{i=1}^p\vert x_i-y_i\vert$$
for some positive constant $C$, 
Lemma $1$ given in the paper yields to $\max_{p+1\leq t\leq T}\Vert s_{3,t}-s_3\left(\frac{t}{T}\right)\Vert=O(1/T)$. 
The conclusion for $s_{2,t}$ follows in the same way. This shows (\ref{bor1}).
\item
Next, we show (\ref{bor2}). 
Let $\lambda(u)$ be the smallest eigenvalue of $\E\left(W_1(u)M_1(u)M_1(u)'\right)$. From Lemma $1$, the application 
$u\mapsto \E\left(W_1(u)M_1(u)M_1(u)'\right)$ is Lipschitz continuous. Moreover, it is easily shown that for all $u\in [0,1]$, the matrix $\E\left(W_1(u)M_1(u)M_1(u)'\right)$ is positive definite. This entails that the application $u\mapsto \lambda(u)$ is continuous and 
positive. This implies (\ref{bor2}).
\item
Since $\sup_{u\in [0,1]}\Vert W_1(u)M_1(u)N_1(u)'\Vert$ is bounded, we deduce from (\ref{bor2}) that $\sup_{u\in[0,1]}\Vert q_2(u)\Vert<\infty$.  
\item
The assertion (\ref{bor4}) easily follows from (\ref{bor1}), (\ref{bor2}) and (\ref{bor3}). 
\end{itemize}
\item
Since 
$$C=\sup_{\omega,T,i}W_i(\omega)\Vert M_i(\omega)\Vert <\infty,$$
we have 
$$\Vert \E\left(S_{1,t}\right)\Vert\leq \sum_{i=p+1}^T k_{t,i}\Vert\E\left(W_iM_iX_i^2\right)\Vert\leq C\sup_{T\geq p+1}\max_{1\leq t\leq T}\E\left(X_t^2\right).$$
The same kind of inequality holds for $\Vert \E\left(S_{2,t}\right)\Vert$. \\
It remains to prove that 
\begin{equation}\label{borne}
\sup_{T\geq p+1}\max_{p+1\leq t\leq T}\Vert \E^{-1}(S_{3,t})\Vert<\infty.
\end{equation}
If $x\in\R^j$, with $\Vert x\Vert=1$, we have using (\ref{bor1}), 
$$x'\E(S_{3,t})x\geq \inf_{p+1\leq t\leq T}x's_{3,t}x\geq \inf_{u\in[0,1]}x's_3(u)x-\frac{C}{T},$$
for a suitable constant $C>0$. Then, using (\ref{bor2}), there exists $\lambda>0$ such that $x'\E(S_{3,t})x\geq \lambda-\frac{C}{T}$. We deduce that if $T$ is large enough, 
the smallest eigenvalue of $\E(S_{3,t})$ is bounded from below. This means that there exists an integer $T_0\geq p+1$ such that 
$$\sup_{T\geq T_0}\max_{p+1\leq t\leq T}\Vert \E^{-1}(S_{3,t})\Vert<\infty.$$
But since each of the matrices $\E(S_{3,t})$ is easily shown to be positive definite for $p+1\leq t\leq T$ and $T\geq T_0$, (\ref{borne}) easily follows.

\item
For a Lipschitzian function $f$ defined over $[0,1]$, the assumptions made on the kernel $K$ implies that 
$$\max_{1\leq t\leq T}\vert f(t/T)-\sum_{i=1}^Tk_{t,i}f(i/T)\vert=O(b).$$
We only prove that $\max_{p+1\leq t\leq T}\Vert \eta_{1,t}\Vert=O(b)$, the proof for $\eta_{2,t}$ is similar.
We use the decomposition
\begin{equation}\label{dec1} 
\eta_{1,t}=s_{3,t}^{-1}\left(s_{3,t}-\E(S_{3,t})\right)\E^{-1}(S_{3,t})s_{1,t}+\E^{-1}(S_{3,t})\left(\E\left(S_{1,t}\right)-s_{1,t}\right).
\end{equation}
From the proof of the two first points of the present Lemma, it is easily seen that 
\begin{equation}\label{dec2}
\sup_{T\geq p+1}\max_{p+1\leq t\leq T}\left\{\Vert \E^{-1}(S_{3,t})\Vert+\Vert s_{1,t}\Vert +\Vert s_{3,t}^{-1}\Vert\right\}<\infty.
\end{equation}
Moreover
$$\max_{p+1\leq t\leq T}\Vert s_{3,t}-\E\left[W_1(t/T)M_1(t/T)M_1'(t/T)\right]\Vert=O(1/T).$$
Since $\E\left(W_tM_tX_t^2\right)=\E\left(W_tM_t\sigma_t^2\right)$, the choice of the weights $W_t$ entails also
$$\max_{p+1\leq t\leq T}\Vert s_{1,t}-\E\left[W_1(t/T)M_1(t/T)X_1^2(t/T)\right]\Vert=O(1/T).$$
Now, since the two applications $$u\mapsto d(u)=\E\left(W_1(u)M_1(u)X_1^2(u)\right)\quad u\mapsto e(u)=\E\left(W_1(u)M_1(u)M_1(u)'\right)$$ are Lipschitz continuous, we get
\begin{equation}\label{dec3}
\max_{p+1\leq t\leq T}\left\{\Vert\E\left(S_{1,t}\right)-s_{1,t}\Vert +\Vert \E(S_{3,t})-s_{3,t}\Vert\right\}=O(b).
\end{equation}
Then, the result announced follows easily from (\ref{dec1}), (\ref{dec2}) and (\ref{dec3}).
\item
We use the decomposition $S_{3,t}=\bar{S}_{3,t}+\E(S_{3,t})$. From the previous points, we have 
$$\max_{p+1\leq t\leq T}\left\{\Vert \E(S_{3,t})\Vert +\Vert \E^{-1}(S_{3,t})\Vert\right\}=O(1).$$
Moreover for $\epsilon>0$, we have using point $3$ of Proposition $3$,
\begin{eqnarray*}
\P\left(\max_{p+1\leq t\leq T}\Vert\bar{S}_{3,t}\Vert>\epsilon\right)&\leq& \frac{1}{\epsilon^4}\sum_{t=p+1}^T\E\Vert \bar{S}_{3,t}\Vert^4\\
&\leq& \frac{C}{\epsilon^4 Tb^2}.
\end{eqnarray*}
Then we conclude that $\max_{p+1\leq t\leq T}\Vert \bar{S}_{3,t}\Vert=o_{\P}(1)$. Then (\ref{borbor}) easily follows.$\square$ 
\end{enumerate}

\section{Proof of Theorem $1$}\label{preuve1}

Setting 
$$\hat{\beta}=D_T^{-1}L_T,\quad L_T=\sum_{t=p+1}^TW_t\hat{O}_t\hat{V}_t,\quad D_T=\sum_{t=p+1}^TW_t\hat{O}_t\hat{O}'_t,$$
we have $\hat{\beta}-\beta=D_T^{-1}\sum_{t=p+1}^TW_t\hat{O}_t\left(\hat{V}_t-\hat{O}_t'\beta\right)$.
Using the two relations
$\hat{V}_t=V_t-M_t'\left(\hat{q}_{1,t}-q_{1,t}\right),\quad \hat{O}_t=O_t-\left(\hat{q}_{2,t}-q_{2,t}\right)'M_t,$
we obtain
$$\widetilde{\beta}-\beta=D_T^{-1}\sum_{t=p+1}^TW_t\left(O_t-(\hat{q}_{2,t}-q_{2,t})'M_t\right)\cdot\left((\xi_t^2-1)\sigma_t^2-M_t'(\hat{q}_{1,t}-q_{1,t})+M_t'(\hat{q}_{2,t}-q_{2,t})\beta\right).$$
We also set $U_t=W_tO_tM_t'$. Observe that $\E U_t=0$.
This yields to the following decomposition.
$$\hat{\beta}-\beta=D_T^{-1}\left(L_{1,T}\beta-L_{2,T}\beta-L_{3,T}+L_{4,T}+L_{5,T}-L_{6,T}\right),$$
where
$$L_{1,T}=\sum_{t=p+1}^TU_t\left(\hat{q}_{2,t}-q_{2,t}\right),\quad L_{2,T}= \sum_{t=p+1}^TW_t\left(\hat{q}_{2,t}-q_{2,t}\right)'M_tM_t'\left(\hat{q}_{2,t}-q_{2,t}\right),$$
$$L_{3,T}= \sum_{t=p+1}^TU_t\left(\hat{q}_{1,t}-q_{1,t}\right),\quad L_{4,T}= \sum_{t=p+1}^TW_t\left(\hat{q}_{2,t}-q_{2,t}\right)'M_tM_t'\left(\hat{q}_{1,t}-q_{1,t}\right),$$
$$L_{5,T}=\sum_{t=p+1}^TW_tO_tZ_t\sigma^2_t,\quad L_{6,T}= \sum_{t=p+1}^TW_t\left(\hat{q}_{2,t}-q_{2,t}\right)'M_tZ_t\sigma^2_t.$$
$L_{5,T}$ is the main term in the asymptotic expansion of the numerator $L_T$. 
We will use the formula 
\begin{equation}\label{approxmat}
B^{-1}A-b^{-1}a=b^{-1}(A-a)-b^{-1}(B-b)b^{-1}(A-a)-b^{-1}(B-b)b^{-1}a+B^{-1}(B-b)b^{-1}(B-b)b^{-1}A.
\end{equation}
Now for $j=1,2,3$, we set $\bar{S}_{j,t}=S_{j,t}-\E\left(S_{j,t}\right)$. 
Using (\ref{approxmat}), we have for $j=1,2$,

\begin{eqnarray}\label{frdec}
\hat{q}_{j,t}-q_{j,t}&=&\E^{-1}\left(S_{3,t}\right)\E\left(S_{j,t}\right)-q_{j,t}+S_{3,t}^{-1}\bar{S}_{3,t}\E^{-1}\left(S_{3,t}\right)\bar{S}_{3,t}\E^{-1}\left(S_{3,t}\right)S_{j,t}\nonumber\\
&+&\E^{-1}\left(S_{3,t}\right)\bar{S}_{j,t}-\E^{-1}\left(S_{3,t}\right)\bar{S}_{3,t}\E^{-1}\left(S_{3,t}\right)\bar{S}_{j,t}
-\E^{-1}\left(S_{3,t}\right)\bar{S}_{3,t}\E^{-1}\left(S_{3,t}\right)\E\left(S_{j,t}\right)\nonumber\\
\end{eqnarray}

To prove Theorem $1$, we will prove that 
\begin{equation}\label{etap1}
\frac{1}{\sqrt{T}}L_{5,T}\stackrel{\mathcal{D}}{\rightarrow}\mathcal{N}_n\left(0,\Sigma_2\right),
\end{equation}
\begin{equation}\label{etap2}
\frac{1}{\sqrt{T}}L_{j,T}=o_{\P}(1),\quad j\in\{1,2,3,4,6\},
\end{equation}
\begin{equation}\label{etap3}
\lim_{T\rightarrow \infty}\frac{1}{T}D_T=\Sigma_1,\mbox{ a.s.}.
\end{equation}

The proofs of (\ref{etap1}), (\ref{etap2}) and (\ref{etap3}) are established in the following subsections.

\paragraph{Proof of assertion (\ref{etap1})}

To prove (\ref{etap1}), we use the central limit theorem for triangular arrays of martingale differences (see \citet{Pol} Chapter VIII.$1$, Theorem $1$). Using the Cramer-Wold device, it is enough to prove that
\begin{equation}\label{reel}
\frac{1}{\sqrt{T}}\sum_{t=p+1}^TW_t\sigma^2_tx'O_txZ_t\stackrel{\mathcal{D}}{\rightarrow}\mathcal{N}\left(0,x'\Sigma_2x\right),
\end{equation}
for each vector $x\in\R^n$. 
We set 
$$A_t=W_t\sigma^2_tx'O_txZ_t,\quad \widetilde{A}_t=W_t\sigma^2_tx'\left(N_t-q_2\left(\frac{t}{T}\right)'M_t\right)xZ_t.$$
Then, if $\mathcal{F}_t=\sigma\left(\xi_s:s\leq t\right)$, the two families $\left\{\left(A_t,\mathcal{F}_t\right):1\leq t\leq T\right\}$ and $\left\{\left(\widetilde{A}_t,\mathcal{F}_t\right):1\leq t\leq T\right\}$ form a martingale difference. Their corresponding partial sums are asymptotically equivalent because the quantity $q_{2,t}$ is simply replaced by $q_2(t/T)$ in the expression of $\widetilde{A}_t$. Indeed, we have
\begin{eqnarray*}
\vert A_t-\widetilde{A}_t\vert\leq \Vert x\Vert^2\cdot \vert Z_t\vert\cdot\Vert q_{2,t}-q_2\left(t/T\right)\Vert\cdot\Vert M_tW_t\sigma_t^2\Vert. 
\end{eqnarray*}
Using the fact that $M_tW_t\sigma_t^2$ is bounded uniformly in $t$ and 
Lemma \ref{biais}, $1.$, we deduce that 
$\frac{1}{\sqrt{T}}\sum_{t=1}^T(A_t-\widetilde{A}_t)\rightarrow 0$,
in probability. As a consequence, it is sufficient to prove (\ref{reel}) for $\widetilde{A}_t$ instead of $A_t$.
Moreover, 
$$\E\left({\widetilde{A}_t}^2\vert \mathcal{F}_{t-1}\right)
=\E\vert \xi_0^2-1\vert^2\cdot x'W_t^2\sigma_t^4\left(N_t-q_2\left(\frac{t}{T}\right)'M_t\right)\left(N_t-q_2\left(\frac{t}{T}\right)'M_t\right)'x.$$

The process $G$ defined by
$G_t=W_t^2\sigma_t^4\left(N_t-q_2\left(\frac{t}{T}\right)'M_t\right)\left(N_t-q_2\left(\frac{t}{T}\right)'M_t\right)'$
is (coordinatewise) a process of type $I$ and  Proposition \ref{cumulantcontrol} (point $4$) leads to 
\\$\lim_{T\rightarrow \infty}\frac{1}{T}\sum_{t=p+1}^T \overline{G}_t=0$ in probability. 
Moreover, using Lemma \ref{approxstat}, Lemma \ref{biais} (\ref{bor2}) and some Lipschitz properties, one can show that
$$\lim_{T\rightarrow \infty}\frac{1}{T}\sum_{t=p+1}^T\E\vert\xi_0^2-1\vert^2\cdot\E(G_t)=\lim_{T\rightarrow\infty}\frac{1}{T}\sum_{t=p+1}^T\E\vert\xi_0^2-1\vert^2\E\left(G_1(t/T)\right)
= x'\Sigma_2x.$$
Then we get 
$\frac{1}{T}\sum_{t=p+1}^T\E\left({\widetilde{A}_t}^2\vert \mathcal{F}_{t-1}\right)\rightarrow x'\Sigma_2x$,
in probability. Next, we check the Lindberg condition. If $\epsilon>0$, we have 
$$\E\left({\widetilde{A}_t}^2\mathds{1}_{\vert \widetilde{A}_t\vert>\epsilon \sqrt{T}}\vert \mathcal{F}_{t-1}\right)\\
\leq \frac{1}{\epsilon^{\delta}T^{\frac{\delta}{2}}}\E\vert \xi_0^2-1\vert^{2+\delta}W_t^{2+\delta}\sigma_t^{4+2\delta}\vert x'\left(N_t-q_2\left(\frac{t}{T}\right)'M_t\right)x\vert^{2+\delta}.$$
We easily deduce that 
$\frac{1}{T}\sum_{t=p+1}^T \E\left({\widetilde{A}_t}^2\mathds{1}_{\vert \widetilde{A}_t\vert>\epsilon}\vert\mathcal{F}_{t-1}\right)\rightarrow 0$
in probability. This proves (\ref{reel}) and (\ref{etap1}) follows.

\paragraph{Proof of $\frac{1}{\sqrt{T}}N_{j,T}=o_{\P}(1)$ for $j\in\{1,3,6\}$}\label{covbiais}.

We only prove the result for $j=3$. The two other cases can be treated in the same way.
We set $\eta_{1,t}=\E^{-1}(S_{3,t})\E\left(S_{1,t}\right)-q_{1,t}\in\mathcal{M}_{j,1}$. 
The proof of the result follows from the following points.
\begin{enumerate}
\item
We first prove that $\frac{1}{\sqrt{T}}\sum_{t=p+1}^TU_t\eta_{1,t}=o_{\P}(1)$. 
It is enough to prove that for $\left(w,y\right)\in \llbracket 1,n\rrbracket\times\llbracket 1,m\rrbracket$,
$$z_T=\frac{1}{\sqrt{T}}\sum_{t=p+1}^T U_t(w,y)\eta_{1,t}(y,1)=o_{\P}(1).$$
But this assertion follows from Lemma \ref{biais} ($3.$) and Lemma \ref{inter}, since $\E z_T^2$ can be bounded by $b^2$ (up to a positive constant).

\item
Now we prove that 
\begin{equation}\label{ccal}
\frac{1}{\sqrt{T}}\sum_{t=p+1}^T U_t \E^{-1}(S_{3,t})\bar{S}_{1,t}=o_{\P}(1).
\end{equation}
In order to prove (\ref{ccal}), it is enough to prove that for a given vector
$(w,y,z)\in\llbracket 1,n\rrbracket\times\llbracket 1,m\rrbracket\times\llbracket 1,m\rrbracket$, 
$$\frac{1}{\sqrt{T}}\sum_{t=p+1}^TU_t(w,y)\E^{-1}(S_{3,t})(y,z)\bar{S}_{1,t}(z,1)=o_{\P}(1).$$
But the result is a direct consequence of Lemma \ref{biais} (point $2$) and Proposition \ref{cumulantcontrol} (point $2$ applied with $h=1$).

\item
Using the same arguments as for point $2$, one can easily show that
$$\sum_{t=p+1}^TU_t \E^{-1}(S_{3,t})\bar{S}_{3,t}\E^{-1}(S_{3,t})\bar{S}_{1,t}=o_{\P}\left(\sqrt{T}\right),\quad 
\sum_{t=p+1}^TU_t \E^{-1}(S_{3,t})\bar{S}_{3,t}\E^{-1}(S_{3,t})\E\left(S_{1,t}\right)=o_{\P}\left(\sqrt{T}\right).$$
\item
Finally, we prove that 
\begin{equation}\label{ccal2}
\frac{1}{\sqrt{T}}\sum_{t=p+1}^TU_tS_{3,t}^{-1}\bar{S}_{3,t}\E^{-1}(S_{3,t})\bar{S}_{3,t}\E^{-1}(S_{3,t})S_{1,t}=o_{\P}(1).
\end{equation}
Since $U_t$ is uniformly bounded in $\omega,t,T$, there exists $C>0$ such that
\begin{eqnarray*} 
&&\Vert\frac{1}{\sqrt{T}}\sum_{t=p+1}^TU_tS_{3,t}^{-1}\bar{S}_{3,t}\E^{-1}(S_{3,t})\bar{S}_{3,t}\E^{-1}(S_{3,t})S_{1,t}\Vert\\
&\leq& \frac{C}{\sqrt{T}}\max_{p+1\leq i\leq T}\Vert S_{3,i}^{-1}\Vert\cdot\sum_{t=p+1}^T \Vert \bar{S}_{3,t}\Vert^2\cdot\Vert S_{1,t}\Vert\\
&\leq& C\max_{p+1\leq i\leq T}\Vert S_{3,i}^{-1}\Vert\cdot \sqrt{\frac{1}{T}\sum_{t=p+1}^T \Vert \bar{S}_{3,t}\Vert^4}\cdot\sqrt{\sum_{t=p+1}^T \Vert S_{1,t}\Vert^2}.
\end{eqnarray*}  
From Proposition \ref{cumulantcontrol}, we have the bounds
$$\max_{p+1\leq t\leq T}\E\Vert\bar{S}_{3,t}\Vert^4=O\left(\frac{1}{T^2b^2}\right),\quad \max_{p+1\leq t\leq T}\E\Vert\bar{S}_{1,t}\Vert^2=O\left(\frac{1}{Tb}\right).$$
 Then using the point $2$ of Lemma \ref{biais},
we conclude that $\frac{1}{T}\sum_{t=p+1}^T \Vert S_{1,t}\Vert^2=O_{\P}(1)$. Finally, using Lemma \ref{biais} ($4.$), we conclude that
$$\Vert\frac{1}{\sqrt{T}}\sum_{t=p+1}^TU_tS_{3,t}^{-1}\bar{S}_{3,t}\E^{-1}(S_{3,t})\bar{S}_{3,t}\E^{-1}(S_{3,t})S_{1,t}\Vert=O_{\P}\left(\frac{1}{\sqrt{T}b}\right).$$
Hence, (\ref{ccal2}) follows using the assumption $b\sqrt{T}\rightarrow \infty$.
\end{enumerate}

\paragraph{Proof of $\frac{1}{\sqrt{T}}N_{j,T}=o_{\P}(1)$ for $j\in\{2,4\}$}.

We only prove the result for $j=4$, the proof for the case $j=2$ is similar.
Here, we only use the basic decompositions
\begin{equation}\label{decsimple}
\hat{q}_{k,t}-q_{k,t}=\eta_{k,t}+\E^{-1}(S_{3,t})\bar{S}_{k,t}-S_{3,t}^{-1}\bar{S}_{3,t}\E^{-1}(S_{3,t})S_{k,t},
\end{equation}
for $k=1,2$ and with
$$\eta_{k,t}=\E^{-1}(S_{3,t})\E\left(S_{k,t}\right)-q_{k,t}.$$
Then we have
\begin{eqnarray*}
\Vert N_{4,T}\Vert&\leq& \sum_{t=p+1}^T\Vert \hat{q}_{2,t}-q_{2,t}\Vert\cdot\Vert W_tM_tJ'_t\Vert\cdot\Vert \hat{q}_{1,t}-q_{1,t}\Vert\\
&\leq& \frac{C_1}{2}\sum_{t=p+1}^T\left\{\Vert \hat{q}_{2,t}-q_{2,t}\Vert^2+\Vert \hat{q}_{1,t}-q_{1,t}\Vert^2\right\},
\end{eqnarray*}
where $C_1=\sup_{\omega,T,t}\Vert W_t M_tM_t'\Vert$. 
It remains to prove that for $k=1,2$,
\begin{equation}\label{conv+}
\frac{1}{\sqrt{T}}\sum_{t=p+1}^T\Vert \hat{q}_{k,t}-q_{k,t}\Vert^2=o_{\P}(1).
\end{equation}
We only prove (\ref{conv+}) for $k=1$, the proof for $k=2$ being the same.
The proof easily follows from the three following points.
\begin{enumerate}
\item
From Lemma \ref{biais} ($3.$) and our bandwidth condition, we have $\frac{1}{\sqrt{T}}\sum_{t=p+1}^T\frac{\Vert \eta_{1,t}\Vert^2}{\sqrt{T}}=o(1).$
\item
Since $$C_2=\sup_{\substack{T\geq p+1\\p+1\leq t\leq T}}\Vert \E^{-1}(S_{3,t})\Vert$$ is finite using Lemma \ref{biais} (point $2$),
we use the inequality $$\frac{1}{\sqrt{T}}\sum_{t=p+1}^T\Vert \E^{-1}(S_{3,t})\bar{S}_{1,t}\Vert^2\leq \frac{C_2^2}{\sqrt{T}}\sum_{t=p+1}^T \Vert \bar{S}_{1,t}\Vert^2.$$
But we know from Proposition \ref{cumulantcontrol} that $\max_{p+1\leq t\leq T}\E\left(\Vert\bar{S}_{1,t}\Vert^2\right)=O\left(\frac{1}{Tb}\right)$. 
Then condition $b\sqrt{T}\rightarrow \infty$ entails that 
$$\frac{1}{\sqrt{T}}\sum_{t=p+1}^T\E\left(\Vert\bar{S}_{1,t}\Vert^2\right)=o(1).$$
This shows that $\frac{1}{\sqrt{T}}\sum_{t=p+1}^T\Vert \E^{-1}(S_{3,t})\bar{S}_{1,t}\Vert^2=o_{\P}(1)$. 
\item
Finally, we show that
\begin{equation}\label{ccalc3}
\frac{1}{\sqrt{T}}\sum_{t=p+1}^T \Vert S_{3,t}^{-1}\bar{S}_{3,t}m_t^{-1}S_{1,t}\Vert^2=o_{\P}(1).
\end{equation}
We have $\Vert S_{3,t}^{-1}\bar{S}_{3,t}\E^{-1}(S_{3,t})S_{1,t}\Vert^2\leq C_2^2\max_{p+1\leq i\leq T}\Vert S_{3,i}^{-1}\Vert^2 \Vert \bar{S}_{3,t}\Vert^2
\Vert S_{1,t}\Vert^2$,
where $C_2>0$ is defined in the previous point. We have $\max_{p+1\leq i\leq T}\Vert S_{3,i}^{-1}\Vert^2=O_{\P}(1)$ (see Lemma \ref{biais}, point $4$) and 
$$\E\left(\Vert \bar{S}_{3,t}\Vert^2\cdot\Vert S_{1,t}\Vert^2\right)\leq \E^{\frac{\delta}{1+\delta}}\left(\Vert \bar{S}_{3,t}\Vert^{2\frac{1+\delta}{\delta}}\right)\cdot\E^{\frac{1}{1+\delta}}\left(\Vert S_{1,t}\Vert^{2(1+\delta)}\right).$$
Without loss of generality, one can assume that $\frac{1+\delta}{\delta}$ is an integer.
Using Proposition \ref{cumulantcontrol} (point $3$), we have
$\max_{p+1\leq t\leq T}\E\left(\Vert \bar{S}_{3,t}\Vert^{2\frac{1+\delta}{\delta}}\right)=O\left((Tb)^{-\frac{1+\delta}{\delta}}\right).$
Moreover, using convexity, the moment assumption on the noise and the fact that $W_i M_i \sigma_i^2$ is bounded,
$$\max_{p+1\leq t\leq T}\E\Vert S_{1,t}\Vert^{2(1+\delta)}\leq \max_{p+1\leq i\leq T}\E \Vert W_i M_i X_i^2\Vert^{2(1+\delta)}=O(1).$$
Then assertion (\ref{ccalc3}) easily follows from the condition $\sqrt{T}b\rightarrow \infty$.
\end{enumerate}
\subsection{Proof of assertion \ref{etap3}}\label{denominateur}
Recalling that $\hat{O}_t=O_t-\left(\hat{q}_{2,t}-q_{2,t}\right)'M_t$, we have 
\begin{eqnarray*}
D_T&=& \frac{1}{T}\sum_{t=p+1}^TW_tO_tO_t'+\frac{1}{T}\sum_{t=p+1}^TW_t\left(\hat{q}_{2,t}-q_{2,t}\right)'M_tM_t'\left(\hat{q}_{2,t}-q_{2,t}\right)\\
&-& \frac{1}{T}\sum_{t=p+1}^TW_t\left(\hat{q}_{2,t}-q_{2,t}\right)'M_t O_t'-\frac{1}{T}\sum_{t=p+1}^TW_tO_tM_t'\left(\hat{q}_{2,t}-q_{2,t}\right).
\end{eqnarray*}
We have already shown that the three last terms in the previous decomposition are $o_{\P}(1)$. Then, it remains to show that
\begin{equation}\label{ccalc4}
\lim_{T\rightarrow\infty}\frac{1}{T}\sum_{t=p+1}^TW_tO_tO_t'=\Sigma_1,
\end{equation}
in probability. One can obtain (\ref{ccalc4}) using the same arguments as for deriving the limit of 
$$\frac{1}{T}\sum_{t=p+1}^T\E\left({\widetilde{A}_t}^2\vert \mathcal{F}_{t-1}\right)$$ in the proof of  assertion (\ref{etap1}), .$\square$

\section{Auxiliary results for the proof of Theorem $2$}
\begin{cor}\label{smallsecondmoment}
Let $(Y_t)_{1\leq t\leq T}$ be a process of type $II$. We denote by $(\bar{S}_t)_{1\leq t\leq T}$ the corresponding smoothing. Then, under the assumptions of Theorem $2$, we have
$$\max_{1\leq t\leq T}\vert \bar{S}_t\vert=o_{\P}\left(T^{-\frac{1}{4}}\right).$$
Consequently, if for $i\in\{1,2\}$, $\bar{S}^{(i)}$ is a smoothing then 
$$\max_{1\leq t_1,t_2\leq T}\vert \bar{S}_{t_1}^{(1)} \bar{S}_{t_2}^{(2)}\vert=o_{\P}\left(\frac{1}{\sqrt{T}}\right).$$
\end{cor}
  
\paragraph{Proof of Corollary \ref{smallsecondmoment}}

We set $j=2 \left(\left[\frac{1}{\tau}\right]+1\right)$. Then using the point $4$ of Proposition $3$, we have 

\begin{eqnarray*}
&&\P\left(\max_{p+1\leq t\leq T}\vert \bar{S}_t\vert >\epsilon T^{-\frac{1}{4}}\right)\\
&\leq& \frac{C T^{\frac{j}{4}}}{\epsilon^j}\sum_{t=p+1}^T \E\vert\bar{S}_t\vert^j\\
&\leq& C\frac{T^{\frac{j}{4}}}{(Tb)^{\frac{j}{2}}}.
\end{eqnarray*}
where $C>0$ denotes a generic constant and $\epsilon>0$ is arbitrary. Since, $\frac{T^{\frac{j}{4}+1}}{(Tb)^{\frac{j}{2}}}\leq\frac{1}{\left(T^{\frac{1}{2}-\tau}b\right)^{\frac{j}{2}}}$, the result follows from the assumptions made on $b$.$\square$

Finally, we state a Lemma which will be useful for the proof of Theorem $2$. For a matrix-valued process $(H_t)_t$, measurable with the sigma field $\mathcal{F}_T$, we will use the equality $H_t=o_{\P}\left(\frac{1}{\sqrt{T}}\right)$ when $\max_{p+1\leq t\leq T}\Vert H_t\Vert=o_{\P}\left(\frac{1}{\sqrt{T}}\right)$.  
We also introduce additional notations. For $j=1,2,3$, we define $S^{*}_{j,t}$
as $S_{j,t}$, replacing $W_t$ with $\frac{1}{\sigma_t^4}$.
\begin{lem}\label{lemmeutile}
Assume that the assumptions of Theorem $2$ hold.
\begin{enumerate}
\item
For $j=1,2$, we have
$$\hat{q}_{j,t}-q_{j,t}=\eta_{j,t}+\E^{-1}(S_{3,t})\bar{S}_{j,t}-\E^{-1}(S_{3,t})\bar{S}_{3,t}\E^{-1}(S_{3,t})\E(S_{j,t})+o_{\P}\left(\frac{1}{\sqrt{T}}\right),$$
and 
$$\Vert \hat{q}_{j,t}-q_{j,t}\Vert^2=o_{\P}\left(\frac{1}{\sqrt{T}}\right).$$
\item
We have 
$$\hat{W}^{*}_t=\frac{1}{\sigma_t^4}\left\{1-\frac{2}{\sigma^2_t}\left(M_t'L_t+O'_t\left(\hat{\beta}-\beta\right)\right)+o_{\P}\left(\frac{1}{\sqrt{T}}\right)\right\},$$
with
$$L_t=\hat{q}_{1,t}-q_{1,t}-\left(\hat{q}_{2,t}-q_{2,t}\right)\beta.$$
In particular, we have $\hat{W}^{*}_t=\frac{1}{\sigma_t^4}\left(1+E_t\right)$ with $E_t^2=o_{\P}\left(\frac{1}{\sqrt{T}}\right)$. 
\item
We have $\max_{p+1\leq t\leq T}\Vert \hat{S}_{3,t}^{-1}\Vert =O_{\P}(1).$
\item
For $j=1,2,3$, there exists two matrices $a_{j,t}$ and $c_{j,t}$ such that
\begin{equation}\label{oubli}
\hat{s}^{*}_{j,t}-\E\left(S_{j,t}^{*}\right)=\overline{S^{*}}_{j,t}+R_{j,t}+o_{\P}\left(\frac{1}{\sqrt{T}}\right),
\end{equation}
where 
$$R_{j,t}=a_{j,t}\left(\hat{\beta}-\beta\right)+\sum_{i=p+1}^Tk_{t,i}c_{j,i}L_i=o_{\P}\left(T^{-1/4}\right).$$
Moreover, $\Vert \hat{s}^{*}_{j,t}-\E\left(S_{j,t}^{*}\right)\Vert^2=o_{\P}\left(\frac{1}{\sqrt{T}}\right)$.
\item

For $j=1,2$, we have $\hat{q}^{*}_{j,t}-q_{j,t}^{*}=F_{j,t}^{*}+\Delta_{j,t}$, where 
$$F_{j,t}^{*}=\eta_{j,t}^{*}+\E^{-1}(S_{3,t}^{*})\overline{S^*}_{j,t}-\E^{-1}(S_{3,t}^{*})\overline{S^*}_{3,t}\E^{-1}(S_{3,t}^{*})\E\left(S_{j,t}^{*}\right),$$
$$\Delta_{j,t}=\E^{-1}(S_{3,t}^{*})R_{j,t}-\E^{-1}(S_{3,t}^{*})R_{3,t}\E^{-1}(S_{3,t}^{*})\E\left(S_{j,t}^*\right)+\o.$$
Moreover, $\Vert \hat{q}^{*}_{j,t}-q_{j,t}^{*}\Vert^2=\o$.
\end{enumerate}
\end{lem}

\paragraph{Proof of Lemma \ref{lemmeutile}}
\begin{enumerate}
\item
Under the assumptions of Theorem $2$, we recall that $\bar{S}_t=o_{\P}\left(T^{-1/4}\right)$ for all univariate smoothing (see Corollary \ref{smallsecondmoment}). Then, applying this property coordinatewise, we have $\Vert \overline{S}_{j,t}\Vert=o_{\P}\left(T^{-1/4}\right)$ for $j=1,2,3$. The announced decomposition is then a consequence of decomposition $(6.4)$ given in the paper and of Lemma \ref{biais}, point $4$. Next, the second assertion follows from the condition
$b^2\sqrt{T}\rightarrow 0$ and Lemma \ref{biais} (point $3$). 
\item
We use the decomposition
\begin{equation}\label{poidsdec}
\hat{W}^{*}_t=\frac{1}{\sigma_t^4}\left\{1+\frac{\sigma_t^4-\hat{\sigma}_t^4-\nu_T}{\sigma_t^4}+\frac{\left(\sigma_t^4-\hat{\sigma}_t^4-\nu_T\right)^2}{\sigma_t^4\left(\hat{\sigma}_t^4+\nu_T\right)}\right\}.
\end{equation}
Now, using the fact that $L_t^2+\Vert \hat{\beta}-\beta\Vert^2=o_{\P}\left(\frac{1}{\sqrt{T}}\right)$, we get 
$$\hat{\sigma}^4_t=\sigma_t^4\left(1+\frac{2}{\sigma_t^2}\left(O'_t\left(\hat{\beta}-\beta\right)+M_t'L_t\right)+o_{\P}\left(\frac{1}{\sqrt{T}}\right)\right).$$
From this decomposition, we deduce that $\frac{1}{\hat{\sigma}_t^4+\nu_T}=O_{\P}(1)$. 
Then we get 
$$\hat{W}^{*}_t=\frac{1}{\sigma_t^4}\left\{1-\frac{2}{\sigma_t^2}\left(O'_t\left(\hat{\beta}-\beta\right)+M_t'L_t\right)+o_{\P}\left(\frac{1}{\sqrt{T}}\right)\right\},$$
which also yields to the approximation $\hat{W}^{*}_t=\frac{1}{\sigma_t^4}\left(1+E_t\right)$ with $E_t^2=o_{\P}\left(\frac{1}{\sqrt{T}}\right)$.

\item
We have
\begin{eqnarray*}
\hat{s}^{*}_{3,t}&=&S_{3,t}^{*}+\sum_{i=p+1}^Tk_{t,i}\frac{M_iM_i'}{\sigma_i^4}E_i\\
&=& S_{3,t}^{*}+o_{\P}(1).
\end{eqnarray*}
Indeed, $\max_{p+1\leq t\leq T}\frac{\Vert M_iM_i'\Vert}{\sigma_i^4}$ is bounded uniformly in $(\omega,T)$ and from point $2$, 
$\max_{p+1\leq t\leq T}\Vert E_t\Vert=o_{\P}(1)$. Then the result follows from the fact that 
$\max_{p+1\leq t\leq T}\Vert (S_{3,t}^{*})^{-1}\Vert=O_{\P}(1)$.  
\item
We only consider the case $j=1$, the cases $j=2$ or $j=3$ being similar. 
Using the point $2$, we have the decomposition 
$$\hat{s}^{*}_{1,t}=S_{1,t}^{*}-2\sum_{i=p+1}^Tk_{t,i}\frac{O'_i\left(\hat{\beta}-\beta\right)+M_i'L_i}{\sigma_i^6}M_iX_i^2+o_{\P}\left(\frac{1}{\sqrt{T}}\right).$$
Now we set $a_{1,t}=-2\sum_{i=p+1}^Tk_{t,i}\E\left(\frac{M_iO'_iX_i^2}{\sigma_i^6}\right)$
and $c_{1,t}=-2\sum_{i=p+1}^Tk_{t,i}\E\left(\frac{M_iM_i'X_i^2}{\sigma_i^6}\right)$. 
Moreover, we have from Lemma \ref{smallsecondmoment}
$$-2\sum_{i=p+1}^Tk_{t,i}\frac{M_iN_i'X_i^2}{\sigma_i^6}-a_{1,t}=o_{\P}\left(T^{-\frac{1}{4}}\right),$$
which leads to 
$$\left(-2\sum_{i=p+1}^Tk_{t,i}\frac{M_iN_i'X_i^2}{\sigma_i^6}-a_{1,t}\right)\cdot\left(\hat{\beta}-\beta\right)=o_{\P}\left(\frac{1}{\sqrt{T}}\right).$$
Then it remains to show that 
\begin{equation}\label{mom+}
\sum_{i=p+1}^Tk_{t,i}\left(\frac{M_iM_i'X_i^2}{\sigma_i^6}-c_{1,i}\right)L_i=o_{\P}\left(\frac{1}{\sqrt{T}}\right).
\end{equation}
Considering the decomposition given in point $1$, assertion (\ref{mom+}) will follow if we show the two following assertions. For all real-valued sequences $(c_t)$ such that 
$\max_{p+1\leq t\leq T}c_t=O(b)$ and all real-valued processes $(Y_t)$, $(G_t)$ of type $I$ or $II$, 
\begin{equation}\label{mom++}
\sum_{i=p+1}^Tk_{t,i}c_i \bar{Y}_i=o_{\P}\left(\frac{1}{\sqrt{T}}\right),
\end{equation}
and
\begin{equation}\label{mom+++}
\sum_{i=p+1}^Tk_{t,i}\bar{Y}_i\bar{S}_i=o_{\P}\left(\frac{1}{\sqrt{T}}\right),
\end{equation}
where 
$$\bar{S}_t=\sum_{i=p+1}^Tk_{t,i}\left(G_i-\E(G_i)\right).$$
The assertion (\ref{mom+++}) can be proved as follows. For $\epsilon>0$ and an even integer $h>0$, we have
\begin{eqnarray*}
\P\left(\max_{p+1\leq t\leq T}\vert \sum_{i=p+1}^T k_{t,i} \bar{Y}_i\bar{S}_i\vert >\epsilon \frac{1}{\sqrt{T}}\right)&\leq&
\frac{T^{\frac{h}{2}}}{\epsilon^h}\sum_{t=p+1}^T\E\vert \sum_{i=p+1}^Tk_{t,i}\bar{Y}_i\bar{S}_i\vert^h\\
&\leq& \frac{T^{\frac{h}{2}+1}}{\epsilon^h}\sum_{p+1\leq i_1,j_1,\ldots,i_h,j_h\leq T}p_{i_1}p_{j_1}\cdots p_{i_h}p_{j_h}\vert \E\left(\bar{Y}_{i_1}\bar{G}_{j_1}\cdots \bar{Y}_{i_h}\bar{G}_{j_h}\right)\vert,\\
\end{eqnarray*}
where $p_i=\max_{p+1\leq t\leq T}k_{t,i}$. 
Using Lemma $4$, we deduce that the right hand side of the previous inequality is $O\left(\frac{1}{T^{\frac{h}{2}-1}b^h}\right)$
But $\frac{h}{2}-1=\frac{h}{2}\left(1-\frac{2}{h}\right)\geq \frac{h}{2}\left(1-\tau\right)$ when $h\geq \frac{2}{\tau}$. 
Using the bandwidth conditions, we get (\ref{mom+++}).
Next, using Corollary \ref{smallsecondmoment} and the bandwidth conditions, we get
$$\max_{p+1\leq t\leq T}\vert \sum_{i=p+1}^Tk_{t,i}c_i\bar{S}_i\vert\leq \max_{p+1\leq t\leq T}\vert c_t\vert\times \max_{p+1\leq t\leq T}\vert\bar{S}_t\vert=o_{\P}\left(bT^{-\frac{1}{4}}\right)=o_{\P}\left(\frac{1}{\sqrt{T}}\right).$$
This proves (\ref{mom++}).
Finally, since $\Vert L_t\Vert^2=o_{\P}\left(\frac{1}{\sqrt{T}}\right)$ and $\hat{\beta}-\beta=O_{\P}\left(\frac{1}{\sqrt{T}}\right)$, the decomposition (\ref{oubli}) holds true. Using some arguments discussed before, we easily get 
$\Vert \hat{s}^{*}_{1,t}-\E\left(S_{1,t}^{*}\right)\Vert^2=o_{\P}\left(\frac{1}{\sqrt{T}}\right)$. 
\item
Using the equality
\begin{eqnarray*}
\hat{q}^{*}_{j,t}-q_{j,t}^{*}&=&\eta_{j,t}^{*}+\E^{-1}(S_{3,t}^{*})\left(\hat{s}^{*}_{j,t}-\E\left(S_{j,t}^{*}\right)\right)\\
&-&\E^{-1}(S_{3,t}^{*})\left(\hat{s}^{*}_{j,t}-\E\left(S_{j,t}^{*}\right)\right)\E^{-1}(S_{3,t}^{*})\E\left(S_{j,t}^{*}\right)\\
&+&\left(\hat{s}^{*}_{3,t}\right)^{-1}\left(\hat{s}^{*}_{3,t}-\E\left(S_{j,t}^{*}\right)\right)\E^{-1}(S_{3,t}^{*})\left(\hat{s}^{*}_{3,t}-\E\left(S_{3,t}^{*}\right)\right)\E^{-1}(S_{3,t}^{*})\hat{s}^{*}_{j,t},
\end{eqnarray*}
the result of point $5$ is an easy consequence of the previous points.$\square$
\end{enumerate}

\section{Proof of Theorem $2$}
We recall that for a triangular array $\{H_t=H_{t,T}\}$ of matrices, we will denote $H_t=o_{\P}\left(\frac{1}{\sqrt{T}}\right)$ when 
$$\max_{p+1\leq t\leq T}\Vert H_t\Vert=o_{\P}\left(\frac{1}{\sqrt{T}}\right).$$
\paragraph{Notations.} Let us also recall the following notations.
$$S_{1,t}^{*}=\sum_{t=p+1}^Tk_{t,i}\frac{M_i X_i^2}{\sigma_i^4},\quad S_{2,t}^{*}=\sum_{t=p+1}^Tk_{t,i}\frac{M_i N_i'}{\sigma_i^4},$$
$$S_{3,t}^{*}=\sum_{t=p+1}^Tk_{t,i}\frac{M_i M_i'}{\sigma_i^4},$$
$$\eta_{1,t}^{*}=\E^{-1}(S^{*}_{3,t})\E\left(S^{*}_{1,t}\right)-q_{1,t}^{*},$$
$$\eta_{2,t}^{*}=\E^{-1}(S^{*}_{3,t})\E\left(S^{*}_{2,t}\right)-q_{2,t}^{*},$$
$$q_{1,t}^{*}=\E^{-1}\left(\frac{M_tM_t'}{\sigma_t^4}\right)\E\left(\frac{M_t X_t^2}{\sigma_t^4}\right),\quad q_{2,t}^{*}=\E^{-1}\left(\frac{M_tM_t'}{\sigma_t^4}\right)\E\left(\frac{M_t N_t'}{\sigma_t^4}\right),$$
$$O_t^{*}=N_t-\left(q_{2,t}^{*}\right)'M_t.$$
The proof of Theorem $2$ uses the same decomposition as for Theorem $1$. More precisely, we have
$$\hat{\beta}_{*}-\beta=\hat{D}_T^{-1}\left(\hat{L}_{1,T}\beta-\hat{L}_{2,T}\beta-\hat{L}_{3,T}+\hat{L}_{4,T}+\hat{L}_{5,T}-\hat{L}_{6,T}\right),$$
where
\begin{eqnarray*}
\hat{L}_{1,T}&=&\sum_{t=p+1}^T\hat{W}^{*}_tO_t^{*}M_t'\left(\hat{q}^{*}_{2,t}-q^{*}_{2,t}\right),\\
\hat{L}_{2,T}&=& \sum_{t=p+1}^T\hat{W}^{*}_t\left(\hat{q}^{*}_{2,t}-q^{*}_{2,t}\right)'M_tM_t'\left(\hat{q}^{*}_{2,t}-q^{*}_{2,t}\right),\\
\hat{L}_{3,T}&=& \sum_{t=p+1}^T\hat{W}^{*}_tO_t^{*}M_t'\left(\hat{q}^{*}_{1,t}-q^{*}_{1,t}\right),\\
\hat{L}_{4,T}&=& \sum_{t=p+1}^T\hat{W}^{*}_t\left(\hat{q}^{*}_{2,t}-q^{*}_{2,t}\right)'M_tM_t'\left(\hat{q}^{*}_{1,t}-q^{*}_{1,t}\right),\\
\hat{L}_{5,T}&=&\sum_{t=p+1}^T\hat{W}^{*}_tO^{*}_tZ_t\sigma^2_t,\\
\hat{L}_{6,T}&=& \sum_{t=p+1}^T\hat{W}^{*}_t\left(\hat{q}^{*}_{2,t}-q^{*}_{2,t}\right)'M_tZ_t\sigma^2_t.\\
\hat{D}_T=\sum_{t=p+1}^T\hat{W}^{*}_tO^{*}_t\left(O^{*}_t\right)'-\hat{L}_{1,T}-\hat{L}_{1,T}'+\hat{L}_{2,T}.
\end{eqnarray*}
Using Lemma \ref{lemmeutile}, it is easy to get 
$$\Vert \hat{L}_{2,T}\Vert+\Vert \hat{L}_{4,T}\Vert=o_{\P}\left(\sqrt{T}\right).$$
Then the proof of Theorem $2$ will follow from the three following points.
\begin{enumerate}
\item
We first show that $L_{j,T}=o_{\P}\left(\sqrt{T}\right)$ for $j=1,3,6$. 
We only consider the case $j=1$, the proofs for the cases $j=3$, or $j=6$ being similar.
Using the notations introduced in Lemma \ref{lemmeutile} (points $2$ and $5$), we have
$$\hat{L}_{1,T}=\sum_{t=p+1}^T\left(1+E_t\right)\cdot\frac{O_t^{*}M_t'}{\sigma_t^4}\cdot\left(F_{2,t}^{*}+\Delta_{2,t}\right).$$
The proof will be a consequence of the following points.
\begin{itemize}
\item
Since $\frac{O_t^{*}M_t'}{\sigma_t^4}$ is centered, we get $\frac{1}{\sqrt{T}}\sum_{t=p+1}^T\frac{O_t^{*}M_t'}{\sigma_t^4}F_{2,t}^{*}=o_{\P}(1)$
The proof is similar to the proof given in Theorem $1$ (see the proof of $\frac{1}{\sqrt{T}}N_{j,T}=o_{\P}(1)$ for $j=1,3,6$). 
\item
Next, we prove that $\frac{1}{\sqrt{T}}\sum_{t=p+1}^T\frac{O^{*}_tM_t'}{\sigma_t^4}\Delta_{2,t}=o_{\P}(1)$. 
First we will show that
$$\frac{1}{\sqrt{T}}\sum_{t=p+1}^T\frac{O_t^{*}M_t'}{\sigma_t^4}\E^{-1}(S_{3,t}^{*})R_{2,t}=o_{\P}(1).$$
Using Lemma \ref{lemmeutile}, we have the expression 
$$R_{2,t}=a_{2,t}\left(\hat{\beta}-\beta\right)+\sum_{i=p+1}^Tk_{t,i}c_{2,i}L_i.$$
Since $\frac{O_t^{*}M_t'}{\sigma_t^4}$ is centered, we have 
$$\frac{1}{\sqrt{T}}\sum_{t=p+1}^T\frac{W_t^{*}M_t'}{\sigma_t^4}\E^{-1}(S_{3,t}^{*})a_{2,t}\left(\hat{\beta}-\beta\right)=o_{\P}(1).$$
Then, it remains to show that 
\begin{equation}\label{notsocomplicated}
\frac{1}{\sqrt{T}}\sum_{t=p+1}^T\frac{O_t^{*}M_t'}{\sigma_t^4}\E^{-1}(S^{*}_{3,t})\sum_{i=p+1}^Tk_{t,i}c_{2,i}L_i=o_{\P}(1).
\end{equation}
The assertion (\ref{notsocomplicated}) will follow if we work with each entry of different matrices and if we prove the following property.  If $(Y_t)$ is a real-valued centered process of type I or II, $\left(\bar{S}_t\right)$ is a real-valued smoothing of type $I$ or $II$ and  
$(c_t)$, $(\widetilde{c}_t)$ are deterministic sequences satisfying $\max_{p+1\leq t\leq T}\vert c_t\vert=O(b)$ and $\max_{p+1\leq t\leq T}\vert \widetilde{c}_t\vert=O(1)$, then 
\begin{equation}\label{notsocomplicated1}
\frac{1}{\sqrt{T}}\sum_{t=p+1}^TY_tc_t=o_{\P}(1),\quad \frac{1}{\sqrt{T}}\sum_{t=p+1}^TY_t\sum_{i=p+1}^Tk_{t,i}\widetilde{c}_i\bar{S}_i=o_{\P}(1).
\end{equation}
The first assertion in (\ref{notsocomplicated}) has been already proved using the bound for the covariance function of the process $(Y_t)$ (see the proof of $\frac{1}{\sqrt{T}}L_{3,T}=o_{\P}(1)$ in the proof of Theorem $1$).
Then it remains to prove the second assertion in (\ref{notsocomplicated1}). 
Writing $\bar{S}_i=\sum_{j=p+1}^Tk_{i,j}\overline{Y'}_j$, we have 
$$\quad \frac{1}{\sqrt{T}}\sum_{t=p+1}^TY_t\sum_{i=p+1}^Tk_{t,i}\widetilde{c}_i\bar{S}_i=\frac{1}{\sqrt{T}}\sum_{t,j=p+1}^Tp_{t,j}\bar{Y}_t\overline{Y'}_i,$$
with $p_{t,j}=\sum_{i=p+1}^T\widetilde{c}_ik_{t,i}k_{i,j}$ satisfies $\max_{p+1\leq j,t\leq T}\vert p_{t,j}\vert=O\left(\frac{1}{Tb}\right)$.
Moreover, using Lemma $2$, we have for a generic constant $C>0$,
\begin{eqnarray*}
\E\vert\frac{1}{\sqrt{T}}\sum_{t,j=p+1}^Tp_{t,j}\bar{Y}_t\overline{Y'}_i\vert^2&=&\frac{1}{T}\sum_{t_1,t_2,j_1,j_2=p+1}^Tp_{t_1,j_1}p_{t_2,j_2}\vert \E\left(\bar{Y}_{t_1}\overline{Y'}_{t_2}\bar{Y}_{j_1}\overline{Y'}_{j_2}\right)\vert\\
&\leq& \frac{C}{T^3b^2}\sum_{t_1,t_2,j_1,j_2=p+1}^T\vert \E\left(\bar{Y}_{t_1}\overline{Y'}_{t_2}\bar{Y}_{j_1}\overline{Y'}_{j_2}\right)\vert\\
&\leq& \frac{C}{Tb^2}.
\end{eqnarray*}
Using the assumption $T b^2\rightarrow \infty$, we deduce the second assertion in (\ref{notsocomplicated1}). This proves (\ref{notsocomplicated}).

\item
Finally, there exists a constant $C>0$ such that 
$$\Vert \frac{1}{\sqrt{T}}\sum_{t=p+1}^T E_t\cdot\frac{O_t^{*}M_t'}{\sigma_t^4}\cdot\left(\hat{q}^{*}_{2,t}-q_{2,t}^{*}\right)\Vert\leq 
\frac{1}{\sqrt{T}}\sum_{t=p+1}^T\frac{\Vert O_t^{*}M_t'\Vert}{\sigma_t^4}\cdot\max_{p+1\leq t\leq T}\vert E_t\vert\times \max_{p+1\leq t\leq T}\Vert \hat{q}^{*}_{2,t}-q_{2,t}^{*}\Vert.$$
Then using Lemma \ref{lemmeutile}, we deduce that 
$$\frac{1}{\sqrt{T}}\sum_{t=p+1}^T E_t\cdot\frac{O_t^{*}M_t'}{\sigma_t^4}\left(\hat{q}^{*}_{2,t}-q_{2,t}^{*}\right)=o_{\P}(1).$$ 
\end{itemize}
\item
Next, we prove the following convergence in distribution: $\frac{1}{\sqrt{T}}\hat{L}_{5,T}\rightarrow \mathcal{N}_n\left(0,\V\left(\xi_1^2\right)\Sigma\right)$.
We have 
\begin{eqnarray*}
\frac{1}{\sqrt{T}}\hat{L}_{5,T}&=&\frac{1}{\sqrt{T}}\sum_{t=p+1}^T\frac{O_t^{*}}{\sigma_t^2}Z_t+\frac{1}{\sqrt{T}}\sum_{t=p+1}^T\frac{E_tO_t^{*}}{\sigma_t^2}Z_t\\
&=& A_T+B_T.
\end{eqnarray*}
The convergence of $(A_T)$ is obtained as in the proof of Theorem $1$. Moreover, using the expression of $E_t$ given in Lemma \ref{lemmeutile} and arguments which are now familiar, the convergence $B_T=o_{\P}(1)$ can be obtained.

\item
Finally we show that $\frac{\hat{D}_T}{T}\rightarrow \Sigma$ a.s. It just remains to prove that 
$$\frac{1}{T}\sum_{t=p+1}^T\frac{E_t}{\sigma_t^4}O^{*}_t\left(O^{*}_t\right)'=o_{\P}(1),\quad \frac{1}{T}\sum_{t=p+1}^T\frac{O^{*}_t\left(O^{*}_t\right)'}{\sigma_t^4}\rightarrow \Sigma.$$ 
Using the fact that $\max_{p+1\leq t\leq T}\Vert E_t\Vert=o_{\P}(1)$ and $\frac{O^{*}_t}{\sigma_t^2}$ is uniformly bounded in $t,T,\omega$, the first assertion follows. The second assertion has been already shown in the proof of Theorem $1$ with general weights $W_t$.$\square$

\end{enumerate}

\section{Proof of Theorem $3$}
The proof uses the arguments given in \citet{Fryz} (see Proposition $3$ of that paper).
We have the decomposition
\begin{equation}\label{decdec}
\hat{\alpha}_t-\alpha(u)=-\hat{q}_{2,b',t}\left(\hat{\beta}-\beta\right)+S_{3,b',t}^{-1}\left(A_t(u)-A^{\#}_t(u)\right)+S_{3,b',t}^{-1}A^{\#}_t(u).
\end{equation}
Note that 
$$A^{\#}_t(u)=\sum_{i=p+1}^Tw_{t,i}(b')W_i(u)M_i(u)\sigma^2_i(u)Z_i$$
and using the central limit theorem for martingale differences, we obtain as in \citet{Fryz}, 
$$\sqrt{Tb'}A^{\#}_t(u)\rightarrow \mathcal{N}_j\left(0,\V(\xi_1^2)\cdot\int K(x)^2dx\cdot \E\left(W_1(u)\sigma_1(u)^4M_1(u)M_1(u)'\right)\right).$$
From Theorem $1$, we have $\hat{\beta}-\beta=O_{\P}\left(\frac{1}{\sqrt{T}}\right)$. Moreover, using Lemma \ref{biais} (point $4$), and the fact that $S_{2,b',t}$ is uniformly bounded in $t,T,\omega$, we get
$$\Vert \hat{q}_{2,b',t}\Vert\leq \max_{p+1\leq t\leq T}\Vert S_{3,b',t}^{-1}\Vert\cdot \max_{p+1\leq t\leq T}\Vert S_{2,b',t}\Vert=O_{\P}(1).$$
This leads to
\begin{equation}\label{presque} 
\hat{q}_{2,b',t}\left(\hat{\beta}-\beta\right)=o_{\P}\left(\frac{1}{\sqrt{Tb'}}\right).
\end{equation}
Moreover we have
\begin{equation}\label{conconv}
\lim_{T\rightarrow \infty}S_{3,b',t}=\E\left(W_1(u)M_1(u)M_1(u)'\right) \mbox{ a.s.}
\end{equation}
To justify (\ref{conconv}), 
it is sufficient to show that for every real-valued process $(Y_t)_{p+1\leq t\leq T}$ of the form $Y_t=f\left(\frac{t}{T},X_{t-1}^2,\ldots,X_{t-p}^2\right)$ with $f$ bounded and Lipschitz, we have  $\sum_{i=p+1}^Tk_{t,i}(b')Y_i\rightarrow \E\left(Y_1(u)\right)$. Using Lemma $1$, the Lipschitz property of $f$ and the properties of the kernel $K$, it is easy to get 
$$\Vert\sum_{i=p+1}^Tk_{t,i}(b')\E\left(W_iM_iM_i'\right)-\E\left(W_1(u)M_1(u)M_1(u)'\right)\Vert=O\left(b'+\frac{1}{T}\right).$$
For the stochastic part, we have, using Proposition $1$, $\sum_{i=p+1}^Tk_{t,i}(b')\bar{Y}_i=O_{\P}\left(\frac{1}{\sqrt{Tb'}}\right)$. 
This justifies the convergence (\ref{conconv}). 

The convergence in distribution announced in Theorem $3$ now easily follows from (\ref{presque}), (\ref{conconv}) and decomposition (\ref{decdec}).$\square$   

\section{Proof of Theorem $4$}
In this proof, we will set 
$$\hat{a}_t=\begin{pmatrix}\hat{\alpha}_t\\\hat{\beta}\end{pmatrix},\quad \mathcal{X}_i=\begin{pmatrix} M_i\\N_i\end{pmatrix},\quad \mathcal{X}_i(u)=\begin{pmatrix} M_i(u)\\N_i(u)\end{pmatrix}.$$
Then, setting $a_t=\begin{pmatrix}\alpha_t\\\beta\end{pmatrix}$, we have the expressions 
$$\hat{\sigma}^2_{t,i}=\mathcal{X}'_i\hat{a}_t,\quad \hat{\sigma}^2_{t,i}(u)=\mathcal{X}_i(u)'\hat{a}_t.$$
As in the proof of Theorem $3$, we have the decomposition 
$$\hat{\alpha}_{*,t}-\alpha(u)=\left(\check{s}_{3,b',t}\right)^{-1}\check{s}_{2,b',t}\left(\beta-\hat{\beta}\right)+\left(\check{s}_{3,b',t}\right)^{-1}\left(\sum_{i=p+1}^Tk_{t,i}(b')\check{W}^{*}_i(u)\sigma_i(u)^2M_i(u)Z_i+B_t(u)-B^{\#}_t(u)\right).$$
\begin{itemize}
\item
We first show that $\check{s}_{3,b',t}\rightarrow \E\left(\frac{M_1(u)M_1(u)'}{\sigma_1(u)^4}\right)$ in probability.
We set $t_{-}=\max(1,t-Tb')$ and $t^{+}=\min(T,t+Tb')$. 
Define the following event 
$$A_{t,T}=\cap_{i=t_{-}}^{t^{+}}\left\{\hat{\sigma}^2_{t,i}>\frac{1}{2}\sigma_i^2\right\}.$$
Note that we have,
$$\frac{\hat{\sigma}_{t,i}^2}{\sigma_i^2}=1+\frac{\mathcal{X}_i'}{\sigma_i^2}\left(\hat{a}_t-a_i\right)$$
and $\frac{\Vert \mathcal{X}_i\Vert}{\sigma^2_i}$ is a random variable bounded by a constant $C>0$.
Then we have the inclusion 
$$\left\{\hat{\sigma}^2_{t,i}\leq \frac{1}{2}\sigma_i^2\right\}\subset \left\{\max_{t_{-}\leq i\leq t^{+}}\Vert \hat{a}_t-a_i\Vert\geq \frac{1}{2C}\right\}.$$
But
$$\max_{t_{-}\leq i\leq t^{+}}\Vert \hat{a}_t-a_i\Vert\leq \Vert \hat{a}_t-a(u)\Vert +\max_{t_{-}\leq i\leq t^{+}}\Vert a(u)-a_i\Vert\leq C_1\left(\Vert \hat{a}_t-a(u)\Vert+b'+\frac{1}{T}\right)$$
for a suitable constant $C_1>0$ and from Theorem $1$ and Theorem $3$, $\Vert\hat{a}_t-a(u)\Vert=O_{\P}\left(b'+\frac{1}{\sqrt{Tb'}}\right)$.
This yields to $\max_{t_{-}\leq i\leq t^{+}}\Vert\hat{a}_t-a_i\Vert=O_{\P}\left(b'+\frac{1}{\sqrt{Tb'}}\right)$. 
 Then we conclude that $\P\left(A_{t,T}^c\right)\rightarrow_{T\rightarrow \infty}0$. 
On the event $A_{t,T}$, we have 
\begin{eqnarray*}
\vert\check{W}^{*}_{t,i}-\frac{1}{\sigma_i^4}\vert&=&\frac{\vert\sigma_i^4-\hat{\sigma}_{t,i}^4-\mu_T\vert}{\sigma_i^4\left(\hat{\sigma}^4_{t,i}+\mu_T\right)}\\
&\leq& 4 \frac{\vert\sigma_i^4-\hat{\sigma}_{t,i}^4-\mu_T\vert}{\sigma_i^8}. 
\end{eqnarray*}
But since $\frac{\hat{\sigma}_{t,i}^2}{\sigma_i^2}-1=\frac{\mathcal{X}_i'}{\sigma_i^2}\left(\hat{a}_t-a_i\right)$, we conclude that 
$$\check{W}^{*}_{t,i}=\frac{1}{\sigma_i^4}\left(1+F_{t,i}\right),$$
with $\max_{t_{-}\leq i\leq t^{+}}\vert F_{t,i}\vert=o_{\P}(1)$. Using the properties of the kernel $K$, this yields to 
$$\check{s}_{3,b',t}=\sum_{i=p+1}^Tk_{t,i}(b')\frac{M_iM_i'}{\sigma_i^4}+o_{\P}(1).$$
Using the arguments used in the proof of Theorem $3$, we have $\sum_{i=p+1}^Tk_{t,i}(b')\frac{M_iM_i'}{\sigma_i^4}\rightarrow \E\left(\frac{M_1(u)M_1(u)'}{\sigma_1(u)^4}\right)$ in probability and the last expectation is also the limit of $\check{s}_{3,b',t}$. 
\item
Next we show the convergence in distribution 
$$\sqrt{Tb'}\sum_{i=p+1}^Tk_{t,i}(b')\check{W}^{*}_{t,i}(u)M_i(u)\sigma_i(u)^2Z_i\rightarrow \mathcal{N}_m\left(0,\mathcal{V}_{*}(u)\right).$$
As shown in the proof of Theorem $3$, this convergence holds if $\check{W}^{*}_{t,i}(u)$ is replaced with $\frac{1}{\sigma_i(u)^4}$ in the last expression and it remains to show that
\begin{equation}\label{oublié} 
\sqrt{Tb'}\sum_{i=p+1}^Tk_{t,i}(b')\left(\check{W}^{*}_{t,i}(u)-\frac{1}{\sigma_i(u)^4}\right)M_i(u)\sigma_i(u)^2Z_i\rightarrow 0
\end{equation}
in probability. 
We will use the following equality 
$$\check{W}^{*}_i(u)-\frac{1}{\sigma_i(u)^4}=\sum_{\ell=1}^{\ell_0}\frac{\left(\sigma_i(u)^4-\hat{\sigma}_{t,i}(u)^4-\mu_T\right)^{\ell}}{\sigma_i(u)^{4\ell+4}}+\frac{\left(\sigma_i(u)^4-\hat{\sigma}_{t,i}(u)^4-\mu_T\right)^{\ell_0+1}}{\sigma_i(u)^{4\ell_0+4}\hat{\sigma}_{t,i}(u)^4}.$$
First, it is easy to show that for $\ell=1,\ldots,\ell_0$, 
\begin{equation}\label{eq+}
\sqrt{Tb'}\sum_{i=p+1}^Tk_{t,i}(b')\frac{\left(\sigma_i(u)^4-\hat{\sigma}_{t,i}(u)^4-\mu_T\right)^{\ell}}{\sigma_i(u)^{4\ell+4}}M_i(u)\sigma_i(u)^2\left(\xi_i^2-1\right)=o_{\P}(1).
\end{equation}
Indeed, using the equality $$\frac{\hat{\sigma}^4_{t,i}(u)}{\sigma_i(u)^4}-1=\left(\hat{a}_t-a_i\right)'\frac{\mathcal{X}_i(u)\mathcal{X}_i(u)'}{\sigma_i(u)^4}\left(\hat{a}_t-a_i\right)+2\frac{\mathcal{X}_i(u)'}{\sigma_i(u)^2}\left(\hat{a}_t-a_i\right),$$
developing $\left(\frac{\hat{\sigma}^4_{t,i}(u)}{\sigma_i(u)^4}-1\right)^{\ell}$ and using the decomposition $\hat{a}_t-a_i=\hat{a}_t-a(u)+a(u)-a_i$ it is clear that $$\sum_{i=p+1}^Tk_{t,i}(b')\frac{\left(\sigma_i(u)^4-\hat{\sigma}_{t,i}(u)^4\right)^{\ell}}{\sigma_i(u)^{4\ell+4}}M_i(u)\sigma_i(u)^2Z_i$$
is composed of terms which write as a product of factors involving some coordinates of the vector $\hat{a}_t-a_t$, $\mu_T$ and of one factor of type 
$\sum_{i=p+1}^Tk_{t,i}(b')O_{u,t,i}\left(\xi_i^2-1\right)$ where $O_{u,t,i}$ is bounded and measurable w.r.t $\sigma\left(\xi_{i-s}: s\geq 1\right)$. Such a term is $o_{\P}\left(\frac{1}{\sqrt{Tb'}}\right)$ because $\hat{a}_t-a(u)=o_{\P}(1)$ and since $\left(O_{u,t,i}Z_i\right)_{p+1\leq i\leq T}$ is a sequence of martingale differences, we have $$\sum_{i=p+1}^Tk_{t,i}(b')O_{u,t,i}Z_i=O_{\P}\left(\frac{1}{\sqrt{Tb'}}\right).$$  
This proves (\ref{eq+}). 

Now we consider the remainder term. Recall that $\max_{t_{-}\leq i\leq t^{+}}\Vert\hat{a}_t-a_i\Vert=O_{\P}\left(b'+\frac{1}{\sqrt{Tb'}}\right)$. 
Using the assumption on $b'$, this entails 
\begin{equation}\label{eq++}
\displaystyle\max_{t_{-}\leq i\leq t^{+}}\left(1-\frac{\hat{\sigma}_{t,i}(u)^4}{\sigma_i(u)^4}-\frac{\mu_T}{\sigma_i(u)^4}\right)^{\ell_0+1}=O_{\P}\left(\left(b'+\frac{1}{\sqrt{Tb'}}\right)^{\ell_0+1}\right)=o_{\P}\left(\frac{1}{\sqrt{Tb'}}\right).
\end{equation}
If we define
$$A_{t,T}(u)=\cap_{i=t_{-}}^{t^{+}}\left\{\hat{\sigma}^2_{t,i}(u)>\frac{1}{2}\sigma_i(u)^2\right\},$$
we have $\P\left(A_{t,T}(u)^c\right)\rightarrow_{T\rightarrow \infty}0$ (the proof is the same as for $A_{t,T}$).
Moreover we have, using (\ref{eq++}),
\begin{eqnarray*}
&&\mathds{1}_{A_{t,T}(u)}\left\Vert\sum_{i=p+1}^Tk_{t,i}(b')M_i(u)\sigma^2_i(u)\frac{\left(\sigma_i(u)^4-\hat{\sigma}_{t,i}(u)^4-\mu_T\right)^{\ell_0+1}}{\sigma_i(u)^{4\ell_0+4}\hat{\sigma}_{t,i}(u)^4}\left(\xi_i^2-1\right)\right\Vert\\
&\leq& 4\max_{t_{-}\leq i\leq t^{+}}\left\{\left\vert1-\frac{\hat{\sigma}_{t,i}(u)^4}{\sigma_i(u)^4}-\frac{\mu_T}{\sigma_i(u)^4}\right\vert^{\ell_0+1}\cdot\frac{\Vert M_i(u)\Vert}{\sigma_i(u)^2}\right\}\cdot\sum_{i=p+1}^Tk_{t,i}(b')\vert \xi_i^2-1\vert\\
&=&o_{\P}\left(\frac{1}{\sqrt{Tb'}}\right).
\end{eqnarray*} 
Then, we conclude that 
$$\sum_{i=p+1}^Tk_{t,i}(b')\frac{\left(\sigma^4_i(u)-\hat{\sigma}^4_{t,i}(u)-\mu_T\right)^{\ell_0+1}}{\sigma_i(u)^{4\ell_0+4}\hat{\sigma}^4_{t,i}(u)}M_i(u)\sigma_i(u)^2Z_i=o_{\P}\left(\frac{1}{\sqrt{Tb'}}\right).$$
This shows (\ref{oublié}).
\item 
Finally, we show that $B_t(u)-B^{\#}_t(u)=O_{\P}(b')$.
We consider the event $\widetilde{A}_{t,T}=A_{t,T}\cap A_{t,T}(u)$. We have of course $\lim_{T\rightarrow\infty}\P\left(\left(\widetilde{A}_{t,T}\right)^c\right)=0$. 
We will just show that 
$$H_t=\check{s}_{1,b',t}-\sum_{i=p+1}^Tk_{t,i}(b')\check{W}^{*}_{t,i}(u)M_i(u)X_i(u)^2=O_{\P}(b'),$$
the control of the difference being the same for the two other quantities $\check{s}_{2,b',t}$ and $\check{s}_{3,b',t}$.
Since $\P\left(\widetilde{A}_{t,T}\right)\rightarrow 0$, we have $H_t\mathds{1}_{\left(\widetilde{A}_{t,T}\right)^c}=O_{\P}(b')$. 
Now, on the event $\widetilde{A}_{t,T}$, we have 
$$\check{W}^{*}_{t,i}M_iX^2_i-\check{W}^{*}_{t,i}(u)M_i(u)X^2_i(u)=\frac{M_i}{\sqrt{\hat{\sigma}_{t,i}^4+\mu_T}}F_i+\frac{X_i(u)^2}{\sqrt{\hat{\sigma}_{t,i}^4(u)+\mu_T}}G_i,$$
where 
$$F_i=\frac{X_i^2}{\sqrt{\hat{\sigma}_{t,i}^4+\mu_T}}-\frac{X_i(u)^2}{\sqrt{\hat{\sigma}_{t,i}(u)^4+\mu_T}},$$
$$G_i=\frac{M_i}{\sqrt{\hat{\sigma}_{t,i}^4+\mu_T}}-\frac{M_i(u)}{\sqrt{\hat{\sigma}_{t,i}(u)^4+\mu_T}}.$$
On $\widetilde{A}_{t,T}$ we have for a suitable constant $C>0$,
\begin{eqnarray*}
&&\vert F_i\vert\\
&\leq& \frac{\vert X_i^2-X_i(u)\vert}{\sigma_i^2}+\xi_i^2\cdot\frac{\vert \hat{\sigma}^2_{t,i}-\hat{\sigma}_{t,i}^2(u)\vert}{\sigma^2_i}\cdot \frac{\hat{\sigma}^2_{t,i}+\hat{\sigma}^2_{t,i}(u)}{\sigma^2_i+\sigma^2_i(u)}\\
&\leq& C \left(\vert X_i^2-X_i(u)^2\vert+\Vert \hat{a}_t\Vert^2\xi_i^2\sum_{\ell=1}^p\vert X^2_{i-\ell}-X^2_{i-\ell}(u)\vert\right). 
\end{eqnarray*}
Moreover, $\mathds{1}_{\widetilde{A}_{t,T}}\frac{M_i}{\sqrt{\hat{\sigma}^4_{t,i}+\mu_T}}$ is bounded uniformly in $t,i,T,\omega$. Then, we obtain for a suitably chosen constant $D>0$,
\begin{eqnarray*}
&&\mathds{1}_{\widetilde{A}_{t,T}}\left\Vert\sum_{i=p+1}^Tk_{t,i}(b')\frac{M_i}{\sqrt{\hat{\sigma}_{t,i}^4+\mu_T}}F_i\right\Vert\\
&\leq& D\left\{\sum_{i=p+1}^Tk_{t,i}(b')\vert X_i^2-X_i(u)^2\vert +\Vert \hat{a}_t\Vert^2\sum_{\ell=1}^p\sum_{i=p+1}^Tk_{t,i}(b')\xi_i^2\vert X_{i-\ell}^2-X^2_{i-\ell}(u)\vert\right\}.
\end{eqnarray*}
Using the fact that $\Vert\hat{a}_t\Vert=O_{\P}(1)$, the support condition on the kernel $K$ and Lemma $1$, we conclude that 
$$\mathds{1}_{\widetilde{A}_{t,T}}\sum_{i=p+1}^Tk_{t,i}(b')\frac{M_i}{\sqrt{\hat{\sigma}_{t,i}^4+\mu_T}}F_i=O_{\P}(b').$$
Similarly, we also get 
$$\mathds{1}_{\widetilde{A}_{t,T}}\sum_{i=p+1}^Tk_{t,i}(b')\frac{X_i(u)^2}{\sqrt{\hat{\sigma}_{t,i}^4(u)+\mu_T}}G_i=O_{\P}(b').$$
This proves that $H_t\mathds{1}_{\widetilde{A}_{t,T}}=O_{\P}(b')$ and then $H_t=O_{\P}(b')$. The proof is now complete.$\square$ 
\end{itemize}

\section{Auxiliary results for the proof of Theorem $5$}
In this section, we assume that the assumptions of Theorem $5$ hold true. The quantities studied in the following lemma are defined in (\ref{decomp}).
\begin{lem}\label{auxtools}
We set $\mathcal{B}=[0,b]\cup [1-b,1]$.
\begin{enumerate}
\item
We have $\max_{u\in[0,1]}\Vert \E(B_u)\Vert=O(b)$ and for all positive integer $q$, $\max_{u\in[0,1]}\E\left(\Vert \bar{B}_u\Vert^{2q}\right)=O\left(\frac{b^q}{T^q}\right)$. 
\item
$\max_{u\in[0,1]}\E\left(\Vert \Delta_u\Vert^4\right)=O\left(\left(Tb\right)^{-2}\right)$. 
\item
We have for all positive integer $q$,
$$\max_{u\in \mathcal{B}}\E\Vert S_u-\kappa_u\Vert^{2q}=O(1),\quad \max_{u\in[b,1-b]}\E\Vert S_u-\kappa_u\Vert^ {2q}=O\left(b^{2q}+(Tb)^{-q}\right).$$
\item
We set $R_u=\widetilde{R}_uS_u^{-1}$. We have $\max_{u\in[0,1]}\Vert S_u^{-1}\Vert=O_{\P}(1)$, $\max_{u\in[b,1-b]}\E\left(\Vert \widetilde{R}_u\Vert^{2q}\right)=O\left(b^{4q}+(Tb)^{-2q}\right)$ and  $\max_{u\in \mathcal{B}}\E\left(\Vert \widetilde{R}_u\Vert^{2q}\right)=O(1)$. 
\end{enumerate}
\end{lem}

\paragraph{Proof of Lemma \ref{auxtools}} 
\begin{enumerate}
\item
The first assertion is obvious. For the second assertion, we use Lemma $2$. It is sufficient to show that for some subscripts $y,w$,
$$\max_{u\in[0,1]}\E\vert \sum_{i=p+1}^Tz_i(u) \overline{W_i\mathcal{X}_{i,y}\mathcal{X}_{i,w}}\vert^{2q}=O\left(\frac{b^q}{T^q}\right),$$
where $z_i(u)=e_i(u)\left(a_{w}(i/T)-a_{w}(u)\right)$ satisfies 
$$\max_{p+1\leq i\leq T,u\in[0,1]}\vert q_i(u)\vert=O\left(\frac{1}{T}\right)\mbox{ and }\max_{u\in[0,1]}\sum_{i=p+1}^Tz_i(u)=O(b).$$ Then the result follows from Lemma $2$ if we write the moment of order $2q$ as a multiple sum.
\item
The result follows from Proposition $3$.
\item
Using Proposition $3$, we have $\max_{u\in[0,1]}\E\Vert \bar{S}_u\Vert^{2q}=O\left(\frac{1}{(Tb)^q}\right)$.
Moreover for the bias part 
$$\E(S_u)-\kappa_u=\sum_{i=p+1}^Te_i(u)\left(\E\left(W_i\mathcal{X}_i\mathcal{X}_i'\right)-\kappa_u\right)+\left(\sum_{i=p+1}^Te_i(u)-1\right)\kappa_u.$$
The first term is bounded by $C\left(b+\frac{1}{T}\right)$ where $C>0$ does not depend on $u,T$ and 
$$\max_{u\in I}\big\vert \sum_{i=p+1}^Te_i(u)-1\big\vert$$ is $O\left(\frac{1}{Tb}\right)$ or $O(1)$ when $I=[b,1-b]$ or $I=[0,1]\setminus[b,1-b]$ respectively.
\item
Since $S_{t/T}=\frac{1}{Tb}\sum_{i=p+1}^T K\left(\frac{t-i}{Tb}\right)W_i\mathcal{X}_i\mathcal{X}_i'$, note that if $\vert u-\frac{t}{T}\vert\leq \frac{1}{T}$, then $\Vert S_u-S_{t/T}\Vert\leq \frac{C}{Tb^2}$ where $C>0$ does not depend on $u,t$ and $T$ (we recall that $W_i\mathcal{X}_i\mathcal{X}_i'$ is bounded). Then it is sufficient to show that $\max_{p+1\leq t\leq T}\Vert S_{t/T}^{-1}\Vert=O_{\P}(1)$. This can be proved as in Lemma \ref{biais}, point $4$. Details are omitted. The other assertions are a consequence of the previous point and of the bound $\Vert \widetilde{R}_u\Vert\leq C\Vert S_u-\kappa_u\Vert^2$.$\square$
\end{enumerate}

\begin{lem}\label{tlcUstat}
Let $\left(\widetilde{\Gamma}(u)\right)_{u\in[0,1]}$ be a family of positive definite matrices such that $u\mapsto \widetilde{\Gamma}(u)$ is Lipschitz. Set $\mathcal{V}_T=\int_0^1\Lambda_u'\widetilde{\Gamma}(u)\Lambda_u du$. Then 
$$T\sqrt{b}\left(\mathcal{V}_T-\E\left(\mathcal{V}_T\right)\right)\rightarrow_{T\rightarrow \infty}\mathcal{N}\left(0,4\Vert K^{*}\Vert_2^2\V^2\left(\xi_1^2\right)\zeta\right),$$
where 
$$\zeta=\int_0^1\mbox{ tr}\left(\left(\widetilde{\Gamma}(u)G(u)\right)^2\right)du,\quad G(u)=\E\left(W_1(u)^2\sigma_1(u)^4\mathcal{X}_1(u)\mathcal{X}_1(u)'\right).$$
Moreover 
$$\E\left(\mathcal{V}_T\right)=\frac{\V\left(\xi_1^2\right)\Vert K\Vert_2^2}{Tb}\cdot\int_0^1\mbox{ tr}\left(\widetilde{\Gamma}(u)G(u)\right)du+o\left(\frac{1}{T\sqrt{b}}\right).$$
\end{lem}

\paragraph{Proof of Lemma \ref{tlcUstat}}
We set for $i,j\in\llbracket p+1,T\rrbracket$, $Q_{i,j}=\int_0^1e_i(u)e_j(u)\widetilde{\Gamma}(u)du$. 
Moreover, let $\mathcal{Y}_i=W_i\sigma_i^2\mathcal{X}_i$.
Then 
$$\mathcal{V}_T=\sum_{i,\ell=p+1}^TZ_{\ell}Z_i\mathcal{Y}'_iQ_{i,\ell}\mathcal{Y}_{\ell}.$$
We use the decomposition $\mathcal{V}_T=2\mathcal{V}_{1,T}+\mathcal{V}_{2,T}$ with
$$\mathcal{V}_{1,T}=\sum_{p+1\leq\ell<i\leq T}Z_{\ell}Z_i\mathcal{Y}_i'Q_{i,\ell}\mathcal{Y}_{\ell},$$
and 
$$\mathcal{V}_{2,T}=\sum_{i=p+1}^TZ_i^2\mathcal{Y}_i'Q_{i,i}\mathcal{Y}_i.$$
Note that $\E\left(\mathcal{V}_T\right)=\E\left(\mathcal{V}_{2,T}\right)$ and $\max_{i,\ell}\Vert Q_{i,\ell}\Vert=O\left(\frac{b}{(Tb)^2}\right)$. 
\begin{itemize}
\item
We first show that 
\begin{equation}\label{zero}
T\sqrt{b}\left(\mathcal{V}_{2,T}-\E\left(\mathcal{V}_{2,T}\right)\right)=o_{\P}(1).
\end{equation}
To show this, we decompose 
$$\mathcal{V}_{2,T}-\E\left(\mathcal{V}_{2,T}\right)=\mathcal{V}_{2,1,T}+\E(Z_1^2)\mathcal{V}_{2,1,T},$$
with 
$$\mathcal{V}_{2,1,T}=\sum_{i=p+1}^T\left(Z_i^2-\E Z_i^2\right)\mathcal{Y}_i'Q_{i,i}\mathcal{Y}_i,$$
$$\mathcal{V}_{2,2,T}=\sum_{i=p+1}^T\overline{\mathcal{Y}_i'Q_{i,i}\mathcal{Y}_i}.$$
Since $\left(\mathcal{Y}_i\right)$ is bounded and $\mathcal{V}_{2,1,T}$ is a sum of uncorrelated random variables, the second moment of $\mathcal{V}_{2,1,T}$ is $O\left(\frac{1}{T^3b^3}\right)$. Then $T\sqrt{b}\mathcal{V}_{2,1,T}=o_{\P}(1)$. 
Moreover, one can decompose $\mathcal{V}_{2,2,T}$ as a finite sum of terms of order $O_{\P}\left(\frac{1}{\sqrt{T^3b^2}}\right)$. This can be easily seen taking the second moment of these terms and then using Lemma $2$. 
This leads to (\ref{zero}).
\item
Next, we prove the assertion on $\E\left(\mathcal{V}_{2,T}\right)$. We have 
\begin{eqnarray*}
\E\left(\mathcal{V}_{2,T}\right)&=&\V\left(\xi_1^2\right)\sum_{i=p+1}^T\E\left(\mathcal{Y}_i'Q_{i,i}\mathcal{Y}_i\right)\\
&=& \V\left(\xi_1^2\right)\sum_{i=p+1}^T\int_0^1\frac{1}{(Tb)^2}K^2\left(\frac{u-\frac{i}{T}}{b}\right)\E\left(\mathcal{Y}_1\left(i/T\right)'\widetilde{\Gamma}(u)\mathcal{Y}_1\left(i/T\right)\right)du+O\left(\frac{b}{(Tb)^2}\right)\\
&=&\V\left(\xi_1^2\right)\sum_{i=p+1}^T\int_0^1\frac{1}{(Tb)^2}K^2\left(\frac{u-\frac{i}{T}}{b}\right)du\E\left(\mathcal{Y}_1\left(i/T\right)'\widetilde{\Gamma}\left(\frac{i}{T}\right)\mathcal{Y}_1\left(i/T\right)\right)+O\left(\frac{1}{T}\right).\\ 
\end{eqnarray*}
Next, note that for a Lipschitz function $f:[0,1]\rightarrow \R$, 
$$\frac{1}{Tb}\sum_{i=p+1}^T\int_0^1K^2\left(\frac{u-\frac{i}{T}}{b}\right)f\left(\frac{i}{T}\right)du-\int_0^1f(u)du\Vert K\Vert_2^2=O\left(b+\frac{1}{Tb}\right).$$
This yields to
$$\E\left(\mathcal{V}_{2,T}\right)-\frac{\V\left(\xi_1^2\right)}{Tb}\int_0^1\E\left(\mathcal{Y}_1(u)'\widetilde{\Gamma}(u)\mathcal{Y}_1(u)\right)du\Vert K\Vert_2^2=O\left(\frac{1}{T}+\frac{1}{(Tb)^2}\right).$$
But using the fact that $Tb^{3/2}\rightarrow \infty$, the last quantity is $O\left(\frac{1}{T\sqrt{b}}\right)$.

\item
Finally, we study the convergence in distribution of $\mathcal{V}_{1,T}$. The argument is to apply the central limit theorem for martingales differences. First, we show the Lindberg condition.
We set $\mathcal{P}_i=\sum_{\ell=p+1}^{i-1}Z_{\ell}\mathcal{Y}_{\ell}'Q_{i,\ell}\mathcal{Y}_i$. Then $\mathcal{V}_{1,T}=\sum_{i=p+1}^T\mathcal{P}_i Z_i$.
It is enough to show that 
\begin{equation}\label{Ling}
\left(T\sqrt{b}\right)^4\sum_{i=p+1}^T\E\left(\mathcal{P}_i^4Z_i^4\right)=o(1).
\end{equation}
Using the independence between $Z_i$ and $\mathcal{P}_i$, the fact that $\mathcal{Y}_i$ is bounded and the Burhhölder inequality, we have for a generic constant $C>0$,
\begin{eqnarray*}
\E\left(\mathcal{P}_i^4Z_i^4\right)&\leq& C\E\left(\Vert \sum_{\ell=p+1}^{i-1}Z_{\ell}\mathcal{Y}_{\ell}'Q_{i,\ell}\Vert^4\right)\\
&\leq& C\left(\sum_{\ell=p+1}^{i-1}\Vert Q_{i,\ell}\Vert^2\right)^2\\
&=& O\left(\frac{1}{\left(T^3b^2\right)^2}\right).
\end{eqnarray*}
This leads to (\ref{Ling}), using the condition $Tb^2\rightarrow\infty$. Now, to obtain the convergence mentioned in Lemma \ref{tlcUstat}, it remains to prove that 
\begin{equation}\label{Ling2}
T^2b\E\left(Z_1^2\right)\sum_{i=p+1}^T\mathcal{P}_i^2-\Vert K^{*}\Vert_2^2\E^2\left(Z_1^2\right)\zeta=o_{\P}(1).
\end{equation}
We will use the following decomposition 
$$\sum_{i=p+1}^T\mathcal{P}_i^2=\sum_{i=p+1}^TA_i+2B_1+B_2,$$
where 
\begin{eqnarray*}
A_i&=&\sum_{\ell,m=p+1}^{i-1}Z_{\ell}\mathcal{Y}_{\ell}'Q_{i,\ell}\overline{\mathcal{Y}_i\mathcal{Y}_i'}Q_{i,m}\mathcal{Y}_{m}Z_m\\
&=& \int_0^1\int_0^1\sum_{\ell,m=p+1}^Te_i(u)e_i(v)e_{\ell}(u)e_m(v)B_{i,\ell,m}(u,v)dudv,
\end{eqnarray*}
with $$B_{i,\ell,m}(u,v)=Z_{\ell}\mathcal{Y}_{\ell}'\widetilde{\Gamma}(u)\overline{\mathcal{Y}_i\mathcal{Y}_i'}\widetilde{\Gamma}(v)\mathcal{Y}_mZ_m$$
and
$$B_1=\sum_{p+1\leq\ell<m\leq T}Z_{\ell}\mathcal{Y}_{\ell}'E_{\ell,m}\mathcal{Y}_mZ_m,$$
$$B_2=\sum_{\ell=p+1}^TZ_{\ell}^2\mathcal{Y}_{\ell}'E_{\ell,\ell}\mathcal{Y}_{\ell},$$
$$E_{\ell,m}=\sum_{i=\ell\vee m+1}^TQ_{i,\ell}\E\left(\mathcal{Y}_i\mathcal{Y}_i'\right)Q_{i,m}.$$
We also set $p_{i,\ell,m}(u,v)=e_i(u)e_i(v)e_{\ell}(u)e_m(v)=O\left(\frac{1}{(Tb)^4}\right)$. 
\begin{itemize}
\item
We first show that $T^2b\sum_{i=p+1}^TA_i=o_{\P}(1)$. We have 
\begin{eqnarray*}
&&\E\left(\vert \sum_{i=p+1}^TA_i\vert^2\right)\\
&\leq& \int_0^1\int_0^1\sum_{i,i'=p+1}^T\sum_{p+1\leq \ell,m<i}\sum_{p+1\leq \ell',m'<i'}
p_{i,\ell,m}(u,v)p_{i',\ell',m'}(u,v)\vert \E\left(B_{i,\ell,m}(u,v)B_{i',\ell',m'}(u,v)\right)\vert dudv\\
&\leq& \int_{\vert u-v\vert \leq 2b}\sum_{i,\ell,m\in\mathcal{M}_u}\sum_{i',\ell',m'\in\mathcal{M}_u}p_{i,\ell,m}(u,v)p_{i',\ell',m'}(u,v)\vert \E\left(B_{i,\ell,m}(u,v)B_{i',\ell',m'}(u,v)\right)\vert dudv,
\end{eqnarray*}
where $\mathcal{M}_u=\left\{s\in\llbracket p+1,T\rrbracket: \vert \frac{s}{T}-u\vert\leq 3b\right\}$. Then the last quantity is 
$$O\left(b\frac{1}{(Tb)^8}(Tb)^3\right)=O\left(\frac{1}{T^5b^4}\right).$$ Indeed, $p_{i,\ell,m}(u,v)$ can be bounded by $\frac{1}{(Tb)^4}$ (up to a constant) and we can apply Lemma $2$ with $z_{T,i}=\mathds{1}_{i\in\mathcal{M}_u}$.  Since $\frac{T^4b^2}{T^5b^4}=\frac{1}{Tb^2}\rightarrow 0$, we get the result.
\item
Next, we show that $T^2bB_1=o_{\P}(1)$. Using the fact that $(\mathcal{Y}_i)_i$ is bounded and the Burkhölder inequality, we get 
$$\E\big\vert \sum_{p+1\leq \ell<m\leq T}Z_{\ell}\mathcal{Y}_{\ell}'E_{\ell,m}\mathcal{Y}_mZ_m\big\vert^2\leq C \sum_{m=p+1}^T\sum_{\ell=p+1}^{m-1}\Vert E_{\ell,m}\Vert^2\leq C\sum_{\vert \ell-m\vert\leq 4Tb}\Vert E_{\ell,m}\Vert^2.$$
Since $\max_{\ell,m}\Vert E_{\ell,m}\Vert=O\left(\frac{1}{T^3b}\right)$, we find that 
$$\E\big\vert \sum_{p+1\leq \ell<m\leq T}Z_{\ell}\mathcal{Y}_{\ell}'E_{\ell,m}\mathcal{Y}_mZ_m\big\vert^2=O\left(\frac{1}{T^4b}\right).$$
This proves $T^2bB_1=o_{\P}(1)$.
\item
Finally, we study the convergence of $T^2b B_2$. First, observing that $\max_{p+1\leq \ell\leq T}\Vert E_{\ell,\ell}\Vert=O\left(\frac{b^2}{(Tb)^3}\right)$ and using Lemma $2$, we get 
$$\sum_{\ell=p+1}^T\left\{Z_{\ell}^2\mathcal{Y}_{\ell}'E_{\ell,\ell}\mathcal{Y}_{\ell}-\E\left(Z_{\ell}^2\mathcal{Y}_{\ell}'E_{\ell,\ell}\mathcal{Y}_{\ell}\right)\right\}=O_{\P}\left(\frac{1}{\sqrt{T^5b^2}}\right).$$
Then it remains to show that
\begin{equation}\label{Ling3}
T^2b\sum_{\ell=p+1}^T\E\left(\mathcal{Y}_{\ell}'E_{\ell,\ell}\mathcal{Y}'_{\ell}\right)-\Vert K^{*}\Vert^2\zeta=o_{\P}(1).
\end{equation}
We first note that 
$$\sum_{\ell=p+1}^T\E\left(\mathcal{Y}_{\ell}'E_{\ell,\ell}\mathcal{Y}_{\ell}\right)=\sum_{p+1\leq \ell<i\leq T}\mbox{ tr }\left(\E\left(\mathcal{Y}'_{\ell}\mathcal{Y}_{\ell}\right)Q_{i,\ell}\E\left(\mathcal{Y}_i\mathcal{Y}_i'\right)Q_{i,\ell}\right).$$
Next, one can verify (we skip the details) that one can replace $\E\left(\mathcal{Y}_i\mathcal{Y}_i'\right)$ with $\E\left(\mathcal{Y}_{\ell}\mathcal{Y}_{\ell}'\right)$ and $Q_{i,\ell}$ with $\int_0^1e_i(u)e_{\ell}(u)du \widetilde{\Gamma}\left(\frac{\ell}{T}\right)$ without changing the limit in (\ref{Ling3}). 
Then it remains to study the limit of 
$$T^2b \sum_{\ell=p+1}^Tf_{\ell}s_{T,\ell},$$
where 
$$f_{\ell}=\mbox{ tr }\left(\big[\E\left(\mathcal{Y}_{\ell}\mathcal{Y}_{\ell}'\right)\widetilde{\Gamma}\left(\frac{\ell}{T}\right)\big]^2\right),$$
and 
$$s_{T,\ell}=\sum_{i=\ell+1}^T\left(\int_0^1e_i(u)e_{\ell}(u)du\right)^2=O\left(\frac{1}{T^3b}\right).$$
Note first that 
$$\sum_{p+1\leq \ell\leq Tb}f_{\ell}s_{T,\ell}+\sum_{\ell=T(1-3b)}^Tf_{\ell}s_{T,\ell}=O\left(\frac{1}{T^2}\right).$$
Moreover, if $Tb\leq \ell\leq T(1-3b)$, we have 
$$s_{T,\ell}=\sum_{i=\ell+1}^{\ell+2Tb}\frac{1}{T^4b^2}H\left(\frac{i-\ell}{Tb}\right)^2,$$
where $H(x)=\int_{-1}^1W(v)W(x+v)dv$. Since 
$$\frac{1}{Tb}\sum_{k=1}^{2Tb}H\left(\frac{k}{Tb}\right)^2-\Vert K^{*}\Vert_2^2=O\left(\frac{1}{Tb}\right)$$
and 
$$\frac{1}{T}\sum_{\ell=Tb}^{T(1-3b)}f_{\ell}\rightarrow \int_0^1\mbox{ tr }\left(\left[\E\left(\mathcal{Y}_1(u)\mathcal{Y}_1(u)'\right)\widetilde{\Gamma}(u)\right]^2\right)du=\zeta,$$
we get (\ref{Ling3}).$\square$

\end{itemize}
\end{itemize}

\section{Proof of Theorem $5$}
The proof of Theorem $5$ uses two lemma, Lemma \ref{auxtools} and Lemma \ref{tlcUstat}, which are stated in the next section.
We set $\widetilde{\kappa}_u=A\kappa_u^{-1}$ and $\mathcal{Y}_i=W_i\sigma_i^2\mathcal{X}_i$. Under the null assumption ($\beta$ is non time-varying), we use the following decomposition 
of the difference $\hat{\beta}-\beta$.
\begin{equation}\label{decomp}
\widetilde{\beta}(u)-\beta=\widetilde{\kappa}_u\Lambda(u)+\widetilde{\kappa}_u\bar{B}_u+z_u+\left[\widetilde{\kappa}_u\left(\kappa_u-S_u\right)\kappa_u^{-1}+R_u\right]\cdot\left[B_u+\Lambda_u\right],
\end{equation}
where
$$\Lambda_u=\sum_{i=p+1}^Te_i(u)\mathcal{Y}_iZ_i,\quad B_u=\sum_{i=p+1}^Te_i(u)W_i\mathcal{X}_i\mathcal{X}_i'\left(a(i/T)-a(u)\right),$$
$$R_u=\widetilde{\kappa}_u\left(\kappa_u-S_u\right)\kappa_u^{-1}\left(\kappa_u-S_u\right)S_u^{-1},\quad z_u=\widetilde{\kappa}_u\sum_{i=p+1}^Te_i(u)\left[\E\left(W_i\mathcal{X}_i\mathcal{X}'_i\right)-\kappa_u\right]\cdot\left[a(i/T)-a(u)\right].$$
Note that $\max_{u\in[0,1]}\Vert z_u\Vert=O\left(\frac{b}{T}+b^2\right)=O\left(b^2\right)$.
The first part of Theorem $5$ will follow from Lemma \ref{tlcUstat}, if we show that
\begin{equation}\label{equivalence}
\int_0^1\left(\widetilde{\beta}(u)-\beta\right)'\Gamma(u)\left(\widetilde{\beta}(u)-\beta\right)du-\int_0^1\Lambda_u'\widetilde{\kappa}(u)'\Gamma(u)\widetilde{\kappa}_u\Lambda_udu=o_{\P}\left(\frac{1}{T\sqrt{b}}\right).
\end{equation}
We set $H(u)=\widetilde{\kappa}_u\Gamma(u)\widetilde{\kappa}_u$. To show (\ref{equivalence}), we use the previous decomposition and proceed as follows. 
\begin{enumerate}
\item
We first show that 
\begin{equation}\label{equivalence1}
T\sqrt{b}\int_0^1\bar{B}_u'H(u)\Lambda_udu=o_{\P}(1).
\end{equation}  
We have for a suitable constant $C>0$, 
$$\max_{u\in[0,1]}\E\big\vert \bar{B}_u'H(u)\Lambda_u\big\vert\leq C\max_{u\in[0,1]}\sqrt{\E\left(\Vert \bar{B}_u\Vert^2\right)\cdot \E\left(\Vert\Lambda_u\Vert^2\right)}.$$
Then using Lemma \ref{auxtools}, we get 
$$\max_{u\in[0,1]}\E\big\vert \bar{B}_u'H(u)\Lambda_u\big\vert=O\left(\frac{1}{T}\right)$$
and (\ref{equivalence1}) easily follows. 
\item
Next we show that 
\begin{equation}\label{equivalence2}
T\sqrt{b}\int_0^1\bar{B}_u'H(u)\bar{B}_udu=o_{\P}(1).
\end{equation} 
Using Lemma \ref{auxtools}, we have $\max_{u\in[0,1]}\E\Vert \bar{B}_u\Vert^2=O\left(\frac{b}{T}\right)$ and 
(\ref{equivalence2}) easily follows. 
\item
Next, we prove that 
\begin{equation}\label{equivalence3}
\int_0^1 z_u'\Gamma(u)\widetilde{\kappa}_u\Lambda_udu=o_{\P}\left(\frac{1}{T\sqrt{b}}\right).
\end{equation}
We have
$$\E\big\vert\int_0^1z_u'\Gamma(u)\widetilde{\kappa}_u\Lambda_udu\big\vert=O\left(\max_{u\in[0,1]}\Vert z_u\Vert\cdot \sqrt{\E\left(\Vert \Lambda_u\Vert^2\right)}\right).$$
Then using Lemma \ref{auxtools}, we get 
$$\E\big\vert\int_0^1z_u'\Gamma(u)\widetilde{\kappa}_u\Lambda_udu\big\vert=O\left(\frac{b^2}{\sqrt{Tb}}\right).$$
This leads to (\ref{equivalence3}), using our bandwidth conditions.
\item
Under our bandwidth conditions, we have $T\sqrt{b}\int_0^1z_u'\Gamma(u)z_udu=O\left(Tb^{4.5}\right)=o(1)$.
\item
Next we show that 
\begin{equation}\label{equivalence4}
\int_0^1\Lambda_u'H(u)\left(\kappa_u-S_u\right)\kappa_u^{-1}\left(B_u+\Lambda_u\right)du=o_{\P}\left(\frac{1}{T\sqrt{b}}\right).
\end{equation}
If $I$ is a subinterval of $[0,1]$, we have for a positive constant $C$,
\begin{eqnarray*}
&&\E\big\vert \int_I\Lambda_u'H(u)\left(\kappa_u-S_u\right)\kappa_u^{-1}\left(\bar{B}_u+\Lambda_u\right)du\big\vert\\
&\leq& C\cdot\vert I\vert\cdot \max_{u\in I}\E^{1/3}\left(\Vert \Lambda_u\Vert^3\right)\cdot\E^{1/3}\left(\Vert \kappa_u-S_u\Vert^3\right)\cdot\left(\E^{1/3}\left(\Vert \bar{B}_u\Vert^3\right)+\E^{1/3}\left(\Vert\Lambda_u\Vert^3\right)\right).
\end{eqnarray*}  
When $I=[0,b]\cup[1-b,1]$, the last bound is $O\left(\frac{1}{T}\right)$, using Lemma \ref{auxtools}.   
Then, 
$$\int_I\Lambda_u'H(u)\left(\kappa_u-S_u\right)\kappa_u^{-1}\left(\bar{B}_u+\Lambda_u\right)du=o_{\P}\left(\frac{1}{T\sqrt{b}}\right).$$
When $I=[b,1-b]$, Lemma \ref{auxtools} leads to 
$$\int_I\Lambda_u'H(u)\left(\kappa_u-S_u\right)\kappa_u^{-1}\left(\bar{B}_u+\Lambda_u\right)du=O\left(\frac{1}{Tb}\cdot\left(b+\frac{1}{\sqrt{Tb}}\right)\right).$$
Using our bandwidth conditions, this yields to
$$\int_0^1\Lambda_u'H(u)\left(\kappa_u-S_u\right)\kappa_u^{-1}\left(\bar{B}_u+\Lambda_u\right)du=o_{\P}\left(\frac{1}{T\sqrt{b}}\right).$$
Moreover, if $I=[0,b]\cup[1-b,1]$ or $I=[b,1-b]$, we have 
\begin{eqnarray*}
&&\E\big\vert \int_I\Lambda_u'H(u)\left(\kappa_u-S_u\right)\kappa_u^{-1}\E(B_u)du\big\vert\\
&\leq& Cb\vert I\vert \max_{u\in I}\sqrt{\E\Vert S_u-\kappa_u\Vert^2}\max_{u\in [0,1]}\sqrt{\E\Vert\Lambda_u\Vert^2}\\
&\leq & C\frac{\sqrt{b}}{\sqrt{T}}\vert I\vert \max_{u\in I}\sqrt{\E\Vert S_u-\kappa_u\Vert^2}.
\end{eqnarray*}
Using Lemma \ref{auxtools} and our bandwidth conditions, we also obtain  $\int_0^1\Lambda_u'H(u)\left(\kappa_u-S_u\right)\kappa_u^{-1}\E(B_u)du=o_{\P}\left(\frac{1}{T\sqrt{b}}\right)$ and then (\ref{equivalence4}). 
\item
Next, setting $M_u=\widetilde{\kappa}_u\left(\kappa_u-S_u\right)\kappa_u^{-1}\left(B_u+\Lambda_u\right)$, we show that 
\begin{equation}\label{equivalence5}
T\sqrt{b}\int_0^1M_u'\Gamma(u)M_udu=o_{\P}(1).
\end{equation}
Using Lemma \ref{auxtools}, we have 
$$\E\Vert M_u\Vert^2\leq C \sqrt{\E\left(\Vert \kappa_u-S_u\Vert^4\right)\cdot \E\left(\Vert B_u+\Lambda_u\Vert^4\right)}\leq C
\sqrt{\E\left(\Vert \kappa_u-S_u\Vert^4\right)}\left(b^2+\frac{1}{Tb}\right).$$
Then, studying $\E\big\vert \int_I M_u'\Gamma(u)M_udu\big\vert$ when $I=[0,b]\cup[1-b,1]$ or $I=[b,1-b]$, (\ref{equivalence5}) follows using  the previous bound, Lemma \ref{auxtools} and our bandwidth conditions. 
\item
Next we prove that 
\begin{equation}\label{equivalence6}
T\sqrt{b}\int_0^1\Lambda_u'\Gamma(u)R_u\left(B_u+\Lambda_u\right)du=o_{\P}(1).
\end{equation}
When $I=[b,1-b]$ or $I=[0,b]\cup[1-b,1]$, we use the bound
\begin{eqnarray*}
&&\big\vert\int_0^1\Lambda_u'\Gamma(u)R_u\left(B_u+\Lambda_u\right)du\big\vert^2\\
&\leq& C\int_I\Vert\Lambda_u\Vert^2\left(\Vert B_u\Vert^2+\Vert\Lambda_u\Vert^2\right)du\cdot \int_I\Vert R_u\Vert^2du.
\end{eqnarray*}
Moreover, using Lemma \ref{auxtools}, we have 
$$\E\int_I\Vert\Lambda_u\Vert^2\left(\Vert B_u\Vert^2+\Vert\Lambda_u\Vert^2\right)du\leq 2\vert I\vert \max_{u\in[0,1]}\sqrt{\E\Vert \Lambda_u\Vert^4}\cdot\sqrt{\E\Vert B_u\Vert^4+\E\Vert\Lambda_u\Vert^4}=\vert I\vert\times O\left(\frac{b}{T}+\frac{1}{(Tb)^2}\right).$$
Now, if $I=[0,b]\cup [1-b,1]$, we have $\int_I\Vert R_u\Vert^2du=O_{\P}(b)$ and we obtain from the previous bounds
$$\big\vert\int_0^1\Lambda_u'\Gamma(u)R_u\left(B_u+\Lambda_u\right)du\big\vert^2=O_{\P}\left(\frac{b^3}{T}+\frac{1}{T^2}\right)=o_{\P}\left(\frac{1}{T^2b}\right).$$
Now, if $I=[b,1-b]$, we have, from Lemma \ref{auxtools}, $\int_I\Vert R_u\Vert^2du=O_{\P}\left(b^4+\frac{1}{(Tb)^2}\right)$ and then
$$\big\vert\int_0^1\Lambda_u'\Gamma(u)R_u\left(B_u+\Lambda_u\right)du\big\vert^2=O_{\P}\left(\left(\frac{b}{T}+\frac{1}{(Tb)^2}\right)\cdot\left(b^4+\frac{1}{(Tb)^2}\right)\right).$$
This is clearly $o_{\P}\left(\frac{1}{T^2b}\right)$ under our bandwidth conditions. Then (\ref{equivalence6}) follows.
\item
Finally, setting $M_u=R_u\left(B_u+\Lambda_u\right)$, we show that 
\begin{equation}\label{equivalence7}
T\sqrt{b}\int_0^1M_u'\Gamma(u)M_udu=o_{\P}(1).
\end{equation}
We have 
$$\int_I M_u'\Gamma(u)M_udu\leq \max_{u\in [0,1]}\Vert S_u^{-1}\Vert^2\cdot \int_I\Vert B_u+\Lambda_u\Vert^2\cdot\Vert\widetilde{R}_u\Vert^2du.$$
Moreover,
\begin{eqnarray*}
&&\E\int_I\Vert B_u+\Lambda_u\Vert^2\cdot\Vert\widetilde{R}_u\Vert^2du\\
&\leq& \vert I\vert\cdot \max_{u\in I}\sqrt{\E\Vert \widetilde{R}_u\Vert^4}\cdot\left(\sqrt{\E\Vert B_u\Vert^4}+\sqrt{\E\Vert\Lambda_u\Vert^4}\right)\\
&\leq& \vert I\vert\cdot \max_{u\in I}\sqrt{\E\Vert \widetilde{R}_u\Vert^4}\cdot O\left(b^2+\frac{1}{Tb}\right). 
\end{eqnarray*}
Considering the two cases for $I$, (\ref{equivalence7}) follows from Lemma \ref{auxtools} and the bandwidth conditions.
\end{enumerate}
The first part of the theorem is now complete. Now we prove that the asymptotic of our statistic remains the same if we replace $\beta$ with the estimate $\hat{\beta}$ of Theorem $1$.
Since $\hat{\beta}-\beta=O_{\P}\left(\frac{1}{\sqrt{T}}\right)$, we have 
$$S_T\left(\widetilde{a},\hat{\beta}\right)-S_T\left(\widetilde{a},\beta\right)=-2 I_T+O_{\P}\left(\frac{1}{T}\right),$$
where $I_T=\int_0^1 \left(\hat{\beta}-\beta\right)'\Gamma(u)\left(\widetilde{\beta}(u)-\beta\right)du$ and it remains to prove 
that
\begin{equation}\label{finini} 
I_T=o_{\P}\left(\frac{1}{T\sqrt{b}}\right).
\end{equation}
 Using decomposition (\ref{decomp}), we have already shown that 
$$\int_0^1\left(\widetilde{\beta}(u)-\beta-\widetilde{\kappa}_u\Lambda_u\right)'\Gamma(u)\left(\widetilde{\beta}(u)-\beta-\widetilde{\kappa}_u\Lambda_u\right)du=o_{\P}\left(\frac{1}{T\sqrt{b}}\right).$$
Using Cauchy-Schwarz inequality, it is easily seen that
$$\int_0^1\left(\hat{\beta}-\beta\right)'\Gamma(u)\left(\widetilde{\beta}(u)-\beta-\widetilde{\kappa_u}\Lambda_u\right)du=o_{\P}\left(\frac{1}{T\sqrt{b}}\right).$$  
Then to show (\ref{finini}), it remains to show $\int_0^1\Gamma(u)\widetilde{\kappa}(u)\Lambda(u)du=o_{\P}\left(\frac{1}{\sqrt{Tb}}\right)$. 
We have
\begin{eqnarray*} 
&&\int_0^1\Gamma(u)\widetilde{\kappa}(u)\Lambda(u)du\\
&=&\sum_{i=p+1}^T\int_0^1\Gamma(u)\widetilde{\kappa}(u)e_i(u)du\cdot\mathcal{Y}_iZ_i\\
&=&O_{\P}\left(\frac{1}{\sqrt{T}}\right).
\end{eqnarray*}
The last equality follows after noticing that $\max_{p+1\leq i\leq T}\Vert \int_0^1\Gamma(u)\widetilde{\kappa}(u)e_i(u)du\Vert=O\left(\frac{1}{T}\right)$ and applying Lemma $2$ componentwise. Then we get (\ref{finini}). The proof of Theorem $5$ is now complete.$\square$

\section{Proof of Proposition $3$}
Setting $\hat{H}_t=X_t^2-\hat{d}_t$ and $H_t=X_t^2-\E X_t^2$ for $p+1\leq t\leq T$, we have
$$\hat{\beta}=\left(\hat{a}_1,\ldots,\hat{a}_p\right)'=\left(\sum_{i=p+1}^T\hat{\mathcal{H}}_t\hat{\mathcal{H}}_t'\right)^{-1}\sum_{t=p+1}^T\hat{H}_t\hat{\mathcal{H}}_t,$$
where $\hat{\mathcal{H}}_t=\left(\hat{H}_{t-1},\ldots,\hat{H}_{t-p}\right)'$. 
We use the decomposition 
$$\hat{H}_t=H_t+\E X_t-\sum_{i=p+1}^Tk_{t,i}\E X_i^2-\sum_{i=p+1}^T k_{t,i}H_i.$$
To prove the result, we will show that 
\begin{equation}\label{fol1}
\frac{1}{\sqrt{T}}\sum_{t=p+1}^TH_t\sum_{i=p+1}^Tk_{t-s,i}H_i=o{\P}(1),\quad s\leq t
\end{equation}
and
\begin{equation}\label{fol2}
\frac{1}{\sqrt{T}}\sum_{t=p+1}^T\sum_{i=p+1}^Tk_{t,i}H_i\cdot\sum_{i=p+1}^Tk_{t-s,i}H_i=o_{\P}(1),\quad s\leq t.
\end{equation}
Using (\ref{fol1}) and (\ref{fol2}) and some arguments given in the proof of Theorem $1$, one can show that $\sqrt{T}\hat{\beta}$ has the same asymptotic distribution than the same estimator but with 
$\hat{H}$ replaced by $H$. Then one can deduce that $\sqrt{T}\hat{\beta}$ converges to a $p-$dimensional Gaussian vector with distribution 
$\mathcal{N}_p\left(0,\sigma^2 I_p\right)$ and the result of Proposition $3$ easily follows. Let us prove (\ref{fol1}) and (\ref{fol2}). 
\begin{enumerate}
\item
For (\ref{fol1}), we decompose the sum into two terms,
$$\sum_{t=p+1}^TH_t\sum_{i=p+1}^Tk_{t-s,i}H_i=\sum_{t=p+1}^Tk_{t-s,t}H_t^2+\sum_{i\neq t}k_{t-s,i}H_iH_t.$$
The expectation of the (positive) first term is smaller to $C/b$ where $C$ is a positive constant. Under the assumptions of Proposition $3$, we have $T b^2\rightarrow 0$. Then the latter expectation is $o\left(\sqrt{T}\right)$. The variance of the second term is 
\begin{eqnarray*}
&&\sum_{i\neq t, i'\neq t'}k_{t-s,i}k_{t'-s,i'}\E\left(H_i H_t H_{i'}H_{t'}\right)\\
&=& \sum_{i\neq t} k_{t-s,i}^2\E\left(H_i^2\right)\E\left(H_t^2\right)+\sum_{i\neq t}k_{t-s,i}k_{i-s,t}\E\left(H_i^2\right)\E\left(H_t^2\right).
\end{eqnarray*}
It is easy to show that this variance is of order $O\left(b^{-1}\right)=o(T)$. This shows (\ref{fol1}). 
\item
Next we show (\ref{fol2}). The $H_i's$ are independent. One can use the results given in \citet{ZW}, Lemma A1 and A3, from which we deduce that 
$$\max_{p+1\leq t\leq T}\left\vert \sum_{i=p+1}^Tk_{t,i}H_i\right\vert=\frac{O_{\P}\left(T^{\frac{1}{2(1+\delta)}}+\sqrt{Tb \log(T)}\right)}{Tb}.$$ 
Using our bandwidth conditions, this entails (\ref{fol2}).
\end{enumerate} 
The proof of Proposition $3$ is now complete.$\square$

\section{Asymptotic semiparametric efficiency}
 
\subsection{Proof of Proposition $2$}
We set $g_t=g(t/T)$ and $e_t=\frac{M_t'g_t+L_t'h}{\sigma_t^2}$. Then using the inequalities
$$0\leq 1-\frac{1}{1+x}-x+x^2\leq x^3,\quad \log(1+x)-x+\frac{x^2}{2}\leq x^3,\quad x\geq 0,$$
we have 
\begin{eqnarray*}
\log\frac{d\P_{T,\alpha+\frac{g}{\sqrt{T}},\beta+\frac{h}{\sqrt{T}}}}{d\P_{T,\alpha,\beta}}\left(X_{p+1},\ldots,X_T\right)
&=&\frac{1}{2}\sum_{t=p+1}^T\left\{\frac{X_t^2}{\sigma_t^2}\left[1-\frac{1}{1+\frac{e_t}{\sqrt{T}}}\right]-\log\left(1+\frac{e_t}{\sqrt{T}}\right)\right\}\\
&=&\Delta_{T,g,h}-\frac{1}{2T}\sum_{t=p+1}^T\left(\xi_t^2-\frac{1}{2}\right)e_t^2+r_T,
\end{eqnarray*}
with $\vert r_T\vert\leq \frac{1}{2T\sqrt{T}}\sum_{t=p+1}^T \left(\xi_t^2+1\right)e_t^3$.
Since $(e_t)_t$ is bounded, we have $r_T=o_{\P_{T,\alpha,\beta}}(1)$.
Next, observe that $e_t$ is Lipschitz in $X_{t-1}^2,\ldots,X_{t-p}^2$. Setting $Y_t=Z_t d^2_t+\frac{e_t^2}{2}$, 
we have $\frac{1}{T}\sum_{t=p+1}^T\left(\xi_t^2-\frac{1}{2}\right)e_t^2=\frac{1}{T}\sum_{t=p+1}^TY_t$
and since $(Y_t)_t$ is a processes of type $II$, we get from Lemma \ref{inter},
$$\E\vert \frac{1}{T}\sum_{t=p+1}^T \bar{Y}_t\vert^2=O\left(\frac{1}{T}\right).$$
This entails $\frac{1}{T}\sum_{t=p+1}^T\bar{Y}_t=o_{\P_{T,\alpha,\beta}}(1)$. 
Finally, noticing that $e_t$ is a bounded and Lipschitz function in $\left(\frac{t}{T},X_{t-1}^2,\ldots,X_{t-p}^2\right)$,
we deduce from Lemma \ref{approxstat},
$$\frac{1}{T}\sum_{t=p+1}^T\E\left(Y_t\right)-\frac{1}{2T}\sum_{t=p+1}^T\E\left(e_1\left(\frac{t}{T}\right)^2\right)=o(1),$$
where $e_t(u)=\frac{M_t(u)'g(u)+N_t(u)'h}{\sigma_t(u)^2}$. 
From Lemma \ref{approxstat}, $u\mapsto \E\left(e_1(u)^2\right)$ is continuous. This leads to 
$$\lim_{T\rightarrow \infty}\frac{1}{T}\sum_{t=p+1}^T\left(\xi_t^2-\frac{1}{2}\right)e_t^2
=\int_0^1 \E\left(e_1(u)^2\right)du
=\Vert(g,h)\Vert^2_{\H},$$
where the limit is in $\P_{T,\alpha,\beta}-$probability. This achieves the proof of Proposition $2$.$\square$

\subsection{Proof of Corollary $1$}
To prove the first assertion, it is enough to check the equality 
\begin{equation}\label{adj}
<\dot{\kappa}^{*}v,(g,h)>_{\H}=h'v,
\end{equation}
for all $(g,h,v)\in H\times \R^n$. 
Using the notations $E(u)=\begin{pmatrix} E_1(u)&E_2(u)\\E_2(u)'&E_3(u)\end{pmatrix}$,
with 
$$E_1(u)=\E\left(\frac{M_1(u)M_1(u)'}{\sigma_1(u)^4}\right),\quad E_2(u)=\E\left(\frac{M_1(u)N_1(u)'}{\sigma_1(u)^4}\right),\quad 
E_3(u)=\E\left(\frac{N_1(u)N_1(u)'}{\sigma_1(u)^4}\right),$$
we have $q^{*}_2(u)=E_1(u)^{-1}E_2(u)$ and it is easy to verify the equality
$$\Sigma=\int_0^1 \left[E_3(u)-E_2(u)'E_1(u)E_2(u)\right]du.$$
Then, it is easy to get (\ref{adj}) using the expression of the scalar product on $H$.
For the second assertion, we apply Theorem $3.11.2$ of \citet{VWW}, using the equality $\Vert \dot{\kappa}^{*}v\Vert_{\H}=2v'\Sigma^{-1}v$. $\square$

\section{Numerical experiments for inference/testing}
 
\subsection{Example of semiparametric estimation.}
We first illustrate the methods of parameters inference in the semiparametric model with constant lag coefficients. We consider the noise distributions $\xi_1\sim\mathcal{N}(0,1)$, $t(9)$ (Student distribution with $9$ degrees of freedom) and $t(5)$.
These three distributions satisfy the moment assumption $\E\xi_0^{4(1+\delta)}<\infty$ used 
in Theorem $1$, Theorem $3$ and Theorem $4$. 
The number of lags is fixed to $p=2$ and the intercept function is defined by $a_0(u)=2+\sin(2\pi u)$. 
We compare the estimates obtained using the procedure described in the paper and the plug-in estimates which are asymptotic optimal. Two sample sizes are considered: $T=500$ and $T=1500$.
Only one bandwidth is used and selected by the CV procedure (the same bandwidth is used for estimating the intercept function and plug-in estimates are also computed using this initial bandwidth).
Note that the $t-$distributions do not satisfy the moment assumption for the asymptotic normality of the plug-in estimator of lag coefficients (in Theorem $2$, we assumed that $\xi$ has moments of any order but our assumption is probably not optimal). The plug-in estimator seems to have a smaller RMSE (see Table \ref{RMSE}), even when $T=500$. Observe also that our estimates are less accurate when the noise has fatter tails.
The RMSE for $\hat{a}_0$ is defined
by $\sqrt{\frac{1}{T}\sum_{t=1}^T\E\left(\hat{a}_0(t/T)-a_0(t/T)\right)^2}$.

\begin{table}[H]
\caption{RMSE for parameter estimation (notation $*$ is for the plug-in estimator) \label{RMSE} }
\begin{center}
\begin{tabular}{|c|c|c|c|c|c|c|c|c|c|}
\cline{2-10}
\multicolumn{1}{c|}{}& \multicolumn{3}{|c|}{$\xi_0\sim\mathcal{N}(0,1)$}& \multicolumn{3}{|c|}{$\xi_0\sim t(9)$}&\multicolumn{3}{|c|}{$\xi_0\sim t(5)$}\\\hline
\multirow{4}{*}{$T=500$}& $\hat{a}_0$& $\hat{a_1}$& $\hat{a}_2$& 
$\hat{a}_0$& $\hat{a_1}$& $\hat{a}_2$&$\hat{a}_0$& $\hat{a_1}$& $\hat{a}_2$\\\cline{2-10}
&$0.5446$&    $0.0859$&    $0.0769$& $0.6380$&    $0.1104$&    $0.1000$&$0.7501$&    $0.1732$&    $0.1430$\\\cline{2-10}
&$\hat{a}_{0,*}$& $\hat{a_{1,*}}$& $\hat{a}_{2,*}$& 
$\hat{a}_{0,*}$& $\hat{a_{1,*}}$& $\hat{a}_{2,*}$&$\hat{a}_{0,*}$& $\hat{a_{1,*}}$& $\hat{a}_{2,*}$\\\cline{2-10}
&$0.5068$  &  $0.0750$&    $0.0651$& $0.5606$&    $0.0949$&    $0.0822$& $0.6489$&    $0.1619$ &   $0.1167$\\\hline
\multirow{4}{*}{$T=1500$}& $\hat{a}_0$& $\hat{a_1}$& $\hat{a}_2$& 
$\hat{a}_0$& $\hat{a_1}$& $\hat{a}_2$&$\hat{a}_0$& $\hat{a_1}$& $\hat{a}_2$\\\cline{2-10}
& $0.3335$ &   $0.0473$  &  $0.0440$ &$0.3844$&    $0.0615$&    $0.0557$& $0.5181$    &$0.1012$&    $0.0963$\\\cline{2-10}
&$\hat{a}_{0,*}$& $\hat{a_{1,*}}$& $\hat{a}_{2,*}$& 
$\hat{a}_{0,*}$& $\hat{a_{1,*}}$& $\hat{a}_{2,*}$&$\hat{a}_{0,*}$& $\hat{a_{1,*}}$& $\hat{a}_{2,*}$\\\cline{2-10}
&  $0.3192$ &   $0.0433$&    $0.0385$   &$0.3571$ &   $0.0536$ &   $0.0471$&  $0.4365$  &  $0.0775$&    $0.0727$ \\\hline
\end{tabular}
\end{center}
\end{table}

\begin{figure}[H]
\centering
\includegraphics[width=8cm,height=7cm]{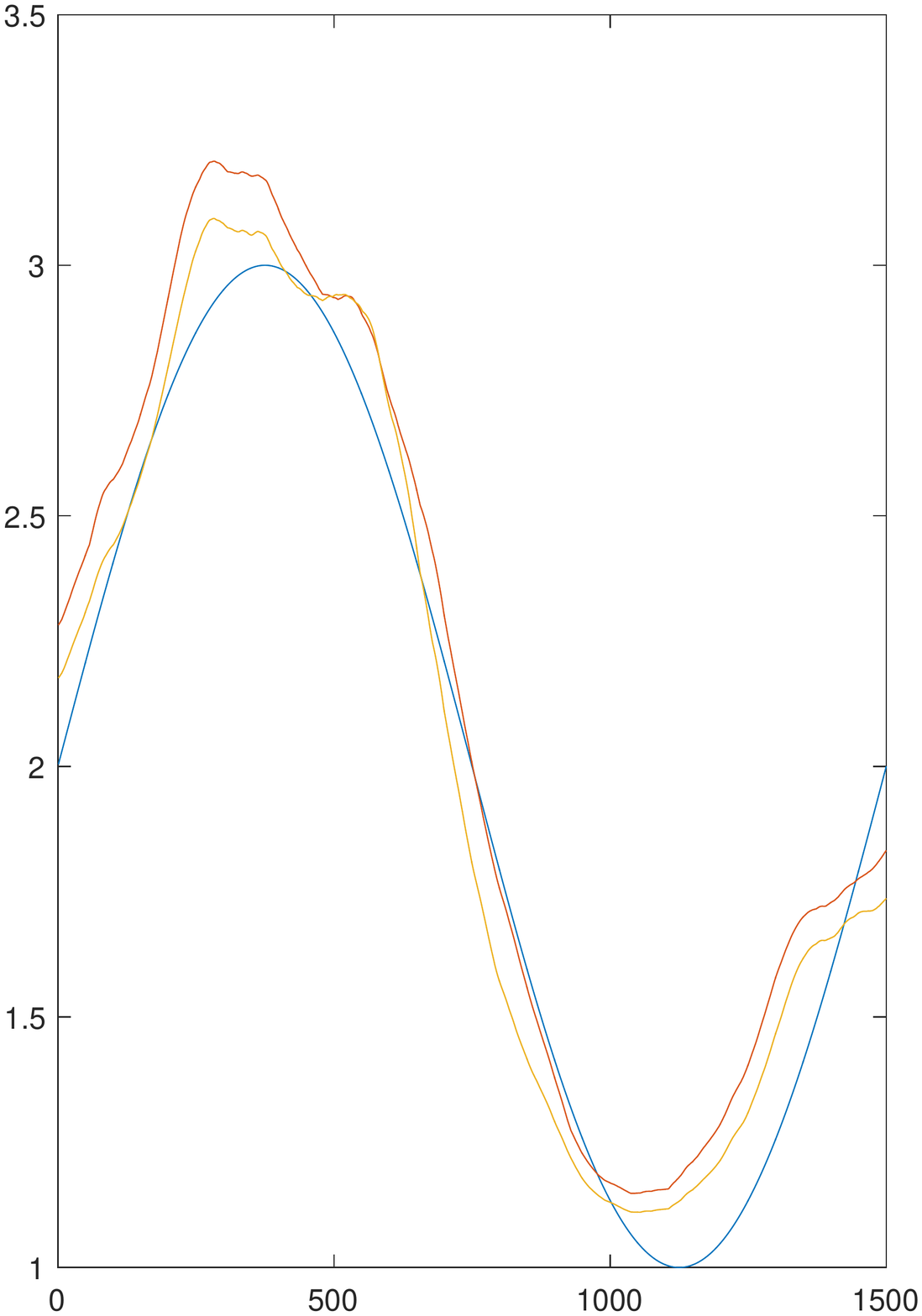}
\caption{Estimation of $a_0(u)=2+sin(2\pi u)$ when $\xi_0\sim\mathcal{N}(0,1)$, $a_1=0.3$, $a_2=0.2$ and $T=1500$ (the red curve is the initial estimate and the yellow curve is the plug-in estimate) \label{estim}}
\end{figure}

\subsection{Testing the constancy of coefficients in a tv(1) process.} 
We consider here a tv-ARCH model with $p=1$ and  $\xi_0\sim\mathcal{N}(0,1)$ or $\xi_0\sim t(9)$.
We consider two setups. In Setup $1$, we have $a_0(u)=2+sin(2\pi u)$ and $a_1(u)=0.5$. In Setup $2$, we have $a_0(u)=1$ and $a_1(u)=0.5+0.25\times \cos\left(2\pi u\right)$. Considering two levels $\alpha=10\%$ and $\alpha=5\%$, we approximate the probability of rejecting $H_0$: $a_0$ constant or $H_0$: $a_1$ constant. Results are reported in Table \ref{paramconst}.
Under $H_0$, this probability has to be close to the level $\alpha$ of the test. One can observe that using a $t(9)-$ distribution for the noise does not create size distortion. However, under the alternative $H_1$, the $t$ distribution entails a smaller power than for the standard Gaussian. This suggests that the power of our tests is impacted by a fat tail noise, which is not surprising. 
Reasonable powers are obtained when $T=2000$, the order of the sample size used in our real data applications.   

\begin{table}[H]
\caption{Approximation of the power for testing parameter constancy in tv$(1)$ processes\label{paramconst} }
\begin{center}
\begin{tabular}{|c|c|c|c|c|c|c|c|c|c|}
\cline{3-10}
\multicolumn{2}{c|}{}& \multicolumn{4}{|c|}{Setup $1$}& \multicolumn{4}{|c|}{Setup $2$}\\
\cline{3-10}
\multicolumn{2}{c|}{}& \multicolumn{2}{|c|}{$T=1000$}& \multicolumn{2}{|c|}{$T=2000$}&\multicolumn{2}{|c|}{$T=1000$}& \multicolumn{2}{|c|}{$T=2000$}\\
\cline{3-10}
\multicolumn{2}{c|}{}&$a_0$&$a_1$&$a_0$&$a_1$&$a_0$&$a_1$&$a_0$&$a_1$\\\hline
\multirow{2}{*}{$\xi_0\sim\mathcal{N}(0,1)$}&$\alpha=5\%$&$0.99$&$0.07$&$1$&$0.07$&$0.08$&$0.54$&$0.07$&$0.91$\\
\cline{2-10}
                                            &$\alpha=10\%$&$0.99$&$0.12$&$1$&$0.13$&$0.13$&$0.68$&$0.12$&$0.96$\\\hline
                                            
\multirow{2}{*}{$\xi_0\sim t(9)$}&$\alpha=5\%$&$0.97$&$0.07$&$1$&$0.06$&$0.06$&$0.34$&$0.06$&$0.69$\\
\cline{2-10}
                                 &$\alpha=10\%$&$0.98$&$0.13$&$1$&$0.12$&$0.13$&$0.47$&$0.11$&$0.8$\\\hline
\end{tabular}
\end{center}
\end{table}

\begin{figure}[H]
\includegraphics[width=5cm,height=5cm]{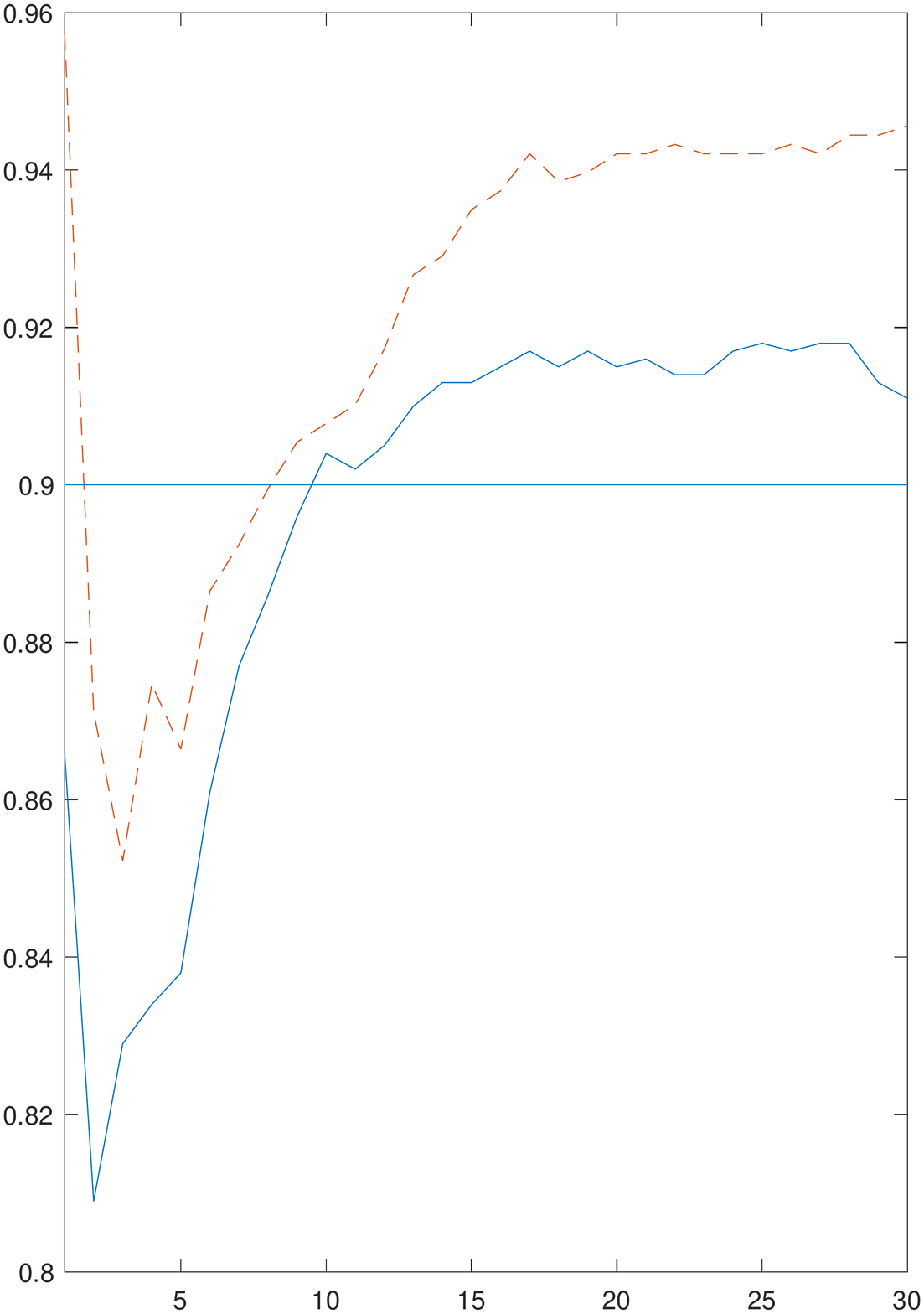}
\includegraphics[width=5cm,height=5cm]{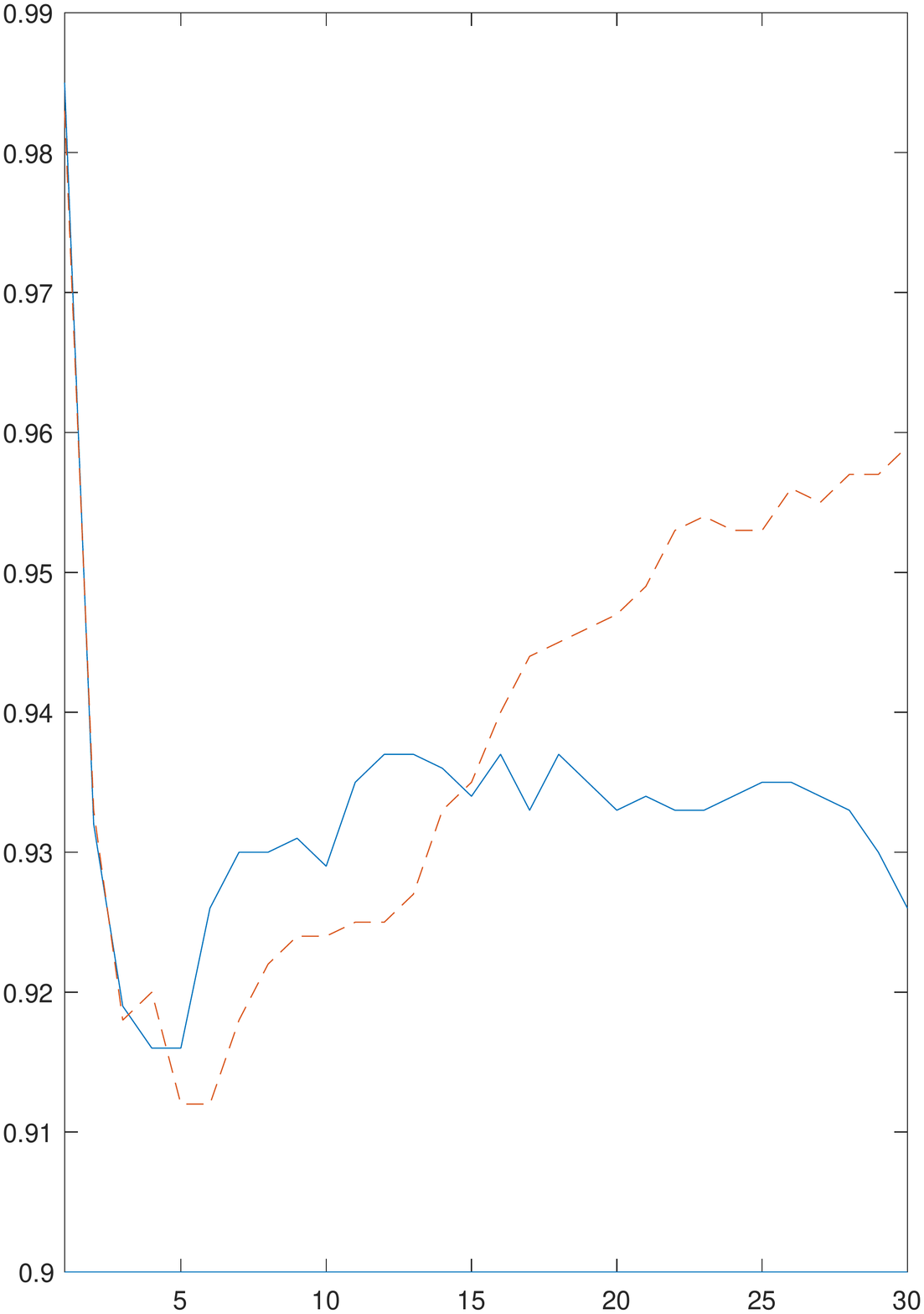}
\includegraphics[width=5cm,height=5cm]{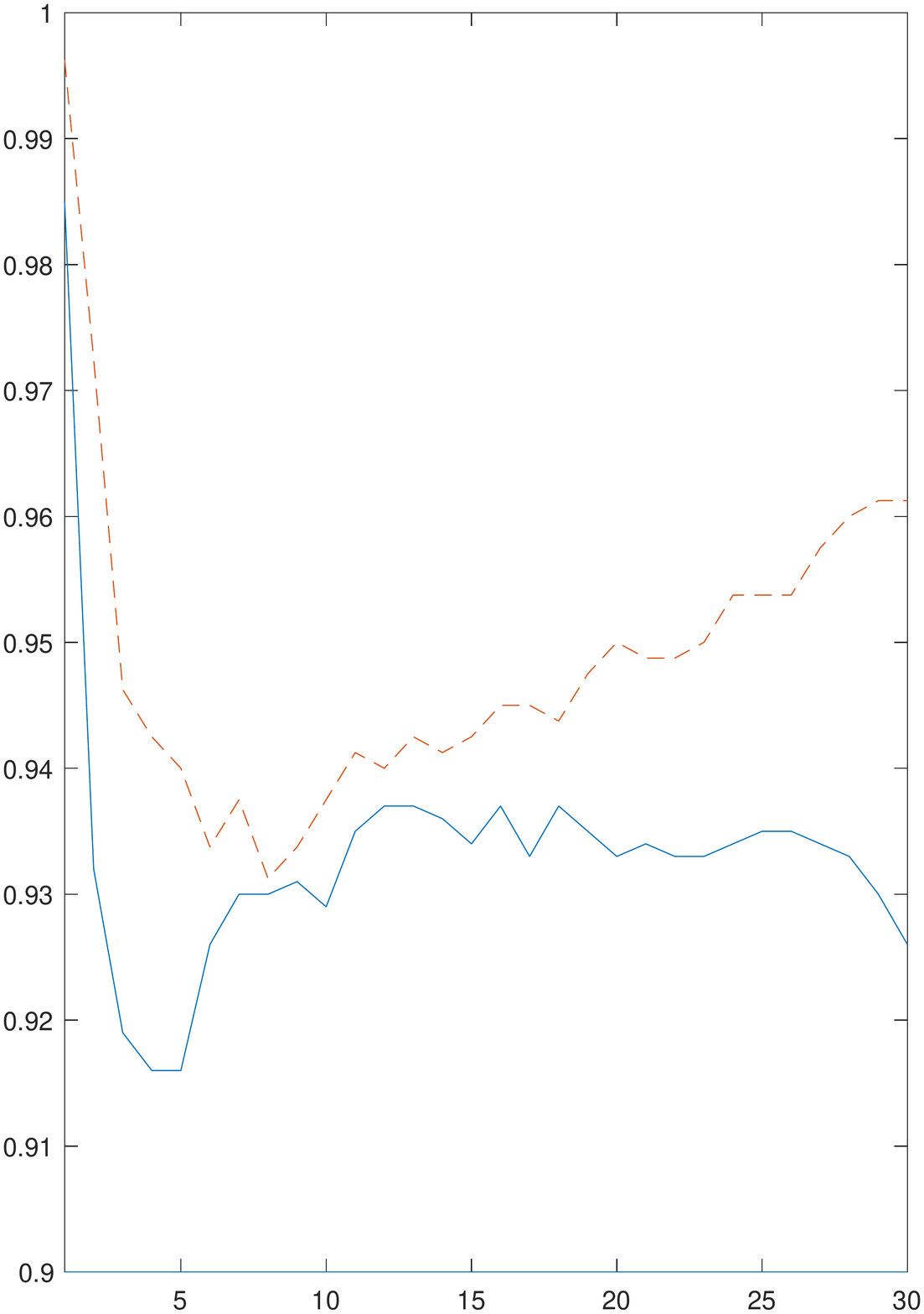}
\caption{Coverage probabilities for Gaussian inputs (left), $t(9)$ inputs (middle) and the difference of two independent random variables following an exponential distribution with parameter $1$ (right). Dashed lines represent the coverage probabilities obtained with the Gaussian quantiles.\label{compa}}
\end{figure}

\paragraph{A comparison with the Gaussian quantiles.} For $T=500$ and $p=1$, we consider the setup $1$. When $\alpha=10\%$ and 
$b=0.01\times \ell$, $1\leq \ell\leq 30$, we compare the coverage probabilities obtained using the Monte Carlo method with the coverage probabilities using the Gaussian quantiles when $\alpha=10\%$ and for testing the constancy of the first lag coefficient. 
In Figure \ref{compa}, one can see that the Monte Carlo method is interesting because the coverage probabilities seem more precise and less sensitive to the bandwidth parameter if we exclude very small bandwidths.  

\subsection{Power curves for testing non time-varying coefficients in a tv(2) process}
\begin{figure}[H]
\begin{center}
\includegraphics[width=8cm,height=7cm]{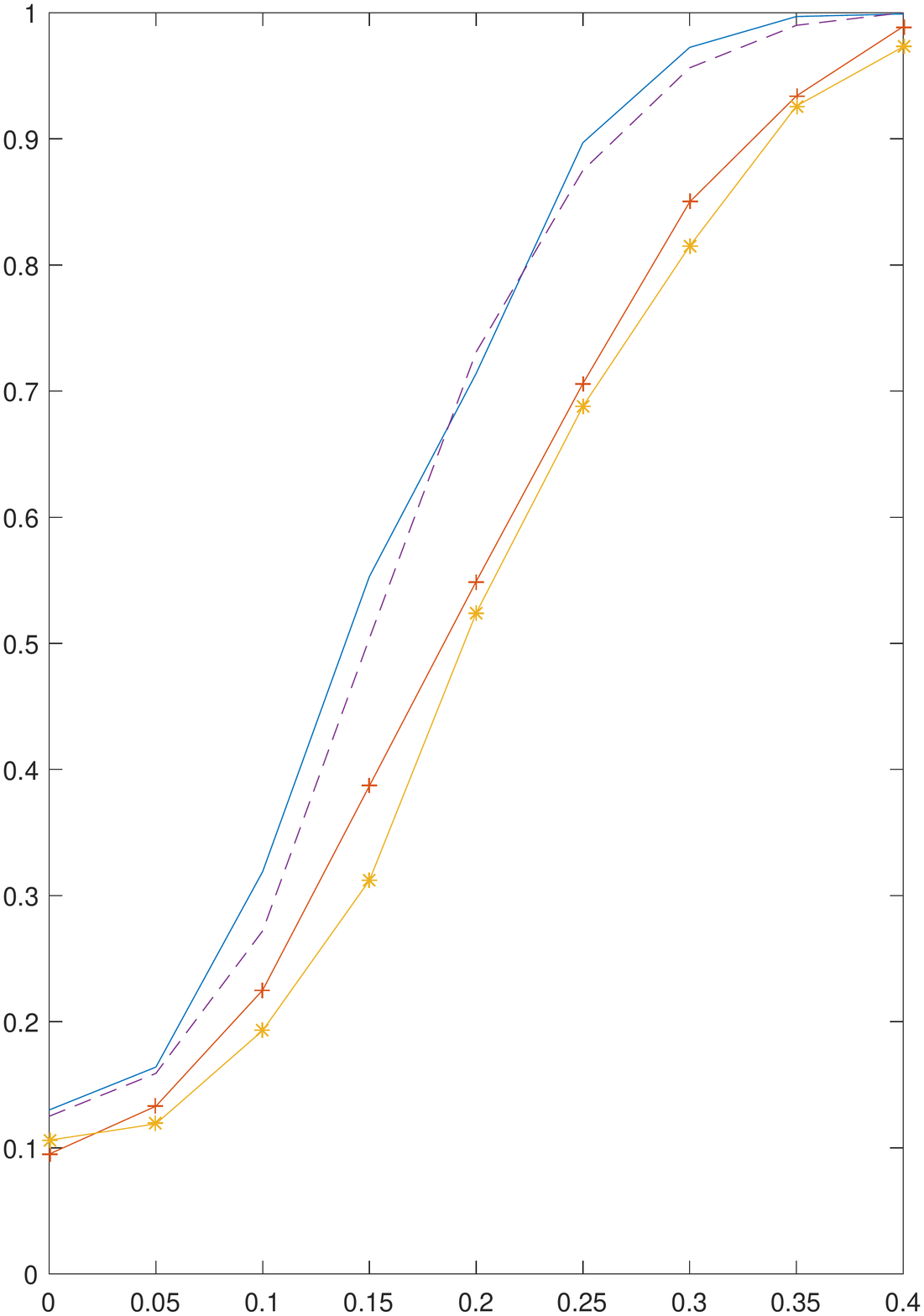}
\includegraphics[width=8cm,height=7cm]{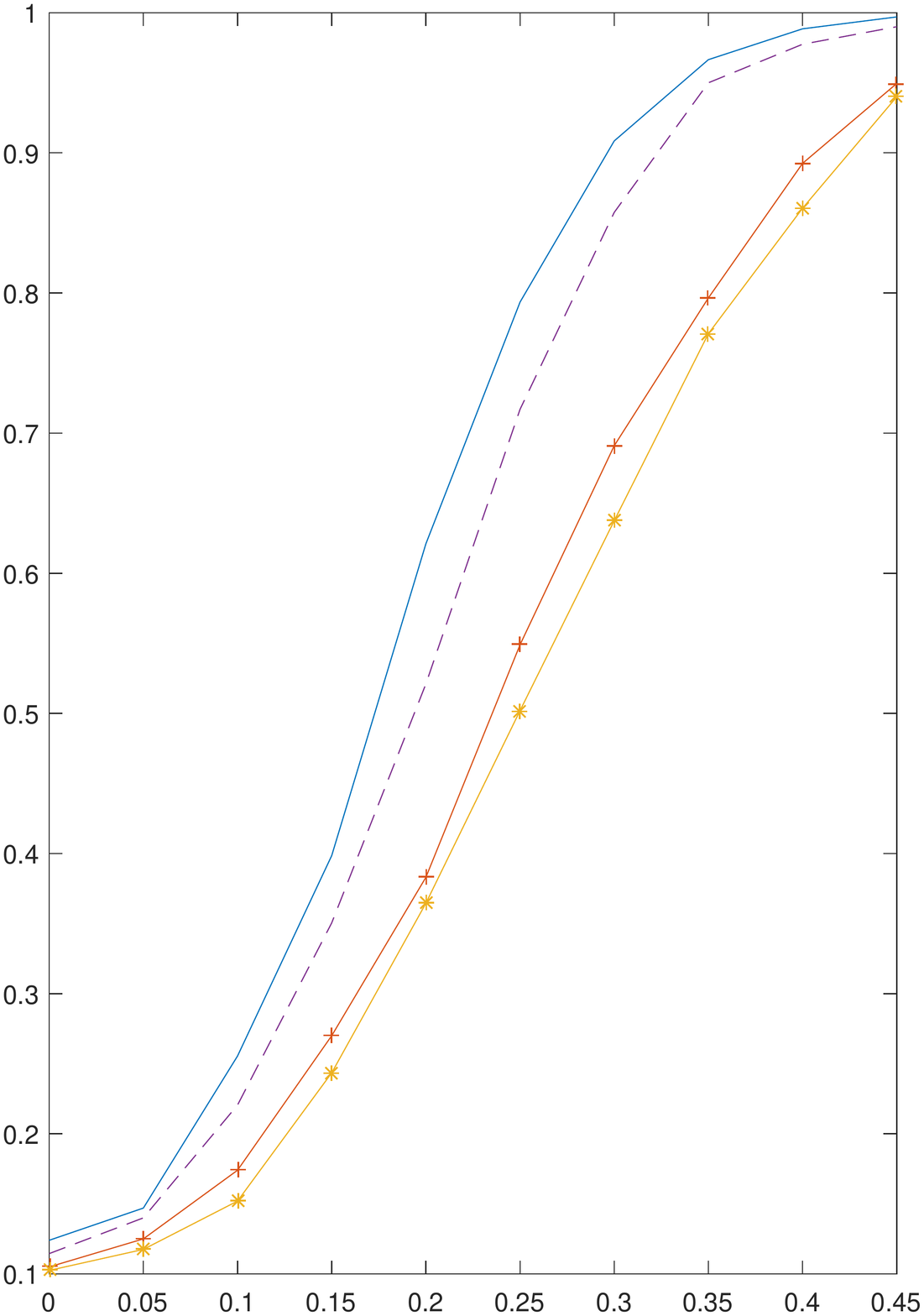}
\end{center}
\caption{Power curves when the noise is Gaussian (on the left) or follows a $t(9)-$distribution (on the right).
The legend for the curves is: $-$ for $a_0$ constant, $--$ for $(a_1,a_2)$ constant, $+$ for $a_1$ constant and $*$ for $a_2$ constant.
\label{powercurve}}
\end{figure}

In this subsection, we simulate approximation of the power for testing $H_0:a_0\mbox{ constant}$ (resp. $a_1$ constant, $a_2$ constant, $(a_1,a2)$ constant) when $a_0(u)=2\left(1+\theta \sin(2\pi u)\right)$, $a_1(u)=0.2+\frac{\theta}{2}\sin(2\pi u)$, $a_2(u)=0.2+\frac{\theta}{2}\cos(2\pi u)$ with $0\leq \theta\leq 0.45$. The noise distribution will be either Gaussian or a student distribution with $9$ degrees of freedom (we remind that Theorem $5$ is only valid when $\E\xi_0^{8(1+\delta)}<\infty$). 
Figure \ref{powercurve} represents an approximation of the power curves when $T=2500$ and $\alpha=10\%$. One can observe that a more fat tail for the noise leads to a slightly smaller power for our tests.

\subsection{Testing the second order dynamic in the semiparametric model}

\begin{figure}[H]
\includegraphics[width=8cm,height=8cm]{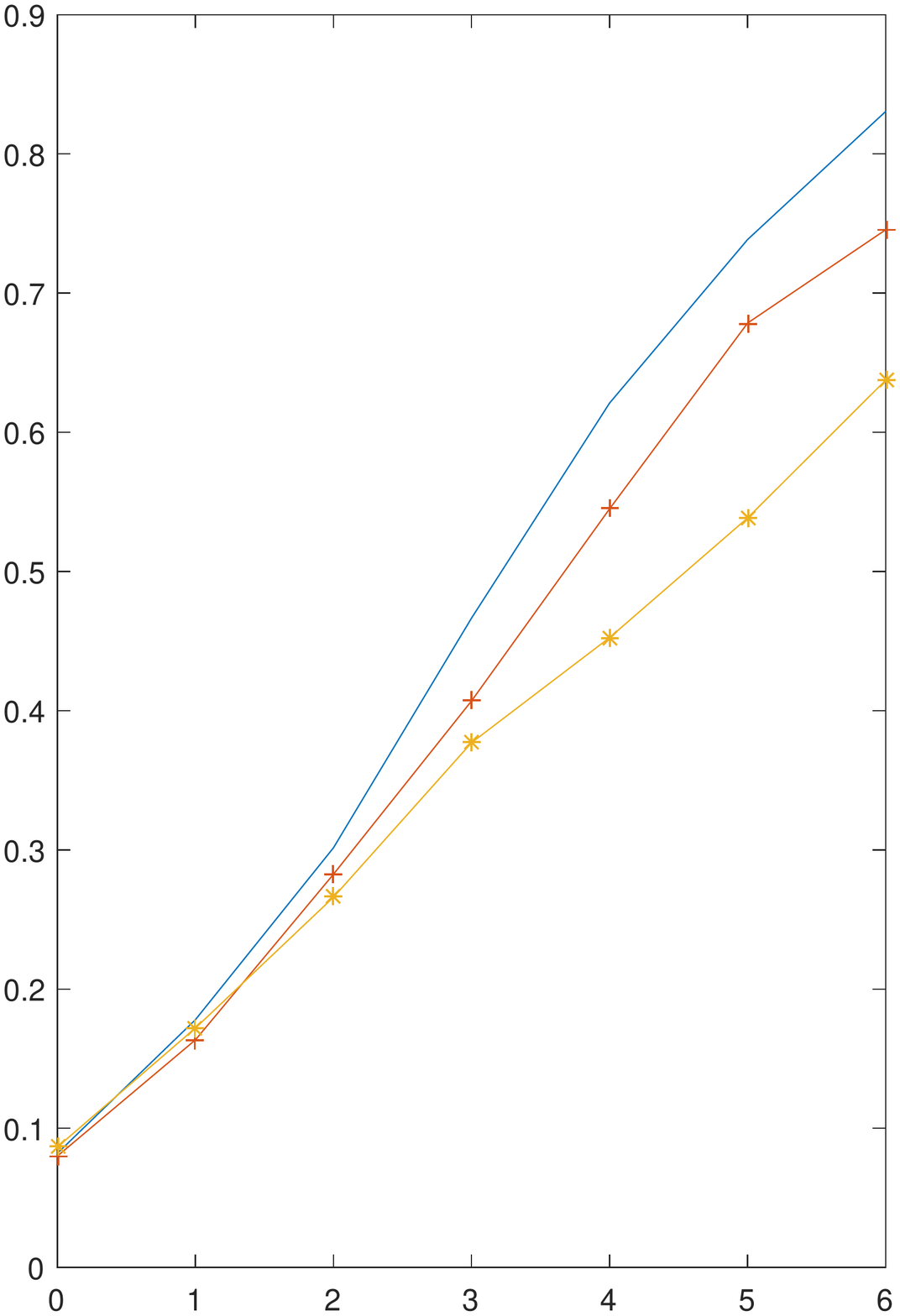}
\includegraphics[width=8cm,height=8cm]{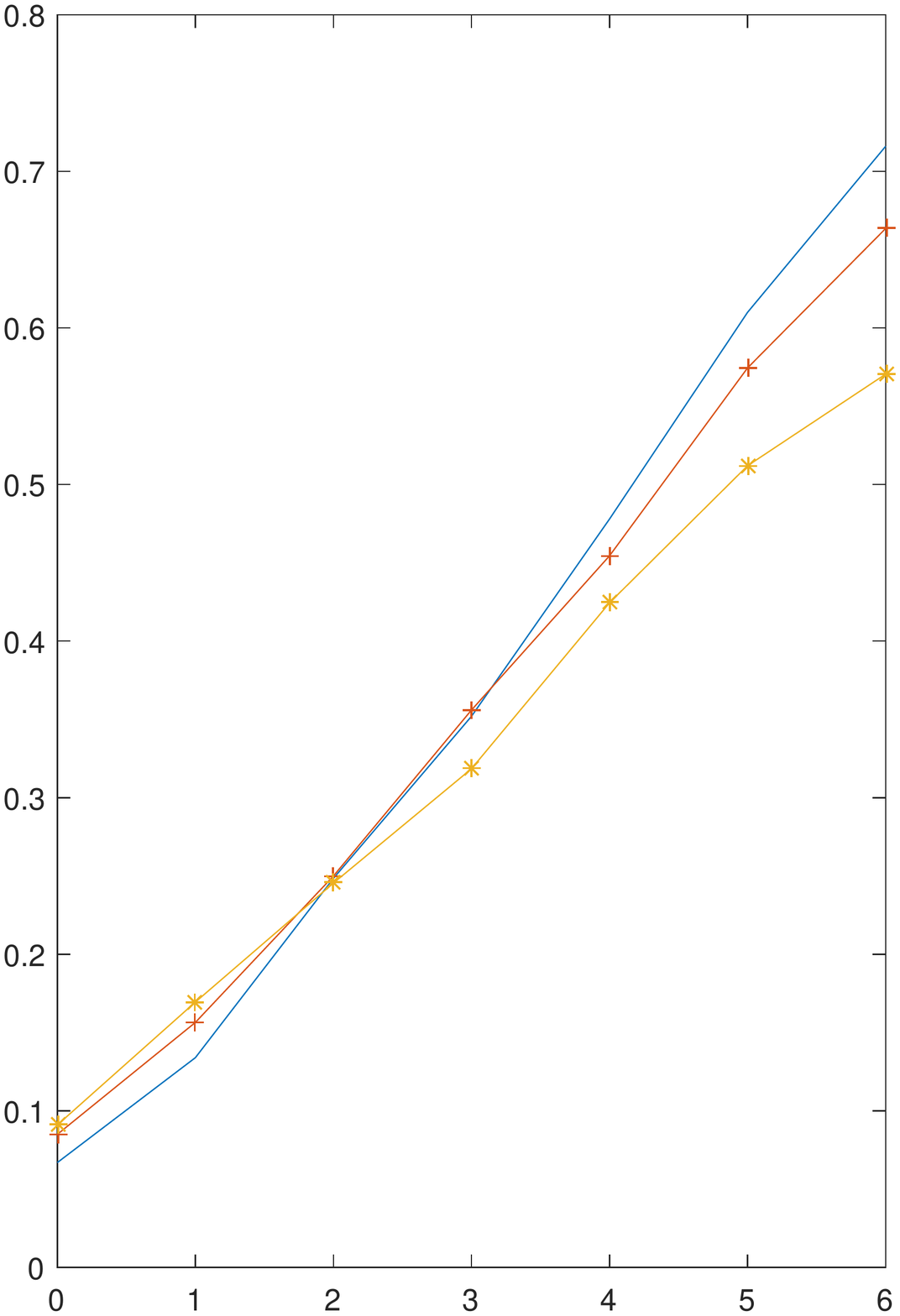}
\includegraphics[width=8cm,height=8cm]{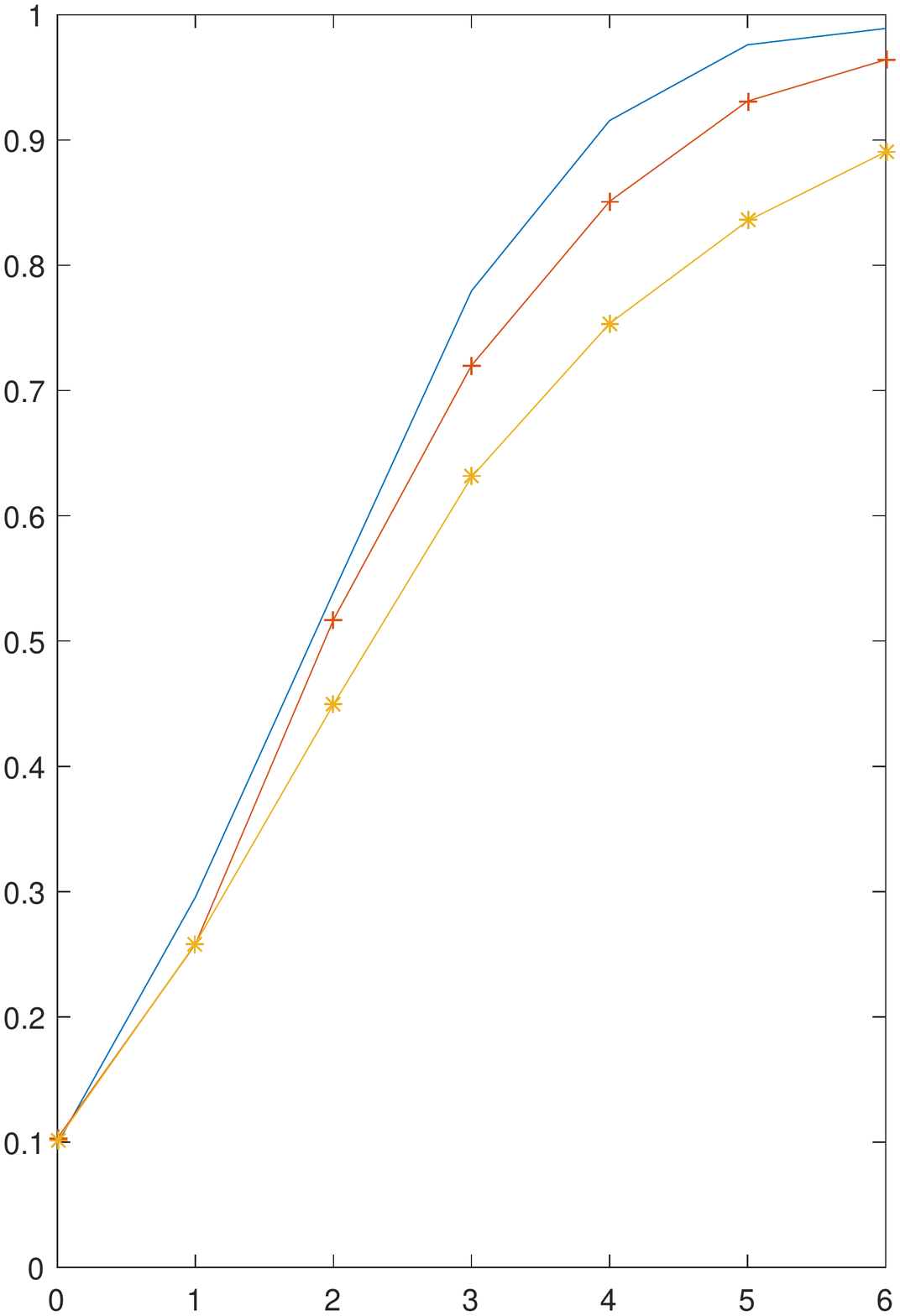}
\includegraphics[width=8cm,height=8cm]{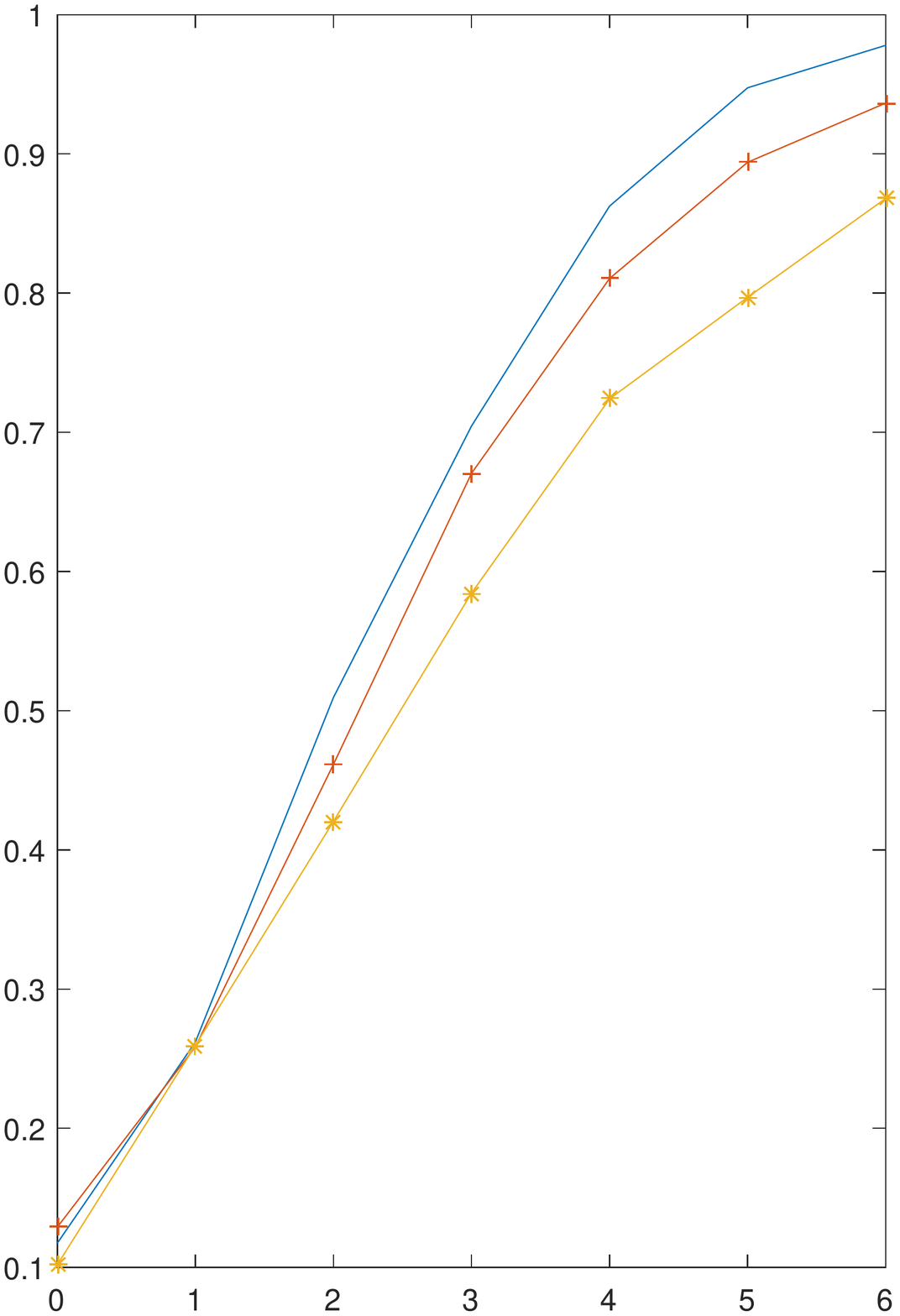}
\caption{Power curves for testing the second order dynamic. Setup $1$ with $T=500$ (top left), setup $2$ with $T=500$ (top right), 
setup $1$ with $T=1000$ (bottom left), setup $2$ with $T=1000$ (bottom right).\label{powercurvesp}} 
\end{figure}

Here we assume that $\beta=\left(a_1,\ldots,a_p\right)'$. The null hypothesis is $\beta=0$. We use the procedure described in the paper after choosing the bandwidth parameter $b$ by cross-validation. We restrict our study to the case $p=2$. 
For the simulation setup, we consider two scenarios. In setup $1$, we consider a constant intercept $a_0(u)=10^{-4}$. In setup $2$, $a_0$ 
is a piecewise affine function such that $a_0(0)=a_0(0.5)=a_0(1)=10^{-4}$ and $a_0(0.25)=a_0(0.75)=4\cdot 10^{-4}$. The noise distribution will be Gaussian, $t(9)$, or $t(5)$. We also consider two sample sizes: $T=500$ and $T=1000$. 
Table \ref{coverdernier1} and Table \ref{coverdernier2} provide approximations of the coverage probabilities. In Figure \ref{powercurvesp},
approximations of some power curves are given under the alternative $a_1=a_1=\theta\times 0.02$, with $\theta=0,\ldots,6$.
The results seem satisfying for the three noise distributions.

\begin{table}[H]
\caption{Approximation of the coverage probabilities when $T=500$\label{coverdernier1}}
\centering
\begin{tabular}{|c|c|c|}
\hline
Setup $1$&  $\alpha=10\%$&   $\alpha=5\%$\\\hline
$\xi_0\sim\mathcal{N}(0,1)$ & $0.92$&    $0.95$\\\hline
$\xi_0\sim t(9)$&  $0.92$&    $0.95$\\\hline
$\xi_0\sim t(5)$&$0.91$& $0.94$\\\hline
\end{tabular}
\hfill
\begin{tabular}{|c|c|c|}
\hline
Setup $2$&  $\alpha=10\%$&   $\alpha=5\%$\\\hline
$\xi_0\sim\mathcal{N}(0,1)$ & $0.93$&    $0.97$\\\hline
$\xi_0\sim t(9)$&  $0.92$&    $0.95$\\\hline
$\xi_0\sim t(5)$&$0.91$& $0.94$\\\hline
\end{tabular}
\end{table}

\begin{table}[H]
\caption{Approximation of the coverage probabilities when $T=1000$\label{coverdernier2}}
\centering
\begin{tabular}{|c|c|c|}
\hline
Setup $1$&  $\alpha=10\%$&   $\alpha=5\%$\\\hline
$\xi_0\sim\mathcal{N}(0,1)$ & $0.90$&    $0.96$\\\hline
$\xi_0\sim t(9)$&  $0.90$&    $0.94$\\\hline
$\xi_0\sim t(5)$&$0.90$& $0.93$\\\hline
\end{tabular}
\hfill
\begin{tabular}{|c|c|c|}
\hline
Setup $2$&  $\alpha=10\%$&   $\alpha=5\%$\\\hline
$\xi_0\sim\mathcal{N}(0,1)$ & $0.88$&    $0.94$\\\hline
$\xi_0\sim t(9)$&  $0.87$&    $0.93$\\\hline
$\xi_0\sim t(5)$&$0.90$& $0.93$\\\hline
\end{tabular}
\end{table}

\subsection{Information criterion for the number of lags in tv-ARCH processes}
In this subsection, we study numerically the performance of the information criterion used for selecting the number of lags in time-varying ARCH processes.
We first consider the case $p=1$, with $a_0(u)=2\left(1+0.4\sin(2\pi u)\right)$ and $a_1(u)=0.3$. Two distributions are considered for the noise: the standard Gaussian and the $t(9)$ distribution. In Table \ref{set1}, we simulate $B=2000$ models for both noise distributions and
three sample sizes: $T=500$, $T=1000$ and $T=2000$. The results are correct for large sample sizes. As for the estimation, one observe 
that the performance of the criterion is sensitive to the tail of the noise distribution.
\begin{table}[H]
\caption{Percentages of models correctly fitted (CF), underfitted (UF) and overfitted (OF) for Setup $1$ \label{set1}}
\centering
\begin{tabular}{|c|c|c|c|c|}
\cline{3-5}
\multicolumn{2}{c|}{}&CF& UF& OF\\\hline
\multirow{3}{*}{$\xi_0\sim \mathcal{N}(0,1)$}& $T=500$&$88$ &$6$ &$6$\\
\cline{2-5}
                                             &$T=1000$&$97$&$2$&$1$\\
\cline{2-5}
                                             &$T=2000$&$99$&$0$&$1$\\\hline 
                                                                                      
\multirow{3}{*}{$\xi_0\sim t(9)$}& $T=500$&$80$ &$14$ &$6$\\
\cline{2-5}
                                   &$T=1000$&$90$&$8$&$2$\\
\cline{2-5}
                                   &$T=2000$&$94$&$5$&$1$\\\hline                                             
\end{tabular}
\end{table}

In a second simulation setup, we consider the case $p=0,1,2$, using the same intercept $a_0$ and setting $a_1(u)=0.2+0.2\cdot\sin\left(2\pi u\right)$,
$a_2(u)=0.2+0.2\cdot \cos\left(2\pi u\right)$. In this case, the lag coefficients can be arbitrary close to zero and the true model more difficult to select. Numerical experiments are reported in Table \ref{set2}.
When $\xi_0\sim t(9)$, large sample sizes are necessary to obtain good results.
Once again, one can explain this behavior by the difficulty of getting accurate estimates with such noise distribution tail when the sample size is not large enough.

\begin{table}[H]
\caption{Percentages of correctly fitted, underfitted and overfitted models for Setup $2$ \label{set2}} 
\begin{tabular}{|c|c|c|c|c|c|c|c|c|c|c|}

\cline{3-11}
\multicolumn{2}{c|}{}&\multicolumn{3}{c|}{$p=0$}& \multicolumn{3}{c|}{$p=1$}& \multicolumn{3}{c|}{$p=2$}\\
\cline{3-11}
\multicolumn{2}{c|}{}&CF& UF& OF&
CF& UF& OF&CF& UF& OF\\\hline
\multirow{3}{*}{$\xi_0\sim \mathcal{N}(0,1)$}& $T=500$&$93$ &$0$ &$7$ & $78$&$16$&$6$&$74$&$22$&$4$\\
\cline{2-11}
                                             &$T=1000$&$93$&$0$&$7$&  $92$& $6$& $2$& $91$&$7$&$2$\\
\cline{2-11}
                                             &$T=2000$&$96$&$0$&$4$&   $99$&$0$&$1$&$99$&$1$&$0$ \\\hline                                              
\multirow{3}{*}{$\xi_0\sim t(9)$}& $T=500$&$91$ &$0$ &$9$ & $66$&$28$&$6$&$58$&$37$&$5$\\
\cline{2-11}
                                              &$T=1000$&$95$&$0$&$5$&  $83$& $14$& $3$& $78$&$20$&$2$\\
\cline{2-11}
                                             &$T=2000$&$96$&$0$&$4$&   $94$&$5$&$1$&$95$&$4$&$1$ \\\hline                                             
\end{tabular}
\end{table}

\end{document}